\documentclass[12pt]{amsart}

\usepackage[margin=1in,marginparwidth=0.7in, marginparsep=0.1in]{geometry}

\pdfoutput=1

\usepackage{mathtools}

\usepackage{etoolbox}

\makeatletter
\let\old@tocline\@tocline
\let\section@tocline\@tocline
\newcommand{\subsection@dotsep}{4.5}
\newcommand{\subsubsection@dotsep}{4.5}
\patchcmd{\@tocline}
  {\hfil}
  {\nobreak
     \leaders\hbox{$\m@th
        \mkern \subsection@dotsep mu\hbox{.}\mkern \subsection@dotsep mu$}\hfill
     \nobreak}{}{}
\let\subsection@tocline\@tocline
\let\@tocline\old@tocline

\patchcmd{\@tocline}
  {\hfil}
  {\nobreak
     \leaders\hbox{$\m@th
        \mkern \subsubsection@dotsep mu\hbox{.}\mkern \subsubsection@dotsep mu$}\hfill
     \nobreak}{}{}
\let\subsubsection@tocline\@tocline
\let\@tocline\old@tocline

\let\old@l@subsection\l@subsection
\let\old@l@subsubsection\l@subsubsection

\def\@tocwriteb#1#2#3{%
  \begingroup
    \@xp\def\csname #2@tocline\endcsname##1##2##3##4##5##6{%
      \ifnum##1>\c@tocdepth
      \else \sbox\z@{##5\let\indentlabel\@tochangmeasure##6}\fi}%
    \csname l@#2\endcsname{#1{\csname#2name\endcsname}{\@secnumber}{}}%
  \endgroup
  \addcontentsline{toc}{#2}%
    {\protect#1{\csname#2name\endcsname}{\@secnumber}{#3}}}%

\newlength{\@tocsectionindent}
\newlength{\@tocsubsectionindent}
\newlength{\@tocsubsubsectionindent}
\newlength{\@tocsectionnumwidth}
\newlength{\@tocsubsectionnumwidth}
\newlength{\@tocsubsubsectionnumwidth}
\newcommand{\settocsectionnumwidth}[1]{\setlength{\@tocsectionnumwidth}{#1}}
\newcommand{\settocsubsectionnumwidth}[1]{\setlength{\@tocsubsectionnumwidth}{#1}}
\newcommand{\settocsubsubsectionnumwidth}[1]{\setlength{\@tocsubsubsectionnumwidth}{#1}}
\newcommand{\settocsectionindent}[1]{\setlength{\@tocsectionindent}{#1}}
\newcommand{\settocsubsectionindent}[1]{\setlength{\@tocsubsectionindent}{#1}}
\newcommand{\settocsubsubsectionindent}[1]{\setlength{\@tocsubsubsectionindent}{#1}}

\renewcommand{\l@section}{\section@tocline{1}{\@tocsectionvskip}{\@tocsectionindent}{}{\@tocsectionformat}}%
\renewcommand{\l@subsection}{\subsection@tocline{2}{\@tocsubsectionvskip}{\@tocsubsectionindent}{}{\@tocsubsectionformat}}%
\renewcommand{\l@subsubsection}{\subsubsection@tocline{3}{\@tocsubsubsectionvskip}{\@tocsubsubsectionindent}{}{\@tocsubsubsectionformat}}%
\newcommand{\@tocsectionformat}{}
\newcommand{\@tocsubsectionformat}{}
\newcommand{\@tocsubsubsectionformat}{}
\expandafter\def\csname toc@1format\endcsname{\@tocsectionformat}
\expandafter\def\csname toc@2format\endcsname{\@tocsubsectionformat}
\expandafter\def\csname toc@3format\endcsname{\@tocsubsubsectionformat}
\newcommand{\settocsectionformat}[1]{\renewcommand{\@tocsectionformat}{#1}}
\newcommand{\settocsubsectionformat}[1]{\renewcommand{\@tocsubsectionformat}{#1}}
\newcommand{\settocsubsubsectionformat}[1]{\renewcommand{\@tocsubsubsectionformat}{#1}}
\newlength{\@tocsectionvskip}
\newcommand{\settocsectionvskip}[1]{\setlength{\@tocsectionvskip}{#1}}
\newlength{\@tocsubsectionvskip}
\newcommand{\settocsubsectionvskip}[1]{\setlength{\@tocsubsectionvskip}{#1}}
\newlength{\@tocsubsubsectionvskip}
\newcommand{\settocsubsubsectionvskip}[1]{\setlength{\@tocsubsubsectionvskip}{#1}}

\patchcmd{\tocsection}{\indentlabel}{\makebox[\@tocsectionnumwidth][l]}{}{}
\patchcmd{\tocsubsection}{\indentlabel}{\makebox[\@tocsubsectionnumwidth][l]}{}{}
\patchcmd{\tocsubsubsection}{\indentlabel}{\makebox[\@tocsubsubsectionnumwidth][l]}{}{}

\newcommand{\@sectypepnumformat}{}
\renewcommand{\contentsline}[1]{%
  \expandafter\let\expandafter\@sectypepnumformat\csname @toc#1pnumformat\endcsname%
  \csname l@#1\endcsname}
\newcommand{\@tocsectionpnumformat}{}
\newcommand{\@tocsubsectionpnumformat}{}
\newcommand{\@tocsubsubsectionpnumformat}{}
\newcommand{\setsectionpnumformat}[1]{\renewcommand{\@tocsectionpnumformat}{#1}}
\newcommand{\setsubsectionpnumformat}[1]{\renewcommand{\@tocsubsectionpnumformat}{#1}}
\newcommand{\setsubsubsectionpnumformat}[1]{\renewcommand{\@tocsubsubsectionpnumformat}{#1}}
\renewcommand{\@tocpagenum}[1]{%
  \hfill {\mdseries\@sectypepnumformat #1}}

\let\oldappendix\appendix
\renewcommand{\appendix}{%
  \leavevmode\oldappendix%
  \addtocontents{toc}{%
    \protect\settowidth{\protect\@tocsectionnumwidth}{\protect\@tocsectionformat\sectionname\space}%
    \protect\addtolength{\protect\@tocsectionnumwidth}{2em}}%
}
\makeatother

\makeatletter
\settocsectionnumwidth{2em}
\settocsubsectionnumwidth{2.5em}
\settocsubsubsectionnumwidth{3em}
\settocsectionindent{1pc}%
\settocsubsectionindent{\dimexpr\@tocsectionindent+\@tocsectionnumwidth}%
\settocsubsubsectionindent{\dimexpr\@tocsubsectionindent+\@tocsubsectionnumwidth}%
\makeatother

\settocsectionvskip{5pt}
\settocsubsectionvskip{0pt}
\settocsubsubsectionvskip{0pt}
    
\settocsectionformat{\bfseries}
\settocsubsectionformat{\mdseries}
\settocsubsubsectionformat{\mdseries}
\setsectionpnumformat{\bfseries}
\setsubsectionpnumformat{\mdseries}
\setsubsubsectionpnumformat{\mdseries}

\let\oldtableofcontents\tableofcontents
\renewcommand{\tableofcontents}{%
  \vspace*{-\linespacing}
  \oldtableofcontents}

\usepackage{soul}

\usepackage{amsfonts,amssymb,latexsym,amsmath,graphicx}
\input xy
\xyoption{all}

\usepackage[backref=page, bookmarks=true, bookmarksopen=true,%
bookmarksdepth=3,bookmarksopenlevel=2,%
colorlinks=true,%
linkcolor=blue,%
citecolor=blue,%
filecolor=blue,%
menucolor=blue,%
urlcolor=blue]{hyperref}

\usepackage{cleveref}

\usepackage[all, knot, poly]{xy}
\xyoption{knot}

\usepackage{times}
\usepackage{mathrsfs}

\usepackage{tikz-cd}
\usepackage{tikz}
\usetikzlibrary{arrows}

\newtheorem{lemma}{Lemma}[section]
\newtheorem{proposition}[lemma]{Proposition}
\newtheorem{corollary}[lemma]{Corollary}
\newtheorem{theorem}[lemma]{Theorem}

\theoremstyle{definition}
\newtheorem{definition}[lemma]{Definition}
\newtheorem{definition-lemma}[lemma]{Definition-Lemma}

\theoremstyle{remark} 
\newtheorem{remark}[lemma]{Remark}

\newtheorem{example}[lemma]{Example}
\newtheorem{example-warning}[lemma]{Example-Warning}

\newcommand{\C}{\mathbb{C}}

\newcommand{\E}{\mathbb{E}}

\renewcommand{\H}{\mathbb{H}}

\newcommand{\N}{\mathbb{N}}

\renewcommand{\P}{\mathbb{P}}

\newcommand{\R}{\mathbb{R}}

\newcommand{\Z}{\mathbb{Z}}

\newcommand{\cA}{\mathcal{A}}
\newcommand{\cB}{\mathcal{B}}
\newcommand{\cC}{\mathcal{C}}
\newcommand{\cD}{\mathcal{D}}
\newcommand{\cE}{\mathcal{E}}
\newcommand{\cF}{\mathcal{F}}

\newcommand{\cH}{\mathcal{H}}

\newcommand{\cK}{\mathcal{K}}
\newcommand{\cL}{\mathcal{L}}
\newcommand{\cM}{\mathcal{M}}
\newcommand{\cN}{\mathcal{N}}
\newcommand{\cO}{\mathcal{O}}
\newcommand{\cP}{\mathcal{P}}
\newcommand{\cQ}{\mathcal{Q}}
\newcommand{\cR}{\mathcal{R}}
\newcommand{\cS}{\mathcal{S}}

\newcommand{\cU}{\mathcal{U}}
\newcommand{\cV}{\mathcal{V}}
\newcommand{\cW}{\mathcal{W}}

\newcommand{\ff}{\mathfrak{f}}
\newcommand{\fg}{\mathfrak{g}}
\newcommand{\fM}{\mathfrak{M}}
\newcommand{\fX}{\mathfrak{X}}
\newcommand{\fY}{\mathfrak{Y}}

\newcommand{\sS}{\mathscr{S}}
\newcommand{\sR}{\mathscr{R}}

\DeclareMathOperator{\Hom}{Hom}

\DeclareMathOperator{\Ad}{Ad}
\DeclareMathOperator{\Loc}{Loc}
\DeclareMathOperator{\Map}{Map}
\DeclareMathOperator{\sh}{sh}
\DeclareMathOperator{\msh}{\mu sh}
\DeclareMathOperator{\Pin}{Pin}

\DeclareMathOperator{\Rp}{\R_+}

\DeclareMathOperator{\Sp}{Sp}
\DeclareMathOperator{\Spc}{Spc}

\newcommand\Ab{ \operatorname{Ab} }

\newcommand\Cat{ \operatorname{Cat} }

\newcommand\colim{ \operatorname{colim} }

\newcommand\Cont{ \operatorname{Cont} }

\newcommand\dT{ {\dot{T} } }
\newcommand\Equiv{ \operatorname{Equiv} }
\newcommand\forcoe[2]{ \prescript{\ }{#1}{#2} } 
\newcommand\fib{ \operatorname{fib} }
\newcommand\Fun{ \operatorname{Fun} }

\newcommand\id{ \operatorname{id} }

\newcommand\Lie{ \operatorname{Lie}}
\newcommand\mcsh{ \operatorname{\mu_\C sh } }
\newcommand\mD{ \operatorname{\mu \mathbb{D}} }
\newcommand\mhom{ \operatorname{\mu hom } }

\newcommand\mRH{ \operatorname{\mu RH } }

\newcommand\Perv{  {  \operatorname{\P erv} } }
\newcommand\pre{ {\operatorname{pre} } }

\newcommand\PrLst{  {  \operatorname{Pr}^{\operatorname{L} }_{st} } }
\newcommand\PrRst{  {  \operatorname{Pr}^{\operatorname{R} }_{st} } }
\newcommand\Pmsh{ \operatorname{\P \msh } }
\newcommand\Ring{ {\operatorname{Ring} } }

\newcommand\rev{ {\operatorname{rev} } }
\newcommand\sHom{ \operatorname{\mathcal{H}om} }
\newcommand\sMod{ \operatorname{\fM od} }
\newcommand\supp{ {\operatorname{supp} } }
\newcommand\Symp{ {\operatorname{Symp} } }
\newcommand\tmhom{ \operatorname{t- \mu hom} }
\newcommand\uboxtimes{ {\operatorname{\underline{\boxtimes} } } }

\renewcommand{\mod}{\operatorname{mod}}

\title[The microlocal Riemann-Hilbert correspondence for complex contact manifolds]{The microlocal Riemann-Hilbert correspondence \\ for complex contact manifolds}

\author{Laurent C\^ot\'e, Christopher Kuo, David Nadler, and Vivek Shende}

\begin{document}

\begin{abstract}
Kashiwara showed in 1996 that the categories of microlocalized D-modules can be canonically glued to give a sheaf of categories over a complex contact manifold.  Much more recently, and by rather different considerations, we constructed a canonical notion of perverse microsheaves on the same class of spaces.  Here we provide a Riemann-Hilbert correspondence. 
\end{abstract}

\maketitle

\renewcommand{\contentsname}{}
\vspace{-15mm}
\renewcommand{\baselinestretch}{0.9}
\small  
\tableofcontents
\normalsize
\renewcommand{\baselinestretch}{1}



\section{Introduction}

The Riemann-Hilbert correspondence, originating in the question of the existence of Fuchsian differential equations whose solutions have prescribed monodromy, takes as its modern form 
the assertions that (1) the derived solution functor carries holonomic $\cD$-modules on an algebraic variety to perverse sheaves \cite{kashiwara-solperverse} and (2) after restricting to  regular holonomic $\cD$-modules, this is an equivalence of categories \cite{kashiwara-RH, mebkhout-RH}. 

From the beginning, this correspondence was understood to have a microlocal character.  In particular, Kashiwara had previously introduced the characteristic variety of a $\cD$-module 
and  generalized the Cauchy-Kowalevskya theorem: the characteristic variety of the $\cD$-module records the obstructions to propagation of its solutions \cite{Kashiwara-masters}.  That is, in the later terminology of Kashiwara and Schapira, the solution functor carries the characteristic variety of a holonomic $\cD$-module to the microsupport of its sheaf of solutions \cite[Theorem 11.3.3]{kashiwara-schapira}.

Both the notion of $\cD$-module and sheaf `microlocalize' over the cotangent bundle.\footnote{As a default, we will use $\cD$-module to mean an  object of the abelian category, but use sheaf to mean an object of the derived category.}  On the one hand, for $\cD$-modules, there is a sheaf $\cE$ of `microlocal differential operators' over the cotangent bundle, whose pushforward to the base recovers $\cD$. One can then introduce the sheaf of $\cE$-modules \cite{sato-kashiwara-kawai, kashiwara-microlocal-D}, whose  pushforward to the base likewise recovers the sheaf of $\cD$-modules.  Moreover, the notions of support coincide: the support of an $\cE$-module is the characteristic variety of the corresponding $\cD$-module \cite[Theorem 7.27]{kashiwara-microlocal-D}.

On the other hand, 
for sheaves, the notion of microsupport gives a formal microlocalization as developed by Kashiwara and Schapira \cite{kashiwara-schapira}.
Namely, given the notion of microsupport $ss(F) \subset T^*X$ for a sheaf $F$ on a manifold $X$, one can form the sheaf of categories $\msh$ over $T^*X$ by sheafifying the presheaf of categories
$$\msh^{pre}(U) := \sh(X)/ \{F\,|\, ss(F) \cap U = \emptyset \}$$
It is not difficult to show that $\msh(T^*X) = \sh(X)$, and more generally,
the pushforward of the sheaf of categories $\msh$ recovers the sheaf of categories $\sh$. Moreover, the notions of support coincide: the support of a microsheaf matches the microsupport of the corresponding sheaf. 

The work of Andronikov \cite{andronikof-microlocal-RH} and Waschkies \cite{waschkies-microlocal-RH} establishes a compatible microlocalization of the Riemann-Hilbert correspondence: the derived solution functor induces an equivalence between
regular holonomic $\cE$-modules with
perverse microsheaves
\cite[Theorem 3.6.5]{waschkies-microlocal-RH}. The new content in this result occurs away from the zero-section, hence on the projectivized cotangent bundle $\P^* X$ viewed as a complex contact manifold.

The main theorem of this paper is a globalization of this microlocal equivalence.

\begin{theorem}[\ref{thm:gloabl mrh}]\label{thm:gloabl mrh intro}
Let $V$ be a complex contact manifold. There is a canonical equivalence
$$\mRH_V: \P erv_V \xrightarrow{\sim} \cE_V - \sMod_{rh}$$ 
between perverse microsheaves on $V$ (in the sense of  \cite{CKNS}) and regular holonomic $\cE_V$-modules (in the sense of  \cite{kashiwara-quantization-contact}) extending the microlocal Riemann--Hilbert correspondence.
\end{theorem}

Here the domain category $\P erv_V$ of perverse microsheaves on $V$ was 
introduced 
in our prior work \cite{CKNS}.
It builds upon the works \cite{shende-microlocal, nadler-shende} which  construct a sheaf of microsheaves $\msh_{V, \xi}$ on any {\em real} contact manifold $V$ equipped with `Maslov data' $\xi$.  
(A couple of remarks: this construction did not proceed via gluing, but rather via high-codimensional embeddings; 
as a major application, it allows for the expression of Fukaya categories of Weinstein manifolds in terms of microlocal sheaf theory \cite{GPS1, GPS2, GPS3}.) 
We subsequently showed  in \cite{CKNS} that on a complex contact manifold $V$, there is a canonical sheaf of perverse microsheaves $\P erv_V$,  locally agreeing with the perverse microsheaves of Andronikov \cite{andronikov-microperverse} and Waschkies \cite{waschkies-microperverse}.  This ultimately relies on the core observation: the simple-connectedness of the {\em complex} symplectic group provides canonical Maslov data $\xi$ for which we can define a canonical abelian subsheaf $\P erv_V \subset \msh_{V, \xi}$.

The target category $\cE_V - \sMod_{rh}$ of
of regular holonomic $\cE_V$-modules
was introduced by 
Kashiwara~\cite{kashiwara-quantization-contact}. It is based on the key insight: 
while there are many possible choices for 
a global theory of $\cE$-modules,  requiring  compatibility with Verdier duality leads to a canonical choice. 
Regular holonomic modules are characterized microlocally, so within all such $\cE_V$-modules, there is a canonical abelian subsheaf 
$\mathcal{E}_V-mod_{rh}$ of  
regular holonomic $\cE_V$-modules.
Already, Kashiwara asked for a description of the ``stack of `perverse sheaves on $Z$', which is equivalent to the stack Reg($Z$) of regular holonomic systems on $Z$'' \cite{kashiwara-quantization-contact} (here Kashiwara's $Z$ is our complex contact manifold $V$).  
Our Theorem~\ref{thm:gloabl mrh intro}
provides an answer. 
\vspace{2mm} 

Let us mention one particularly important specialization of the above result. 
A complex symplectic 
manifold $\fX$ admits a canonical sheaf of categories of deformation quantization modules $\cW_\fX-\mod$ \cite{polesello-schapira}, built locally from the $\cE$-modules on a contactization.  This category is linear over $\C[[\hbar]]$, but when $\fX$ admits a $\C^\times$-action scaling the symplectic form, it 
is possible to construct an equivariant lift $(\cW_\fX, F_\fX)-\mod$, which is now $\C$-linear, 
and recovers the category of $\cD$-modules on $X$ when $\fX = T^*X$  
\cite{kashiwara-rouquier}.  

These categories have the following significance in modern geometric representation theory.  The fundamental result of Beilinson and Bernstein -- comparing $\mathfrak{g}$-modules to $\cD$-modules on the corresponding flag variety $G/B$ \cite{Beilinson-Bernstein-localisation} -- can be phrased now in terms of $\cW$-modules on $T^*(G/B)$ -- aka the resolution of the nilpotent cone. The natural generalizations to other symplectic resolutions in fact produce categories of previous representation-theoretic interest \cite{kashiwara-rouquier, BPLW-I, BPLW-II}.  As many applications of the Beilinson-Bernstein localisation exploit the Riemann-Hilbert correspondence and properties of perverse sheaves, we expect the following result to be of corresponding utility:  

\begin{theorem}[\ref{thm: w-modules-to-mcsh}] \label{thm: w-modules-to-mcsh-intro} 
Let $\fX$ be a complex manifold 
with a  $\C^\times$ action, equipped with a symplectic form of $\C^\times$-weight $k$, $k \in \Z \setminus \{0\}$. Then there is an equivariant structure $F_\fX$ on the category of $\cW_\fX$-modules, and an equivalence
\begin{equation}\label{equation:main-rh-symplectic-intro} (\cW_\fX, F_\fX)-\mod_{rh} = (\mcsh_{\fX, \mathbb{C}-c}(\fX))^\heartsuit. 
\end{equation}
\end{theorem}

\begin{remark}
    The right hand side of \eqref{equation:main-rh-symplectic-intro} is related to Fukaya categories via the main theorem of \cite{GPS3}. Indeed, let $Z$ be the vector field on $\fX$ induced by differentiating the $\mathbb{C}^*$ action at the identity in the positive real direction, and scaling by $1/k$.  In typical cases of interest (e.g. symplectic resolutions, smooth Higgs moduli with proper Hitchin map) the pair $(\fX, \operatorname{re}i_{Z} \omega)$ is a Liouville manifold in the sense of real symplectic geometry. Suppose now $L \subset \fX$ is a conical (i.e. stable under $Z$) Lagrangian containing all bounded trajectories of $Z$.  
Then, composing Theorem \ref{thm: w-modules-to-mcsh-intro} with \cite{GPS3}, we have an embedding of ordinary (not derived) categories:
\begin{equation}\label{fukaya embedding}
(\cW_\fX, F_\fX)-\mod_{rh, L}(L) \cong \mu_\mathbb{C} sh_{\fX, L}(L)^{\heartsuit} \hookrightarrow H^0(\mu sh_{\fX, L}(L)) \cong H^0(Fuk(X, \partial L)).
\end{equation}

One particular example of interest: 
the celebrated category $\cO$ of Bernstein, Gelfand, and Gelfand is given by D-modules on the flag variety with characteristic cycle in the union of conormals to Schubert cells; more generally, for symplectic resolutions $\fX$ there is a ``category $\cO$ skeleton'' supporting categories of $\mathcal{W}$-modules with similar properties \cite{BPLW-II}.   

From \eqref{fukaya embedding}, we deduce an embedding of these categories into Fukaya categories.    
In some cases, such embeddings were previously known \cite{mak2021fukaya, cote2024hypertoric}, in fact at the level of derived categories.  Our  \eqref{fukaya embedding} lifts to the derived level whenever 
$\mu_\mathbb{C} sh_{\fX, L}(L)$ is the derived category of its (perverse) heart.  Note this is not always true, e.g. for $\fX = T^*S^2$ and $L$ the zero section. 
\end{remark}

\vspace{2mm}

In the rest of the introduction, we provide an overview of the strategy and steps of the proof of Theorem~\ref{thm:gloabl mrh intro}.

To pursue gluing local categories of microsheaves and $\cE$-modules, recall for both sheaves and $\cD$-modules, functors  are given by integral transforms. For example, given manifolds $M, N$ and a sheaf $K \in sh(M \times N)$, there is a functor $sh(M) \to sh(N)$ given by $F \mapsto K_! F := \pi_{N !} (K \otimes \pi_M^* F)$.  Moreover, microsupport is naturally compatible with such operations so they induce a functor on microsheaves. In particular, when $K_!$ is invertible, notably with $ss(K)$ (away from the zero section) the graph of a conic symplectomorphism $\phi$, one obtains an equivalence $K_!: \phi^* \msh \simeq \msh$. In this case, one says $K$ is a sheaf quantization of $\phi$. One similarly arrives at the notion of $\cD$-module quantization.

The gluing data between local categories of microsheaves and $\cE$-modules are such
 sheaf and $\cD$-module quantizations.
  Suppose (falsely) that there were {\em unique} such quantizations of any given conic symplectomorphism.\footnote{This generally false supposition does in fact hold for sheaves with coefficients in the 1-periodicization of $\Z/2\Z$-mod.}  Then on any contact manifold $V$, one could obtain a canonical sheaf  of microsheaves $\msh_V$ as follows: choose any cover of $V$ by Darboux charts, identify these with balls in cosphere bundles, pull back microsheaves from the cosphere bundles, and glue these together via the putatively unique quantizations of the change-of-charts contactomorphisms.  One could  similarly construct a canonical sheaf  of  $\cE$-modules.

But in fact, the quantization of a conic symplectomorphism $\phi$ is not unique.  For sheaves, the freedom was classified in \cite{guillermou}; and for $\cD$-modules, in \cite{dagnolo-schapira-quantization}.  In general, a quantization is given by the choice of `brane structure' on the graph of $\phi$.
As a first step in the proof of Theorem~\ref{thm:gloabl mrh intro}, we explain in Section~\ref{sect:maslov data to gluing data}
how the global `Maslov data' of  of \cite{nadler-shende}  translates into the  local gluing data of sheaf quantizations for overlaps of charts.  This turns out to be a delicate exercise, and moreover, must be done in the general setting of twisted sheaves as recalled in 
Section~\ref{section:review-sheaf-theory}. This generality is necessary
 to match  the canonical orientation datum in \cite{CKNS} with the corresponding twisting by half-forms in \cite{kashiwara-quantization-contact}. 
 
After detailing the above global-to-local translation,  we arrive at the key calculation which motivates the rest of the paper: the verification that the resulting microsheaf kernels for gluing $\P erv_V$ commute appropriately  with Verdier duality (\Cref{prop: untwisted convolution versus verdier}).  Once this is established, we are able to proceed directly as follows. First, we  apply the existing microlocal Riemann-Hilbert \cite{andronikof-microlocal-RH, waschkies-microlocal-RH} to the kernels to obtain kernels for
gluing $\cE$-modules, in particular regular holonomic $\cE$-modules. Then we recall a characterization of 
Polesello~\cite[Theorem 3.3]{polesello-uniqueness}, recalled in Theorem~\ref{uniqueness-of-kashiwara's} in the text, that uniquely distinguishes
Kashiwara's quantization~\cite[Theorem 2]{kashiwara-quantization-contact}  based on its self-duality.
Finally, we match this criterion with the Verdier duality we have established for the microsheaf kernels, and thus conclude the gluing data coincide. 
\vspace{2mm}

\vspace{2mm} {\bf Acknowledgements.} We thank David Ayala,  Benjamin Gammage, Sam Gunningham, Justin Hilburn, Mikhail Kapranov,  Motohico Mulase, Kevin Sackel, Pierre Schapira, Ben Webster and Filip Živanović  for helpful conversations.

Part of this work was conducted while L.C. was supported by NSF grant DMS-2305257, and while LC was a member of the Max Planck Institute for Mathematics. C.K. was partially supported by VILLUM FONDEN grant 37814, and also NSF grant DMS-1928930 when in residence at the SLMath during Spring 2024. D.N. was 
supported by NSF grant DMS-2101466.
V.S. was supported by  Villum Fonden Villum Investigator grant 37814, Novo Nordisk Foundation grant NNF20OC0066298, and Danish National Research Foundation grant DNRF157.

\vspace{2mm} {\bf Notation and conventions.} Throughout the paper, $\cC$ shall denote a stable presentable symmetric monoidal category and  $Pic(\cC)$ shall denote the subcategory of $\cC$ consisting of invertible objects with respect to the symmetric monoidal product. Whenever we discuss microlocal sheaf theory, we always implicitly fix coefficients $\cC$, which will often be left unspecified.  For  Riemann-Hilbert, the relevant choice is $\cC = \C-mod$.  

Given a manifold $M$, we let $\dT^*M:= T^*M - 0_M$. 

Because the set theoretic product of two contact manifolds is not contact, when writing $V \times V^\prime$ for two contact manifolds $V$, $V^\prime$, we will mean the contact manifold
$$V \times V^\prime \coloneqq (\tilde{V} \times \tilde{V^\prime}) / \Rp$$
obtained by quotienting the product of their symplectizations along the diagonal $\Rp$-action. A similar notation will be used for morphisms in contact geometry.

\section{Twisted sheaves and microlocalization}\label{section:review-sheaf-theory}

Here we recall (in somewhat updated language) the notion of `twisted sheaves' from \cite{kashiwara-representation}, and the direct generalization to this context of the microlocal sheaf theory of \cite{kashiwara-schapira}.  

Given a topological space $X$, we denote the category of sheaves on $X$ by $\sh(X;\mathcal{C})$, or just $\sh(X)$. The assignment $X \supseteq U \mapsto sh(U)$ is itself a sheaf of categories on $X$ which we denote by $\sh_X$.

\subsection{Linearity}\label{subsection:linearity}

Let $(\cC,\otimes,1_\cC)$ be a symmetric monoidal category and $X$ a topological space  admitting good covers (all opens and all non-empty intersections of finitely-many opens are contractible).
Consider the constant presheaf of categories which assigns to $U \subseteq X$ the category $\cC$. The sheafification of this presheaf is the sheaf $\Loc_X$, which assigns to $U \subseteq X$ the category $\Loc(U)$ of $\cC$-valued local systems on $U$. Explicitly, this is because the presheaf map which is defined by
\begin{align*}
\cC &\rightarrow \Loc(U) \\
c &\mapsto c_U,
\end{align*}
where $c_U$ is the constant local system on $U$ with stalk $c$, is an isomorphism on an open cover by the locally contractible assumption.

Assume further that $\cC$ is presentable. One can form the category $\PrLst$ of $\cC$-linear presentable categories with $\cC$-linear colimit-preserving functors (or equivalently functors which are left adjoints). The category $\PrLst$ admits a symmetric monoidal structure $\otimes$, sometimes referred as the stable presentable product, and the unit object is given by $\cC$. We will suppress this notation and stop emphasizing $\cC$-linearity when it is clear from the context.

We shall be interested in the category $\sh(X; \PrLst)$ of $\PrLst$-valued sheaves. An important object in this category is $\sh_X$, the sheaf of $\cC$-valued sheaves on $X$. The category $\sh(X; \PrLst)$ inherits the tensor product structure from $\PrLst$. That is, for $\cA, \cB \in \sh(X;\PrLst)$, there is a $\PrLst$-valued sheaf $\cA \otimes \cB$ which is given by the sheafification of the presheaf
$$ U \mapsto \cA(U) \otimes \cB(U).$$ 

\begin{example} \label{loc-action}
The discussion from the above paragraph implies that there is an action, for any open $U \subseteq X$,
$$ \Loc(U) \otimes \cA(U) \rightarrow \cA(U)$$
induced by the composition
$$ \Loc(U) \otimes \cA(U) \rightarrow (\Loc_X \otimes \cA)(U) = \cA(A).$$ 
Here the later identification is induced by the sheafifying the presheaf morphism which, on $U$, is given by the $\cC$-linear structure
$ \cC \otimes \cA(U) \xrightarrow{\sim} \cA(U)$ and the fact that sheafification commutes with tensor product. 
\end{example}
 
\begin{example}
When $\cA = \sh_X$, this recovers the usual action of local system on sheaves. That is, for a local system $\cL \in \Loc_X(U)$ and 
section $\cF \in \cA(U)$, we may form the sheaf $\cL \otimes \cF \in \cA(U)$. 
\end{example}

\subsection{Twisting}

\begin{definition}
A twisting $\cP_\eta$ on $X$ is a local system whose stalks are given by $Pic(\cC)$, the subcategory of invertible objects in $\cC$, which is particular is a(n) ($\infty$-)group.
Such local systems correspond to maps of the form
\begin{equation}\label{equation:class}
\eta: X \rightarrow B Pic(\cC)
\end{equation}
and pulling back along $\left(\{*\} \rightarrow B Pic(\cC) \right)$ recovers $\cP_\eta$. 
\end{definition}

An alternative way to view $\cP_\eta$ is that it is a ($\infty$)-principal bundle. 
Recall the situation one category level down: one can twist a rank $k$ vector bundle $E$ with a principal $G$-bundle after fixing a representation $\rho: G \rightarrow GL(k;\R)$ by forming $E \times^\rho G$. 
In our case, the most natural category acted by $Pic(\cC)$ is $\cC$ itself where the action is given by 
\begin{align*}
i: Pic(\cC) &\rightarrow Aut(\cC) \\
c &\mapsto c \otimes (-)
\end{align*}
as invertible objects.

\begin{definition}
The sheaf $\Loc_X^\eta$ of $\eta$-twisted local systems is the sheafification of the presheaf $\left(U \mapsto \cC \times^i \cP_\eta(U) \right)$. Here the category $\cC \times^i \cP_\eta(U)$ is given by the colimit
$$ \colim \left( \dots  \mathrel{\substack{\textstyle\rightarrow\\[-0.7ex]
  		      \textstyle\rightarrow\\[-0.7ex]
                      \textstyle\rightarrow \\[-0.7ex]
                      \textstyle\rightarrow}} \cC \times Pic(\cC) \times Pic(\cC) \times \cP_\eta(U)    \mathrel{\substack{\textstyle\rightarrow\\[-0.6ex]
                      \textstyle\rightarrow \\[-0.6ex]
                      \textstyle\rightarrow}} \cC \times Pic(\cC) \times \cP_\eta(U)  \rightrightarrows \cC \times \cP_\eta(U) \right).$$
\end{definition}

Note that the assumption of local contractibility implies that for small open sets $U$, $\cP_\eta(U) = Pic(\cC)$, so, on such open sets, sheafification is not needed. 
As a result, one can compute $\Loc_X^\eta$ by choosing a cover $\{U_\alpha\}$ with trivializations $f_\alpha: \cP_\eta |_{U_\alpha} = Pic(\cC)_{U_\alpha}$, similarly to the computation in  differential geometry for principal $G$-bundles and associated vector bundles.
Then $\Loc_X^\eta$ can be computed by the limit 

\begin{equation} \label{equation: gluing description for twisting}
\lim \left( \prod_\alpha \Loc_{U_\alpha} \rightrightarrows  \prod_{\alpha \beta} \Loc_{U_{\beta \alpha}}  
\mathrel{\substack{\textstyle\rightarrow\\[-0.6ex]
                      \textstyle\rightarrow \\[-0.6ex]
                      \textstyle\rightarrow}} \cdots \right)
\end{equation}
where the restriction maps are appropriately twisted by the chosen trivialization $f_\alpha$'s. This justifies the name `twisted local systems'.

Analogously to how we can tensor a vector bundle by a line bundle,   twisted local systems act on sheaves of $\cC$-linear presentable categories.

\begin{definition}
For $\cA \in \sh(X;\PrLst)$, we set $\cA^\eta \coloneqq \cA \otimes \Loc_X^\eta$.    
\end{definition}

The following lemma will be used when studying twisted sheaves in the next Subsection \ref{subsection: twisted sheaves}

\begin{lemma}\label{dualizability-after-twisting-general} 

Assume, for any open set $U \subseteq X$, the category $\cA(U)$ is dualizable and, for any inclusion $U \subseteq V$, the restriction map $i_{U,V}^*: \cA(V) \rightarrow \cA(U)$ is limit-preserving in addition to colimit-preversing. Then, for any $U \subseteq X$, the category $\cA^\tau(U)$ is dualizable for all $U \subseteq X$.
\end{lemma}

\begin{proof}
For any cover $\{U_\alpha\}$ of $U$, the category $\cA^\tau(U)$ is computed as the limit $\cA^\tau(U) = \lim A^\tau(U_\alpha)$. The limit-preserving assumption implies that the restriction map $i_{U,V}^*$ has a left adjoint ${i_{U,V}}_!$ and $\cA^\tau(U)$ can be computed as a colimit
$$ \cA^\tau(U) = \colim \cA^\tau(U_\alpha)$$
where the transition maps are now given by the left adjoints $i_!$'s. Here we use the equivalence $\PrRst = \PrLst^{op}$.
Because of  local contractibility, we can assume the $U_\alpha$'s are contractible and thus $\cA^\tau(U_\alpha) = \cA(U_\alpha)$.
The proposition then follows from \cite[Proposition 6.3.4]{gaitsgory-rozenblyum-DAG-I}.
\end{proof}

\begin{example}\label{twistings-classical}
To recover the notion of twists in  \cite[3.14]{kashiwara-representation}, recall that if $R$ is a (discrete) ring and $\cC = R - mod$, then $Pic(\cC) = \Z \times B(R^\times)$. Fix a cover $\{U_\alpha\}$ and a trivialization 
$$f_\alpha: \cP_\eta |_{U_\alpha} \xrightarrow{\sim} U_\alpha \times \Z \times B (R^\times)$$
for each $\alpha$. On the double overlaps $U_{\beta \alpha}$, the composition $g_{\beta \alpha} \coloneqq f_\beta f_\alpha^{-1}$ corresponds to a map which (by abuse of notation) we still denote by
$$g_{\beta \alpha} \in \Map(U_{\beta \alpha}, Pic(\cC)) = H^0(U_{\beta \alpha}; \Z) \times \Map(U_{\beta \alpha}, B(R^\times)).$$ Thus, upon linearization by $\Loc(X)^\eta \coloneqq \Loc(X) \times_\cC \cP_\eta$, the $g_{\beta \alpha}$ twists the local systems by the identification $L_{\beta \alpha} [n_{\beta \alpha}] \otimes (-)$ on the double overlap where $n_{\beta \alpha} \in \Z$ and $L_{\beta \alpha}$ is a local system with stalk $R$. The fact that $g_{\gamma \beta} g_{\beta \alpha} = g_{\gamma \alpha}$ gives the identification $$L_{\gamma \beta} [n_{\gamma \beta}] \otimes L_{\beta \alpha} [n_{\beta \alpha}] \xrightarrow{\sim} L_{\gamma \alpha} [n_{\gamma \alpha}]$$
and the required conditions in \cite[3.13]{kashiwara-representation} can be deduced similarly from the fact that $\Omega_* B(R^\times) = R^\times$ is discrete.

Even more concretely, the space $\Map(X, B(R^\times))$ has its connected components given by $H^1(X,R)$, and, 
when $X = M$ is a manifold, by assuming $\{U_\alpha\}$ to be a good cover, 
one can canonically trivialize the $L_{\beta \alpha}$'s and the identifications on triple overlaps will be simply given by some usual \v{C}ech 3-cocycles $\{c_{\gamma \beta \alpha}\}$, satisfying conditions $ c_{\delta \gamma \alpha} c_{\gamma \beta \alpha} = c_{\delta \gamma \beta} c_{\delta \beta \alpha}$ as illustrated in \Cref{fig: cech-three-cocycle}.
\begin{figure}
\centering
\begin{tikzpicture}


\node at (-5,2) {$\alpha$};
\node at (-1,2) {$\beta$};
\node at (-5,0) {$\delta$};
\node at (-1,0) {$\gamma$};

\draw [->, thick] (-4.5,2) -- (-1.5,2) node [midway, above] {$ $};
\draw [->, thick] (-1.5,0) -- (-4.5,0) node [midway, above] {$ $};

\draw [->, thick] (-5,0.3) -- (-5,1.7) node [midway, left] {$ $}; 
\draw [->, thick] (-1,1.7) -- (-1,0.3) node [midway, left] {$ $};

\draw [->, thick] (-1.3,0.3) -- (-4.7,1.8) node [midway, left] {$ $}; 
\draw [->, thick] (-4.6,1.6) -- (-1.2,0.1) node [midway, left] {$ $};

\draw [double equal sign distance, thick] (-1.5,1.4) -- (-2.4,1) node [midway, right] {$ $};
\draw [->, thick] (-2.4,1) -- (-2.46,0.98);
\node at (-2.4,1.55) {$c_{\gamma \beta \alpha}$};

\draw [double equal sign distance, thick] (-3.2,0.6) -- (-4.2,0.2) node [midway, right] {$ $};
\draw [->, thick] (-4.2,0.2) -- (-4.26,0.18);
\node at (-4.2,0.75) {$c_{\delta \gamma \alpha}$};


\node at (1,2) {$\alpha$};
\node at (5,2) {$\beta$};
\node at (1,0) {$\delta$};
\node at (5,0) {$\gamma$};

\draw [->, thick] (1.5,2) -- (4.5,2) node [midway, above] {$ $};
\draw [->, thick] (4.5,0) -- (1.5,0) node [midway, above] {$ $};

\draw [->, thick] (1,0.3) -- (1,1.7) node [midway, left] {$ $}; 
\draw [->, thick] (5,1.7) -- (5,0.3) node [midway, left] {$ $};

\draw [->, thick] (4.6,1.8) -- (1.3,0.4) node [midway, left] {$ $}; 
\draw [->, thick] (1.4,0.26) -- (4.7,1.64) node [midway, left] {$ $}; 

\draw [double equal sign distance, thick] (2.8,1.4) -- (1.9,1) node [midway, right] {$ $};
\draw [->, thick] (1.9,1) -- (1.84,0.98);
\node at (1.9,1.55) {$c_{\delta \beta \alpha}$};

\draw [double equal sign distance, thick] (4.6,0.6) -- (3.7,0.2) node [midway, right] {$ $};
\draw [->, thick] (3.7,0.2) -- (3.64,0.18);
\node at (3.7,0.75) {$c_{\delta \gamma \beta}$};

\node at (0,1) {$=$};

\end{tikzpicture}
\caption{\v{C}ech 3-cocycle condition} \label{fig: cech-three-cocycle}  
\end{figure}
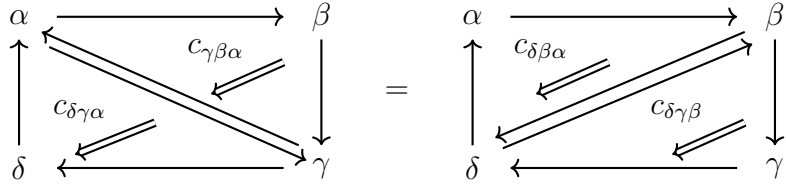

In other words, over $R$, if we decompose $\cA$ to the gluing data $(\cA |_{U_\alpha}, \id_{\cA}, \id_{\id_\cA})$, the sheaf $\cA^\tau$ can be obtained simply by modifying the $2$-morphisms to $(C |_{U_\alpha}, id, c_{\delta \gamma \alpha})$.
Thus an object in $\cA^\tau(X)$ is a family $\{A_\alpha\}$, for $A_\alpha \in \cA(U_\alpha)$, with identifications $g_{\beta \alpha}: A_\alpha [n_{\beta \alpha}] = A_\beta$, on double overlaps, so that $g_{\gamma \beta} \circ g_{\beta \alpha}$ and $g_{\gamma \alpha}$ differ by $c_{\gamma \beta \alpha}$, on triple overlaps.

\end{example}

\begin{remark}
In fact, the description in the above \Cref{twistings-classical} holds more generally.
For two twisting $\eta_1, \eta_2: X \rightarrow B Pic(\cC)$, a homotopy $h: j^* \eta_1 = j^* \eta_2$ on some open set $j: U \hookrightarrow X$, i.e., a path in $Map(U,B Pic(\cC))$ induces an identification of the bundle $h_1: \cP_{\eta_1} |_U = \cP_{\eta_2} |_U$.
Two such homotopies $h_1, h_2$ differs necessarily by a loop $g \coloneqq h_2 \# h_1^\rev \in \Omega_{\eta_1} Map(U, B Pic(\cC)) = Map(U, Pic(\cC))$. The discussion in \Cref{null-homotopies-versus-sections} implies that the auto-equivalence 
$$h_2 \# h_1^\rev: \cP_{\eta_1} = \cP_{\eta_2} = \cP_{\eta_1}$$ is exactly given by tensoring with $g$. The above example is the case when $\eta_2 = *$ is the constant map.
\end{remark}

\subsection{Twisted sheaves} \label{subsection: twisted sheaves}

We specialize to the case $\cA = \sh_X$, the sheaf of sheaves on a locally compact Hausdorff space $X$, and keep the locally contractible assumption. 
For such spaces $X$ and $Y$, colimit-preserving functors from $\sh(Y)$ to $\sh(X)$ are given by sheaf kernels. More precisely, there is an equivalence of categories
\begin{align*} 
\sh(X \times Y) &\xrightarrow{\sim} \Fun^L(sh(X),sh(Y)) \\
K &\mapsto \left(F \mapsto K \circ F \right), 
\end{align*} 
through convolution $K \circ F \coloneqq {\pi_Y}_! (K \otimes \pi_X^* F)$, which follows tracing the computation
$$ \Fun^L(\sh(X),\sh(Y)) = \sh(X)^\vee \otimes \sh(Y) = \sh(X) \otimes \sh(Y) = \sh(X \times Y).$$
Here the first equality follows from abstract properties of dualizable objects, the second equality follows from the fact that $\sh(X)$ is self-dual, and the last equality is  the K\"unneth formula. (See for example  \cite[Corollary 1.3.1.8]{Lurie-SAG}) 

\begin{lemma}[Twisted K\"unneth formula] \label{kunneth}
Let $X$ and $Y$ be locally contractible and locally compact Hausdorff spaces, and $\eta \in \Map(X, BPic(\cC))$ and $\xi \in \Map(Y, B Pic(\cC))$ be twistings.
Then there is an equivalence $\sh^\eta(X) \otimes \sh^\xi(Y) = \sh^{\eta \boxdot \xi}(X \times Y)$
where $\eta \boxdot \xi$ is given by the composition
$$X \times Y \xrightarrow{\eta \times \xi} B Pic(\cC) \times B Pic(\cC) \rightarrow B Pic(\cC)$$
with the last map being the multiplication. 
\end{lemma}

\begin{proof}
The proof is a very standard argument in the set up of $\PrLst$ and can be found in \cite[III.3.1.5]{gaitsgory-rozenblyum-DAG-I}. 
That is, for any choose of over cover $\{U_\alpha \}$, the category $\sh^\tau(X) = \lim \sh^\tau(U_\alpha) = \colim \sh^\tau(U_\alpha)$ is both a limit and a colimit, and thus commutes with $\otimes$. As a result,
$$ \sh^\eta(X) \otimes \sh^\xi(Y)= \lim_{\alpha, \alpha^\prime} \left( \sh^\eta(U_\alpha) \otimes \sh^\xi(V_{\alpha^\prime}) \right)$$
where $\{V_{\alpha^\prime} \}$ is a open cover of $Y$. Assume the covers are given by contractible open sets and the K\"unneth formula follows from the usual one.
\end{proof}

\begin{corollary} \label{twisted-sheaf-convolution}
For locally Hausdorff topological spaces $X$ and $Y$, twists $\eta: X \rightarrow B Pic(\cC)$ and $\xi: Y \rightarrow B Pic(\cC)$, there is an equivalence 
$$\Fun^L(\sh(X)^\eta,\sh(Y)^\xi)  = \sh(X \times Y)^{\eta^{-1} \boxdot \xi}$$
\end{corollary}

\begin{proof}
By \Cref{dualizability-after-twisting-general} and \Cref{kunneth}, it remains to show $\left( \sh^\eta(X) \right)^\vee = \sh^{\eta^{-1}}(X)$. 
To see this one chase through the proof of \cite[Proposition 6.3.4]{gaitsgory-rozenblyum-DAG-I} 
and see that the unit and counit are both given by the diagonal $\Delta: X \hookrightarrow X \times X$.
More precisely, one can twist the inclusion $\sh_\Delta \hookrightarrow \sh_{X \times X}$ of sheaves supported on the diagonal to sheaves on the product by $\eta \boxdot \eta^{-1}$.
Because $\Delta^* (\eta \boxdot \eta^{-1}) = *$ and thus $(\sh_\Delta)^{\eta \boxdot \eta^{-1}} = \sh_\Delta$, we obtain an adjunction
$$  \Delta^*: \sh^{\eta \boxdot \eta^{-1}}(X \times X) \rightleftharpoons \sh(X): \Delta_*,$$
and the unit and counit exhibiting $\left( \sh^{\eta}(X) \right)^\vee = \sh^{\eta^{-1}}(X)$ is given by the two compositions
$$ \cC \xrightarrow{a^*} \sh(X) \xrightarrow{\Delta_*} \sh^{\eta \boxdot \eta^{-1}}(X \times X), \, \text{and} \  \sh^{\eta \boxdot \eta^{-1}}(X \times X) \xrightarrow{\Delta^*} \sh(X) \xrightarrow{a_!} \cC $$
where $a: X \rightarrow \{*\}$ is the projection to the point.
\end{proof}

\begin{example}
Following the notation of \Cref{twistings-classical}, the above identification in the case when $\cC = R - mod$ for a (discrete) ring $R$ is simply given by convolution of the gluing data. That is, for $\{F_\alpha\} \in \sh(X)^{\eta}$ and $\{K_{\alpha^\prime \alpha}\} \in \sh(X \times Y)^{\eta^{-1} \boxdot \xi}$,
one can check directly that the family $\{K_{\alpha^\prime \alpha} \circ F_{\alpha}\}_{\alpha^\prime}$ forms an object in $\sh(Y)^\xi$. 
\end{example}

\begin{example}
Consider a map $f: X \rightarrow Y$ and a twisting $\xi: Y \rightarrow B Pic(\cC)$.
The kernel
$$ 1_{\Gamma_f} \in \sh(X \times Y) = \sh^{f^* \xi \boxdot \xi^{-1}}(X \times Y)$$
gives the adjunction $f^*: \sh^\xi(Y) \rightleftharpoons \sh^{f^* \xi}(X): f_*$.
\end{example}

Note that the same argument as in the proof of \cref{twisted-sheaf-convolution} provides us the notion of tensor product and internal-Hom.
\begin{definition}
Let $\eta_1, \eta_2: X \rightarrow B Pic(\cC)$ be twistings. There is a tensor product 
$$ (-) \otimes (-): \sh^{\eta_1}(X) \otimes \sh^{\eta_2}(X) \rightarrow \sh^{\eta_1 \cdot \eta_2}(X)$$
which is given by the composition $\sh^{\eta_1}(X) \otimes \sh^{\eta_2}(X) \xrightarrow{\sim} \sh^{\eta_1 \boxdot \eta_2}(X \times X) \xrightarrow{\Delta^*} \sh^{\eta_1 \cdot \eta_2}(X)$. By passing to the right adjoint, there is an internal-Hom
$$\sHom(-,-): \sh^{\eta_1}(X)^{op} \otimes \sh^{\eta_2}(X) \rightarrow \sh^{\eta_1^{-1} \cdot \eta_2}(X).$$
In particular, there is a Verdier dual $D_X: \sh^{\eta}(X) \rightarrow \sh^{\eta^{-1}}(X)$ given by Hom-ing into the dualizing sheaf $\sHom(-,\omega_X)$.
\end{definition}

\subsection{Microlocalization for twisted sheaves} \label{twisted-microsheaves}

Let $M$ be a manifold. Since local representatives of a twisted sheaf $F \in \sh(M)^\eta$ differ only by tensoring with a local system, they always have the same microsupport.

\begin{definition}
For $F \in \sh(M)^\eta$, the microsupport of $F$ is defined as the union
$$ss(F) :=  \bigcup_{U, F_U} ss(F_U),$$
where $U$ runs over open sets where $\cP_\eta$ is trivializable and $F_U \in \sh(U)$ runs over local representatives of $F$.
\end{definition}

The notion of microsupport for twisted sheaf locally matches the ordinary notion of microsupport of untwisted sheaves. Therefore, all properties of the microsupport of untwisted sheaves which can be checked locally also hold for twisted sheaves. For example, $ss(F)$ is a closed conic coisotropic subset of $T^* M$.

For a closed conic subset $X \subseteq T^* M$, we use the notation
$$(\sh^\eta)_X(M) \coloneqq \{F \in \sh^\eta(M) | ss(F) \subseteq X \}$$
to denote the subcategory of twisted sheaves microsupported in $X$. 

\begin{definition-lemma} \label{compatibility-twisting-microsupport}
There is an equivalence $(\sh^\eta)_X(M) = (\sh_X)^\eta(M)$. In other words, twisting respects the notion of microsupport and there is an unambiguous meaning of the notation $\sh_X^\eta(M)$.  
\end{definition-lemma}

\begin{proof}
As usual, choose a cover by open balls $\{U_\alpha \}$ and write $(\sh_X)^\eta(M)$ by the limit $\lim \sh_X(U_\alpha)$.
\end{proof}

\begin{example}
By the above definition-lemma, twisted sheaves with microsupport contained in the zero section $\sh^\eta_{0_M}(M) = \Loc^\tau(M)$ are given by twisted local systems. 
\end{example}

For a manifold $M$ and a twist $\eta \in Map\left(M, BPic(\cC) \right)$, we can consider a presheaf of categories
\begin{equation}\label{equation:mupre}
\msh^{\eta, pre}(-): \Omega \mapsto \sh^\eta(M)/\sh_{\Omega^c}^\eta (M).
\end{equation}

\begin{definition} \label{definition-of-microsheaves}
We let $\msh_{T^*M}^\eta $ be the sheafification of the presheaf \eqref{equation:mupre}. This is a sheaf of categories on $T^*M$; the objects of $\msh_{T^*M}^\eta(\Omega)$, for an open set $\Omega \subseteq M$, are called ($\eta$-)twisted microsheaves on $\Omega$. 
The notion of microsupport descends to $\msh_{T^* M}^\eta$ so, for a conic closed subset $X \subseteq T^* M$, we denote by $\msh_X^\tau$ the subsheaf consisting of objects microsupported on $X$. 
\end{definition}
Since the microsupport $ss(F)$ of a sheaf $F \in \sh(M)$ is conic, the restriction of $\msh_{T^* M}^\eta$ to $\dT^* M$ is a  the pullback of a sheaf of categories on the cosphere bundle $S^* M$. We shall denote this later sheaf by $\msh_{S^* M}^\eta$. We use a similar notation $\msh_X^\eta$ for objects microsupported on a closed set $X \subseteq S^* M$. 
Because $T^* M$ retracts to $M$, $\eta$ extends to $T^* M$ uniquely so we abuse the notation to denote the extension by $\eta: T^* M \rightarrow Pic(\cC)$
\begin{lemma} \label{compatibility-twisting-microlocalization} 
There is an equivalence $\msh^\eta_{T^* M} = (\msh_{T^* M})^\eta$. That is, microlocalizing $(\eta-)$twisted sheaves is the same as twisting (by $\eta$) the usual microsheaves.
\end{lemma}

\begin{proof}
Follows from \Cref{compatibility-twisting-microsupport}.    
\end{proof}

We note that by \cite[Proposition 5.4.14]{kashiwara-schapira}, the Verdier dual $D_M: \sh^{\eta}(M) \rightarrow \sh^{\eta^{-1}}(M)^{op}$ microlocalizes to 
\begin{equation} \label{for: mu-verdier}
\mD_M: \msh^\eta_{T^*M} \rightarrow a^* \msh^{\eta^{-1}}_{T^*M}
\end{equation} 
where $a: T^* M \rightarrow T^*M$ is the antipodal map $a(x,\xi) = (x, -\xi)$.
Since Verdier duality restricts to an equivalence
$D_M: \sh_{\R-c}(M)^b = \sh_{\R-c}(M)^{b,op}$ on constructible sheaves with perfect stalks,
by \cite[Remark 7.5.8]{kashiwara-schapira}, its microlocalization $\mD_M$ restricts to an equivalence
$$ \mD_M: (\msh^{\eta}_{T^* M; \R-c})^b =  a^* (\msh^{\eta^{-1}}_{T^* M; \R-c})^{b,op}$$
on microsheaves with Lagrangian microsupport and perfect microstalks.

\subsection{Microlocalizing integral transforms} \label{subsection: microlocalizing integral transform} 

We collect and adapt material from \cite[Section 7.1]{kashiwara-schapira}.  

\begin{definition}[{\cite[Definition 7.1.1]{kashiwara-schapira}}] \label{microlocally-proper-condition}
For manifolds $M$ and $N$ and conic open subsets $\cU \subseteq T^* M$ and $\cV \subseteq T^* N$ one denotes by $N(\cU,\cV)$
the subcategory of $\msh^\pre( T^* M \times \cV)$ consisting of objects $K$ satisfying:
\begin{enumerate}
\item $ss(K) \cap (T^* M \times \cV) \subseteq \cU^a \times \cV$ where $(-)^a$ means the image under the antipodal map $a: (x,\xi) \mapsto (x, - \xi)$.
\item The projection $ss(K) \cap (T^* M \times \cV) \rightarrow \cV$ is proper. 
\end{enumerate}
\end{definition}

\begin{remark}
When $M = \{*\}$, the category $N(\{*\}, \cV) = \msh^\pre(\cV)$.
\end{remark}

\begin{proposition}[{\cite[Proposition 7.1.2, Definition 7.1.3]{kashiwara-schapira}}] \label{microlocal-convolution}
For $i = 1, 2, 3,$ let $M_i$ be a manifold and $\cU_i \subseteq T^* M_i$ be a conic open set.
Let $K \in N(\cU_2,\cU_3)$. Then, $K \circ (-): sh(M_1 \times M_2) \rightarrow sh(M_1 \times M_3)$ which is defined by
$$ K \circ L \coloneqq {p_{13}}_! ( p_{23}^* K \otimes p_{12}^* L)$$
descends to a functor $K \circ (-): \msh^\pre(T* M_1 \times \cU_2) \rightarrow \msh^\pre(T^* M_1 \times \cU_3)$ and the canonical map 
$K \circ L \rightarrow  {p_{13}}_* ( p_{23}^* K \otimes p_{12}^* L)$ is an equivalence.
Furthermore, this functor restricts to $$K \circ (-): N(\cU_1, \cU_2) \rightarrow N(\cU_1, \cU_3).$$ 
In particular, for $K \in N(\cU, \cV)$, there is a functor
$$ K \circ (-): \msh^\pre( \cU) \rightarrow \msh^\pre (\cV),$$
and the canonical morphism
$ K \circ F \rightarrow  {p_2}_* (K \circ \pi_1^* F)$ is an isomorphism for $F \in \msh^\pre(\cU)$.
\end{proposition}

One might want to use the above proposition to define an action of general microsheaf kernels on microsheaves. However, notice that for $\cV_1 \subseteq \cV_2$, the image of $N(\cU, \cV_2)$ under the restriction $\msh^\pre(T^* M \times \cV_1) \rightarrow \msh^\pre(T^* M \times \cV_2)$ is not necessarily contained in $N(\cU, \cV_1)$,  since the intersection of a compact set with an open set is, in general, only closed and condition (2) of \Cref{microlocally-proper-condition} cannot be guaranteed. 

Luckily for us, we will be interested in actions given by microsheaf kernels microsupported on graphs of homogeneous symplectomorphisms and \Cref{microlocal-convolution} can be applied. 
Fix conic open subsets $\cU \hookrightarrow \dT^* M$, $\cV \hookrightarrow \dT^* N$, and a homogenous symplectomorphism $\chi: \cU \xrightarrow{\sim} \cV$.  
Denote by $\Gamma_\chi^a \coloneqq \{\left(x,-\xi,\chi(x,\xi) \right) | (x,\xi) \in \cU\}$ the image of the graph of $\chi$ under $(a \times id)$. 
Consider the subpresheaf $\msh_{\Gamma_\chi^a}^\pre$ of $\msh^\pre$ on $\cU^a \times \cV$ which is given by the assignment
$$ (\Omega \subseteq \cU^a \times \cV) \mapsto \{ F \in \msh^\pre(\Omega) | ss_\Omega(F) \subseteq \Gamma_{\chi}^a \}.$$ 
Here, we use the notation $ss_\Omega(F) \coloneqq ss(F) \cap \Omega$ to denote the microsupport on $\Omega$. 

Denote by $q_i$ the projection of $T^* M \times T^* N$ to the corresponding component. For an open set $\Omega \subseteq \cU^a \times \cV$, set $\cU(\Omega) \coloneqq q_1(\Omega \cap \Gamma_\chi^a)^a$ and $\cV(\Omega) \coloneqq \chi \left(\cU(\Omega) \right)$.
Then observe that, because of the microsupport condition $\Gamma_{\chi}^a$, $\msh_{\Gamma_{\chi}^a}^\pre(\Omega) = \msh_{\Gamma_{\chi}^a}^\pre\left(\cU(\Omega) \times \cV(\Omega)\right)$ and so, combined with the fact that $q_2: \Gamma_\chi^a \rightarrow \cV$ is a diffeomorphism, \Cref{microlocal-convolution} can be applied.
Thus, we obtain an action 
$$\msh_{\Gamma_\chi^a}^\pre( \Omega) \otimes \msh^\pre \left(\cU(\Omega) \right) \rightarrow \msh^\pre \left( \cV(\Omega) \right).$$

\begin{equation}\label{convolution-microsheaf-kernel}
\msh_{\Gamma_\chi^a} \otimes q_1^* \msh_{\cU}  \rightarrow q_2^* \msh_{\cV}
\end{equation}

\begin{proposition}
The morphism $\msh_{\Gamma_\chi^a} \otimes q_1^* \msh_{\cU}  \rightarrow q_2^* \msh_{\cV}$ \ from (\ref{convolution-microsheaf-kernel}) is an equivalence.

\end{proposition}
\begin{proof}
We check at stalks in the sense of sheaves valued in stable categories. That is, for $(x, \xi, y, \eta)$, we need to check that
$\msh_{\Gamma_\chi^a, (x, -\xi, y, \eta)}  \otimes  \msh_{(x,\xi)}  \rightarrow \msh_{(y,\eta)}$ is an equivalence.
We recall that, over stable categories, an objects of $\msh_{(x,\xi)}$ is given by germs of objects in $\msh^\pre(U_0)$ for $U_0 \ni (x,\xi)$. See for example \cite{rozenblyum-filtered}. By \cite[Corollary 7.2.2]{kashiwara-schapira}, on small open sets near $(x, -\xi, y, \eta)$, there exists $\cK$ of microlocal rank $1$, which induces an equivalence $\cK \circ (-): \msh_{(x,\xi)}  \xrightarrow{ \sim } \msh_{(y,\eta)}$ by \cite[Proposition 7.1.10]{kashiwara-schapira}. Thus, the identification of $\cC = \msh_{\Gamma_\chi^a, (x, -\xi, y, \eta)}$ through $1_\cC \mapsto \cK$ identify the desired equivalence as
$$\msh_{\Gamma_\chi^a, (x, -\xi, y, \eta)}  \otimes  \msh_{(x,\xi)} = 1_\cC \otimes  \msh_{(x,\xi)}  = \msh_{(y,\eta)}.$$
\end{proof}

Note that convolution of twisted sheaves from \Cref{twisted-sheaf-convolution} can be microlocalized directly, since microsupport is local in nature as discussed in \Cref{twisted-microsheaves}. In particular, we have the following proposition.
\begin{proposition}
Let $\eta \in Map(M, Pic(\cC))$ and $\xi \in Map(N, Pic(\cC))$ be twistings. Following the same notation as above, there is an action of $(\eta^{-1} \boxdot \xi)$-twisted microsheaf kernel, sending $\eta$-twisted microsheaves on $\cU$ to $\xi$-twisted microsheaves on $\cV$,
$$\msh_{\Gamma_\chi^a}^{\eta^{-1} \boxdot \xi} \otimes q_1^* \msh_{\cU}^\eta  \rightarrow q_2^* \msh_{\cV}^\xi $$
and it is an equivalence.
\end{proposition}

\begin{lemma} \label{microsheaves-on-conormals} 
Let $Z \subseteq M$ be a closed submanifolds and denote by $N^*(Z)$ the conormal bundle of $Z$ in $M$. Then $\msh_{N^*(Z)} = \Loc_{N^*(Z)}$.
\end{lemma}

\begin{proof}
Denote by $\cC_{N^*(Z)}^\pre$ the constant presheaf, i.e., $\cC_{N^*(Z)}^\pre(\Omega) = \cC$ for all $\Omega \subseteq N^*(Z)$. There is a morphism between presheaves whose effect over $\Omega$ is given by 
\begin{align*}
\cC_{N^*(Z)}^\pre(\Omega) = \cC &\rightarrow \msh^\pre_{N^*(\Delta)}(\Omega) \\
C &\mapsto C_Z
\end{align*}
where $C_Z$ is the constant sheaf supported on the submanifold $Z$ with stalk $C$.  
By \cite[Proposition 6.6.1]{kashiwara-schapira}, this map induces equivalences on stalks, and thus becomes an equivalence $\Loc_{N^*(Z)} = \msh_{N^*(Z)}$ upon sheafification.
\end{proof}

\begin{example} \label{action-of-the-diagonal}
We note that $\Gamma_{\id_{T^* M}}^a = N^*(\Delta) \cong T^* M$ where $\Delta \subseteq M \times M$ is the diagonal and the later diffeomorphism can be obtained by dualizing the short exact sequence
$$ 1 \rightarrow TM \xrightarrow{d \Delta} \Delta^* T( M \times M) \rightarrow N(\Delta) \rightarrow 1.$$
By base change, for any sheaf $F \in \sh(M)$, 
$$ C_\Delta \circ F = {p_2}_!( \Delta_! C_M  \otimes p_1^* F) = {p_2}_! \Delta_! (  C_M  \otimes \Delta^* p_1^* F) = C_M \otimes F.$$
That is, the action of $\msh_{N^*(\Delta)}$ on $\msh_{T^* M}$ under the equivalence $\Loc_{T^* M} = \msh_{T^* M}$ is canonically the by local system discussed in \Cref{loc-action}.
\end{example}

\begin{lemma} \label{lem: uniqueness-of-microkernel}
Assume $\cU$ is contractible. Then the category $\msh_{\Gamma_\chi^a}(-\cU \times \cV)$ is non-empty and the functor
\begin{align*}
\msh_{\Gamma_\chi^a}^{\eta^{-1} \boxdot \xi}(-\cU \times \cV) &\rightarrow Fun^{ex} \left(\msh^\eta(\cU),\msh^\xi(\cV) \right) \\
\cK &\mapsto \left(\cF \mapsto \cK \circ \cF \right)
\end{align*}
is conservative. 
\end{lemma}

\begin{proof}
We refer to \Cref{thm: polarization-to-kernel} for the non-emptiness of $\msh_{\Gamma_\chi^a}(-\cU \times \cV)$. We remark that we will need the weaker version discussed in \cite[Theorem 7.2.1]{kashiwara-schapira} and only include this version for completeness.
 
By symmetry, the category $\msh_{\Gamma_{\chi^{-1}}^a}(-\cV \times \cU)$ is non-empty. In fact, the proof of \cite[Theorem 7.2.1]{kashiwara-schapira} implies that the convolution inverse of $\cK$ is given by
$$\cK^* \coloneqq (dv^\vee)^* \mD_{M \times N}(\cK \otimes p_1^* \omega_M),$$
where $v: N \times M = M \times N$ is the coordinates swapping map is the coordinate swapping map.
If two microkernel $\cK_1$, $\cK_2$, induces the same functor, then $\cK_2^* \circ \cK_1$ will acts as the identity. But we know from \Cref{action-of-the-diagonal} that $\cK_2^* \circ \cK_1 = 1_\cU$, since $\cU$ is contractible, and so $\cK_2  = \cK_1$.
\end{proof}

\begin{corollary} \label{cor: loc-when-nonempty}
Let $\cK \in \msh_{\Gamma_\chi^a}(-\cU \times \cV)$ be an invertible object. Then the morphism
\begin{align*}
\Loc_\cU &\rightarrow (a_M \times \chi)^* \msh_{\Gamma_\chi^a}(-\cU \times \cV) \\
L &\mapsto L \otimes \cK
\end{align*}
is an equivalence.
\end{corollary}

\begin{proof}
Since this is morphism between sheaves, it is sufficient to check on an open basis. That is, we only need to evaluate at small enough open balls $\Omega \subseteq \cU$.
In this case, fully-faithfulness is \cite[(7.2.4) in Theorem 7.2.1.]{kashiwara-schapira}. To see it's an surjection, consider any other microkernel $\cK^\prime$.
The above \Cref{lem: uniqueness-of-microkernel} implies that $\cK^* \circ \cK^\prime = c_\cU$ for some $c \in Pic(\cC)$ and thus $\cK^\prime = c \otimes \cK$.
\end{proof}

\section{Integral transforms and Maslov data}
\label{sect:maslov data to gluing data}
We have seen in Section \ref{twisted-microsheaves}, that to a manifold $M$ and a twisting $\eta: M \to Pic(\cC)$, we may associate a sheaf of categories $\msh_{S^*M}^\eta$. More generally, we could have allowed twistings $\eta: S^*M \to Pic(\cC)$.  

In fact, for any contact manifold $\cU$, subject to the vanishing of a certain topological obstruction, the space of twistings has a natural torsor called the space of `Maslov data', and for each Maslov datum $\tau$, there is a corresponding sheaf-of-categories  $\msh_{\cU, \tau}$ on $\cU$ \cite{nadler-shende}.  When $\cU = S^*M$, there is a distinguished Maslov datum $\phi$ arising ultimately from the natural polarization of the cotangent bundle, and an identification of $\msh_{S^*M, \phi}$ with $\msh_{S^*M}$ in the previous  sense  \cite{CKNS}.  

In particular, given a contactomorphism $\chi: \cU \xrightarrow{\sim} \cU'$ along with some Maslov datum $\tau'$ on $\cU'$ and a homotopy $h: \tau \sim \chi^* \tau'$, one obtains an equivalence $h^*: \msh_{\cU, \tau} \xrightarrow{\sim} \chi^* \msh_{\cU', \tau'}$.  
These functors were originally constructed by an abstract argument. 
The purpose of this Section is explain how to more explicitly describe such functors $h^*$ in local charts via integral transforms.

\subsection{Maslov data}
 \label{subsection: maslov-data-and-descent}
We recall results from \cite[Section 11]{nadler-shende}. Some homotopy theory background can be found in \Cref{apn: homotopy-section}.
The central construction is a group homomorphism \cite[(28), Proposition 11.11]{nadler-shende}
\begin{equation}
\fM: U/O \rightarrow B Pic(\cC),
\end{equation}

For any co-oriented contact manifold $(V,\xi)$, we denote by $B Pic(\cC)(\xi)$ the principal $B Pic(\cC)$-bundle classified by the composition 
\begin{equation} \label{classifying-map-contact}
V \xrightarrow{\xi} B U \rightarrow B(U/O) \xrightarrow{B \fM} B^2 Pic(\cC).
\end{equation} 

\begin{definition} 
\
\begin{enumerate}
\item \emph{Maslov data} is a null-homotopy $\tau$ of the map $V \rightarrow B(U/O) \rightarrow B^2 Pic(\cC)$ from (\ref{classifying-map-contact}).
\item A \emph{polarization} is a null-homotopy $\phi$ of the classifying map $V \rightarrow B(U/O)$.
\end{enumerate}
\end{definition}
Evidently a polarization induces Maslov data. If $\phi$ is a polarization, we will routinely abuse notation by also denoting the induced Maslov datum by $\phi$. 

Let $(V, \xi)$ be a contact manifold. Then we can consider the (stable) Lagrangian Grassmannian $\mathfrak{f}: U/O(\xi) \to V$. Its relative cotangent bundle\footnote{If $\pi: E \to B$ is a smooth fiber bundle, then its relative cotangent bundle is the vector bundle $T^*\pi \to E$ defined by the exact sequence $0 \to \pi^*T^*B \to T^*E \to T^*\pi \to 0$.} $T^*\mathfrak{f} \to U/O(\xi)$ is naturally a contact manifold and carries a canonical polarization. Using the large codimensional embedding trick  (reviewed in \Cref{subsubsection:large-codimension} below), one may define microlocal sheaves on $T^*\mathfrak{f}$. The full subcategory of microlocal sheaves supported on the zero section is then a sheaf of categories on $U/O(\xi)$ which we denote by $\msh_{U/O(\xi)}$. As reviewed in \Cref{descending-sheaf-by-equivariance}, this descends to a sheaf of categories on $B Pic(\cC)(\xi)$ which is denoted by $\msh_{B Pic(\cC)(\xi)}$.

\begin{definition}\label{microsheaves-contact}
Given Maslov data $\tau$ on $V$, we define $\msh_{V; \tau} \coloneqq \tau^* \msh_{B Pic(\cC)(\xi)}$. 
\end{definition}

The main feature of the sheaf $\msh_{B Pic(\cC)(\xi)}$ is that it is locally constant in the fiber direction.
In fact, locally on a Darboux chart $\cU$ and along a section, the sheaf $\msh_{B Pic(\cC)(\xi)}$ is non-canonically equivalent to $\msh_{S^* M}$ for some $M$ of the correct dimension. Such a local identification respects the notion of microsupport and we can thus consider, for a closed subset $X \subseteq V$, the subsheaf $\msh_{X;\tau}$ of objects microsupported on $X$.

\begin{proposition}[{\cite[Corollary 4.13]{CKNS}}]\label{microsheaves-by-fiber-polarization} Denote by $\phi_M: S^* M \rightarrow BO$ the (stable) fiber polarization. There is a canonical equivalence \begin{equation}\label{equation:microsh-cosphere}
    \msh_{S^* M; \phi_M} = \msh_{S^* M}
\end{equation}
where  the right hand side is defined as in \Cref{definition-of-microsheaves}. Moreover, the notion of microsupport defined in \Cref{def: microsupport-on-contact} coincides with the usual notion.
\end{proposition}
We shall revisit the proof of \Cref{microsheaves-by-fiber-polarization} in the next subsection, after introducing some more notation. For now we continue our discussion of Maslov data.

Let $\tau_0, \tau_1$ be two choices of Maslov data. A homotopy $h: \tau_0 = \tau_1$ between them is equivalent to a homotopy between the corresponding sections, which we denote it by the same thing. Thus \Cref{microsheaves-contact} implies that there exists an equivalence $\msh_{V; \tau_0} = \msh_{V; \tau_1}$ induced by $h$.

\vspace{2mm}
We can now give a precise formulation of the problem we will solve in this section.  Consider open sets $\cU \subset S^* M$, $\cV \subset S^*N$, and a contactomorphism $\chi: \cU \xrightarrow{\sim} \cV$.
 \Cref{microsheaves-by-fiber-polarization} implies that there is a canonical equivalence $\msh_{\cU; \chi^* \phi_N } = \chi^* \msh_{\cV}$.
In particular, a homotopy of Maslov data $h: \phi_M \simeq \chi^* \phi_N$ induces an equivalence $\msh_\cU = \chi^* \msh_\cV$.

On the other hand, we may use the same homotopy of Maslov data to compare the fiber polarization along the graph of $\chi$ to the fiber polarization along the graph of the diagonal.  Again using \Cref{microsheaves-by-fiber-polarization} along with \Cref{microsheaves-on-conormals}, we find: 
\begin{equation} \label{for: sheaf-quantization-graph}
\msh_{\Gamma_\chi^a} = \msh_{\Gamma_{\id_\cU}^a; \phi_M \boxdot \chi^* \phi_N} \overset{h}{=} \msh_{N^*(\Delta) \cap (\cU^a \times \cU) } = \Loc_\cU.
\end{equation}
Let $\cK \in \msh_{\Gamma_\chi^a}(\cU^a \times \cV)$ be the image of $1_\cU$ under the above equivalence. Then convolution of microsheaf kernels from equation (\ref{convolution-microsheaf-kernel}) induces another equivalence $\msh_\cU = \chi^* \msh_\cV$ by \cite[Proposition 7.1.10]{kashiwara-schapira}. 

The purpose of the present section is to show that these two equivalences agree.  The result is ultimately formulated as \Cref{thm: maslov-to-kernel} below. 

\subsection{Polarization-to-kernel} \label{sec: polarization-to-kernel}

In the present subsection, we shall consider open subsets $\cU \hookrightarrow S^*M$ and $\cV \hookrightarrow S^*N$ and a contactomorphism $\chi: \cU \to \cV$. Let $h: \phi_M \simeq \chi^* \phi_N$ be a homotopy of polarizations. Our goal is to prove that (under the identification \eqref{equation:microsh-cosphere}), the induced equivalence $\msh_{\cU; \phi_M} = \msh_{\cV; \chi_N}$ is induced by convolution with a (twisted) microsheaf kernel, and to characterize this kernel. In the next subsection, we will explain how to generalize this statement to the case where $h$ is merely a homotopy between the Maslov data induced by $\phi_M, \chi^*\phi_N$. 

\subsubsection{Large codimensional embeddings}\label{subsubsection:large-codimension}
The sheaf $\msh_{U/O(\xi)}$ is constructed using the large codimensional embedding method which first appeared in \cite{shende-microlocal} and which we now briefly review. For any co-oriented contact manifold $(V^{2n-1},\xi)$, Gromov's $h$-principle provides contact embeddings
$$i: V \hookrightarrow \R^{2k + 1}$$
for $k \gg 0$ and guarantees that the space of such embeddings can be made as connected as desired by increasing $k$. Such embeddings realize the inverse of the image of $\xi$ in $Map(V, BU)$ as its stable symplectic normal bundle $\nu_V$.  

Recall also that, because of the fiber sequence $O \rightarrow U \rightarrow U/O$, a null-homotopy of $\xi: V \rightarrow BU \rightarrow B(U/O)$ is equivalent to a map $\rho: V \rightarrow BO$ which factorizes $\xi = (BO \rightarrow BU) \circ \rho$. Since $(BO \rightarrow BU)$ is a group, the inverse $\sigma \coloneqq \rho^{-1}$ composes with $BO \rightarrow BU$ to $\nu_V$. Geometrically, this is a Lagrangian sub-bundle $\sigma \subseteq \nu_V$ of the stable symplectic normal bundle (a ``stable normal polarization''). Given such a Lagrangian subbundle, one then \emph{defines} 
\begin{equation}\label{equation:embed-mush}
\msh_{V; \rho}:= \msh_{V^\sigma} |_V.\footnote{Here we once and for all fix an embedding $\R^{2N+1} \hookrightarrow S^*(\R^{N+1})$ to view $\msh_{V^\sigma} |_V$ as the subsheaf  of $\msh_{S^*(\R^{N+1}) }$ on objects microsupported in $V^\sigma$.}
\end{equation} 

A priori, $\msh_{V; \rho}$ depends on the choice of embedding and thickening. To verify that it is well-defined \cite[Lemma 6.3]{nadler-shende} argues by local constancy. 
Namely, let $h: \rho_0 = \rho_1$ be a homotopy between two polarizations. Up to stabilization, this consists of a family of Lagrangian sub-bundles $\sigma_t \subseteq \nu_V$, extending the isotopy of contact embedding $i_t: V \hookrightarrow \R^{2k + 1}$.
The total family $\msh_{V^{\sigma_t} \times [0,1]} |_{V \times [0,1]} $ is constant along the $t$-direction and one thus obtains an equivalence 
$\msh_{V;\rho_0} = \msh_{V;\rho_1}$.

Objects of $\msh_{V;\sigma}$ have a well-defined (micro)support which we typically denote by $ss(-)$.

\begin{definition} \label{def: microsupport-on-contact}
Choose any contact embedding $i: V \hookrightarrow S^* \R^k$, for an object
\begin{equation}\label{equation:microsupport}
    \cF \in \msh_{V; \rho}(V) = \msh_{V^\sigma}\left(i(V) \right) \coloneqq \colim_{\cU \supset i(\cV)} \msh_{V^\sigma}(\cU),
\end{equation}
its (micro)support is the subset $ss(\cF) \coloneqq i(V) \cap ss_{\R^k}(\cF) \subset V$, where $ss_{\R^k}(-)$ refers to the usual notion of microsupport for a microsheaf on $S^* \R^k$.
\end{definition}

\begin{proof}[Proof of \Cref{microsheaves-by-fiber-polarization}]
Pick any embedding $f: M \hookrightarrow \R^k$ for some $k \gg 0$. By \Cref{lemma:exact-normal}
the choice of a section of $f^*: T^* \R^n |_M \rightarrow T^* M$ is contractible and any such section induces an embedding of Liouville manifolds $T^*M \hookrightarrow T^*\mathbb{R}^k$.

Observe also that the complexification of the normal bundle of $M \subset \mathbb{R}^k$ is naturally isomorphic, as a symplectic vector bundle, to the normal bundle of of $T^*M \hookrightarrow T^*\mathbb{R}^k$; hence $-\phi_M = \nu^*$. We let  $S^* M(\nu^*)$ be the thickening of $S^* M$ in the conormal direction. Then the desired map $\msh_{S^* M; \phi_M} \coloneqq \msh_{S^* \R^k; S^* M (\nu^*)} |_{S^* M} \xrightarrow{f^*} \msh_{S^* M}$ is obtained from microlocalizing the equivalence $i^*: sh(\R^k)_{supp(M)} = sh(M)$.

Now to see the statement on microsupport, we choose a local coordinate $(x,v)$ for $\R^k$ such that $M$ is given by $\{v = 0 \}$.
Then $S^* \R^k$ will be given by $(x, v, [\xi, \lambda])$ for $\xi \in T_x^* M$ and $\lambda \in \nu^*_x$. The subset $S^* M (\nu^*) |_{S^* M}$ is then given by $v = 0$ and, which is exactly the image of $\msh_{S^* M; \phi_M} \xrightarrow{\sim}  \msh_{S^* \R^k; S^* M (\nu^*)} |_{S^* M}$.
\end{proof}

$U/O(\xi)$ admits a canonical stable normal polarization and $\msh_{U/O(\xi)}$ is thus defined unambiguously using \eqref{equation:embed-mush}.  For a general contact manifold $V$ equipped with a stable normal polarization $\rho$, we may define $\msh_{V; \rho}$ equivalently either via \eqref{equation:embed-mush} or by pullback $\rho^* \msh_{U/O(\xi)}$ \cite{nadler-shende}.

\subsubsection{GKS quantization} 
We begin by reviewing the GKS quantization method \cite[Theorem 3.7]{guillermou-kashiwara-schapira}. For a contact isotopy $\Phi: S^* M \times I \rightarrow S^* M$, there exists a unique sheaf kernel $K(\Phi) \in \sh(M \times M \times I)$ so that its restriction to the time $0$-slice is $K(\Phi) |_0 = 1_\Delta$ and its microsupport at infinity $ss^\infty(K(\Phi)) \subseteq \Lambda_\Phi$ is contained in the contact movie. A standard microsupport estimate implies that, for any $F \in \sh(M)$, the kernel $K(\Phi)$ moves the microsupport of $F$ by the isotopy, or, when setting $F_t = K(\Phi) |_t \circ F$, we have $ss^\infty(F_t) = \varphi_t( ss^\infty(F))$. A consequence is that, for $t \in I$, convolving with $K(\Phi) |_t$ induces an equivalence $\sh_X(M) = \sh_{\varphi_t(X)}(M)$ for any closed subset $X \subseteq S^* M$.

In particular, suppose that $(\phi_t)_{t \in [0,1]}$ is a family of polarizations for the contact manifolds $V$. Then \cite[Lemma 6.3]{nadler-shende} furnishes a map
\begin{equation}\label{equation:local-constancy}
    \msh_{V; \phi_0} \xrightarrow{\sim} \msh_{V; \phi_1}
\end{equation}
This existence of this map is deduced in loc.\ cit.\ by an abstract argument establishing that $\msh_{V; \phi_t}$ is a locally constant family of categories. In fact, this locally constant family can  be described  concretely via GKS quantization. 

To this end, pick a high-codimensional embedding $V \hookrightarrow \mathbb{R}^{2k+1}$; then the homotopy $\phi_0 \rightsquigarrow \phi_1$ induces a family of thickenings $V_0 \rightsquigarrow V_1$, and we may assume that $V_t = \Phi_t(V_0)$ for a global contact isotopy $\Phi: [0,1] \times S^*\mathbb{R}^{2k+1} \to S^*\mathbb{R}^{2k+1}$; see \Cref{collary:thickening-extend}. We now apply the GKS sheaf quantization, which furnishes a ($[0,1]$-family) sheaf kernel $K(\Phi) \in \sh(\R^k \times \R^k \times [0,1])$. Because the equivalences induced by $K(\Phi)$ on $\sh(M)$ respect microsupport, we obtain an equivalence of sheaves of categories:
\begin{equation}\label{equation:GKS-induced-map-1}
    K(\Phi) \circ (-): \msh_{S^* \R^{N+1}} = \Phi^* \msh_{\Lambda_\Phi \circ S^* M},
\end{equation} where $\Lambda_\Phi \circ S^* \R^k$ is the movie of $S^* M$ under $\Phi$ in $S^* (\R^k \times [0,1])$, and it restricts to 
\begin{equation}\label{equation:GKS-induced-map-2}
   K(\Phi)_1 \circ (-) :\msh_{V_0} = \varphi_1^* \msh_{V_1}
\end{equation}
 where we use the notation $\varphi_t \coloneqq \Phi(-,t)$, $t \in [0,1]$.

\begin{lemma} \label{lem: gks-realizing-constant-family}
    The maps \eqref{equation:local-constancy} and \eqref{equation:GKS-induced-map-2} agree.
\end{lemma}

\begin{proof}
As discussed in \cite[3.4]{guillermou-kashiwara-schapira}, for any $t \in [0,1]$, the restriction
$i_t^*: \sh_{\Lambda_\Phi \circ T^* M}(M \times [0,1]) \rightarrow \sh(M)$ is an equivalence with the inverse of $i_0^*$ given by convolving with $K(\Phi) \circ (-)$.
In fact, \cite[Proposition 3.12]{guillermou-kashiwara-schapira} implies that this equivalence microlocalizes to an equivalence between sheaves of categories
$$K(\Phi) \circ (-): \msh_{S^* \R^k} = \Phi^* \msh_{\Lambda_\Phi \circ S^* \R^k}$$
and the same proposition shows that it restricts to equivalences on subsheaves
$$K(\Phi) \circ (-): \msh_{V_0} = \Phi^* \msh_{\Lambda_\Phi \circ V_0}.$$
Here, the sheaf $\msh_{\Lambda_\Phi \circ V_0}$ is exactly the constant family in which gives rise to the family of sheaves $\{\mu sh_{V; \phi_t}\}_{t \in [0,1]}$ by restricting to the $t$-slice.
Similarly, one can microlocalize $i_t^*$ and obtain 
$$i_t^*: \Phi^* \msh_{\Lambda_\Phi \circ V_0} = \msh_{V_t}.$$ By composing the microlocalization of $K(\Phi)$ with that of $i_1^*$,
we see that the equivalence \eqref{equation:local-constancy} is realized by the equivalence \eqref{equation:GKS-induced-map-2}.
\end{proof}

\begin{remark}
    The proof of \Cref{lem: gks-realizing-constant-family} shows that the action of $\mathcal{K}(\Phi)|_t$ on the sheaf of categories $\mu sh_{V_0}$ is independent of the extension $\Phi$, and only depends on the family of thickenings $V_0 \rightsquigarrow V_1$, and not on the extension to a global contact isotopy. (Indeed, a GKS kernel associated to a contact isotopy which is the identity on $V \subset \mathbb{R}^{2n+1}$ acts by the identity on $\msh_V$). 
\end{remark}

It will also be important to us to consider the case where $V = -\mathcal{U} \times \mathcal{U}$ and we have an isotopy of polariations $\phi_0 \boxdot \phi_0 \rightsquigarrow \phi_0 \boxdot \phi_1$. In this setting, we again have abstractly by \cite{nadler-shende} an induced equivalence
\begin{equation}\label{equation:prov-ns}
     \msh_{-\mathcal{U} \times \mathcal{U}; \phi_0\boxdot \phi_0} \xrightarrow{\sim} \msh_{-\mathcal{U} \times \mathcal{U}; \phi_0 \boxdot \phi_1}
\end{equation}

Meanwhile, there is another way to produce such an equivalence. Namely, consider the family of equivalences 
\begin{equation}\label{equation:prov-prod}
    K(\Phi)|_t \circ (-) : \msh_{-\mathcal{U} \times \mathcal{U}; \phi_0\boxdot \phi_0} \to \msh_{-\mathcal{U} \times \mathcal{U}; \phi_0\boxdot \phi_t}.
\end{equation}

\begin{lemma}\label{lemma:provs-agree}
    \eqref{equation:prov-ns} and \eqref{equation:prov-prod} agree
\end{lemma}
\begin{proof}
    Same argument as the proof of \Cref{lem: gks-realizing-constant-family}, replacing (the microlocalization of) convolution on sheaves
$$K \circ F \coloneqq {p_2}_!( K \otimes \pi_1^* F)$$
by (the microlocalization of) convolution between sheaf kernels
$$K \circ H \coloneqq {p_{13}}_!( \pi_{23}^* K \otimes p_{12}^* H).$$

\end{proof}

We now come to the main result of this section. To set the stage, let $\cU \hookrightarrow S^* M$, $\cV \hookrightarrow S^* N$ be opens, and $\chi: \cU \xrightarrow{\sim} \cV$ be a contactomorphism.
If $h: \phi_M |_\cU = \chi^* (\phi_N |_\cV)$ is a homotopy between polarizations, then we have 
\begin{equation}\label{induced-microsheaf-kernel}
\Loc_{S^* M} = \msh_{N^*(\Delta_M)} = \msh_{N^*(\Delta_M); \phi_M \boxdot \phi_M} = \msh_{N^*(\Delta_M); \phi_M \boxdot \chi^* \phi_N} = \msh_{\Gamma_\chi^a; \phi_M \boxtimes \phi_N} = \msh_{\Gamma_\chi^a}.
\end{equation}

\begin{theorem}\label{thm: polarization-to-kernel}
    Under the identification \eqref{equation:microsh-cosphere}, the induced equivalence $\msh_{\cU;\phi_M} = \msh_{\cV: \chi^* \phi_N}$ is given by microkernel convolution defined by \eqref{convolution-microsheaf-kernel} 
$$\cK(h) \circ (-) : \msh_{\cU} = \chi^* \msh_{\cV}$$
where $\cK(h)$ is the image of $1_{\cU}$ under \eqref{induced-microsheaf-kernel}.
\end{theorem}

\begin{proof}
By \Cref{lemma:provs-agree}, the third equality of \eqref{induced-microsheaf-kernel} can be described via \eqref{equation:prov-prod}. So the image of $1_{\cU}$ under the first three equalities is $K(\Phi)|_1 \circ 1_{\Delta_M}$. Associativity of convolution implies that, for $\cF \in \msh(\cU)$,
$$\cK(h) \circ \cF = (K(\Phi)_1 \circ 1_{\Delta_M}) \circ \cF =K(\Phi)_1 \circ \left(1_{\Delta_M} \circ \cF \right) = K(\Phi)_1 \circ \cF,$$
and we thus see that the equivalence induced by homotopy of polarizations, $\msh_{\cU;\phi_M} = \msh_{\cV: \chi^* \phi_N}$,
is realized by convolving with $\cK(h)$.
\end{proof}

\begin{corollary}\label{cor: polarization-to-kernel-reverse-homotopy}
In the same setting as \Cref{thm: polarization-to-kernel}, denote by $h^\rev: \chi^* (\phi_N |_\cV) =  \phi_M |_\cU$, then
$$\cK(h^\rev) = v^* \mD_{M \times N}(\cK \otimes p_1^* \omega_M)$$
where $v: N \times M = M \times N$ is the coordinates swapping map.
\end{corollary}

\begin{proof}
In the same setting of the proof above. We further notice that, if we write $\Phi^{-1}$ for the isotopy such that $\Phi^{-1}(-,t) = \varphi_t^{-1}$ where $\varphi_t \coloneqq \Phi(-,t)$, then
$$K(\Phi^{-1}) = (v_{\R^k} \times \id_I)^* \sHom(K(\Phi), 1_M \boxtimes \omega_M \otimes 1_I),$$
where $I = [0,1]$, by \cite[Proposition 3.2]{guillermou-kashiwara-schapira}.
\end{proof}

We end this section with a statement concerning the compatibility between Theorem \ref{thm: polarization-to-kernel} and composition.

\begin{lemma} \label{lem: polarization-to-kernel-2-morphism}
In the same setting as \Cref{thm: polarization-to-kernel}, let $h_0, h_1: \phi_M |_\cU = \chi^* (\phi_N |_\cV)$ be two homotopies  between the polarizations. Then, a $2$-homotopy $k: h_0 = h_1$ induces an equivalence $\cK(h_0) = \cK(h_1)$ between microkernels.
\end{lemma}

\begin{proof}
Consider the same high codimensional embedding as in the proof of \Cref{thm: polarization-to-kernel}. The $2$-homotopy $k$ can be realized as a family of Lagrangian thickenings $\sigma_{t,s}$, $t, s \in [0,1]$, such that $\sigma_{t, i}$ is the homotopy between the thickening $S^* M (\nu_M^*) |_\cU$ and $S^* N (\nu_N^*)  |_{\cV}$ given by $h_i$, $i = 0, 1$, and $\sigma_{0, s} = S^* M (\nu_M^*) |_\cU$ and $\sigma_{1,s} = S^* N (\nu_N^*)  |_{\cV}$ for all $s \in [0,1]$. Extend $\sigma_{t,s}$ to an ambient homotopy between isotopies $\Phi: S^* \R^k \times [0,1]_t \times [0,1]_s \rightarrow S^* \R^k$ such that $\Phi_{0, s} = \id_{S^* \R^k}$ and $\Phi_{1,s}$ is a constant contactomorphism for all $s \in [0,1]$. As explained in \cite[Remark 3.9]{guillermou-kashiwara-schapira}, the GKS sheaf quantization applies for any contractible parameter space and we thus have a sheaf kernel $K(\Phi) \in sh(M \times M \times [0,1] \times [0,1])$. By construction and the uniqueness, one recovers $\cK(h_i)$ by setting 
$$ \cK(h_i) \coloneqq K(\Phi) |_{t =1, s = i} \circ 1_{\Delta_M}.$$
However, the constancy condition at $t = 1$ implies that $K(\Phi) |_{t =1} \in \sh(M \times M \times \{1\} \times [0,1]_s)$ is constant on the $s$-direction. Thus $K(\Phi) |_{t =1, s = 0} = K(\Phi) |_{t =1, s = 1}$ canonically and it induces the identification $\cK(h_0) = \cK(h_1)$.
\end{proof}

\begin{corollary} \label{cor: polarization-to-kernel-composition}
Let $\cU_i \hookrightarrow S^* M_i$, $i = 0, 1, 2$ be opens, $\chi_i: \cU_i \xrightarrow{\sim} \cU_{i + 1}$, $i = 0, 1$ be contactomorphisms, $h_i: \phi_{M_i} = \chi_i^* \phi_{M_{i+1}}$, $i = 0, 1$ and $h_2: \phi_{M_0} = (\chi_1 \circ \chi_0)^* \phi_{M_2}$ be homotopy between polarizations. Then a $2$-homotopy $k: h_2 = h_1 \# h_0$ induces an identification
$$ \cK(h_2) = \cK(h_1) \circ \cK(h_0)$$ 
between microkernels and the is identification is compatible with the identification of the corresponding identifications between microsheaf categories induced from identifications between polarizations
$$k: \left(\msh_{\cU_0} \overset{h_0}{=} \msh_{\cU_1} \overset{h_1}{=} \msh_{\cU_2} \right) =  \left( \msh_{\cU_0} \overset{h_2}{=} \msh_{\cU_2} \right).$$
\end{corollary}

\begin{proof}
Consider a similar high codimensional embedding as in the proof of the above \Cref{lem: polarization-to-kernel-2-morphism} and one concludes that $\cK(h_2) = \cK(h_1 \# h_0)$. Call the ambient isotopies which give $h_0$ and $h_1$ by $\Phi$ and $\Psi$. We note that up to a scaling, $\Psi \circ \Phi$ and $\Psi \# \Phi$ are homotopic to each other. Thus, a similar constancy argument at end point implies that the GKS sheaf kernel $K(\Psi \# \Phi) = K(\Psi \circ \Phi) = K(\Psi) \circ K(\Phi)$ where we use the compatibility of GKS sheaf quantization between composition of isotopies for the second quality. The equivalence between microkernels is then induced by applying $(-) \circ 1_{\Delta_M}$. 
\end{proof}

\begin{remark} \label{rmk: polarization-to-kernel-higher-compatibility}
    Higher compatibility, e.g., associativity for compositions hold as well since it follows from the same property for the GKS sheaf quantization. We leave the details to the reader.
\end{remark}

\subsection{Some remarks on principal bundles} \label{remarks principal}

Let $G$ be a group object in spaces and $A$ an abelian group object, i.e., grouplike $E_1$ and $E_\infty$-monoids.
Consider a space $X$ and a map $f: X \rightarrow BG$, i.e., a principal $G$-bundle $G(f) \rightarrow X$.
Assume $\cD$ is a sheaf of categories on $G(f)$ such that $\cD$ is locally constant on the fiber direction and, 
when restricting to any fiber $G \hookrightarrow G(f)$, $\cD |_{G} = L$ for some fixed local system $L$ with stalk category $F$. 
Furthermore, the monodromy of $L$ comes from some representation $\rho: A \rightarrow Aut(F)$. 
In other words, there is an $A$-bundle $A(\alpha) \rightarrow A$, classified by a map $\alpha: G \rightarrow BA$, such that $L = A(\alpha) \times^\rho F$.

Suppose in addition that the map $\alpha$ is in fact a group homomorphism.
In other words, $A(\alpha)$, the fiber of $\alpha$, is a group, as explained in \Cref{kernel-of-space-group}, and we will denote it as $K$. 
One can apply $B(-)$ to $\alpha$ and obtain a map $B(\alpha): BG \rightarrow B^2 A$ in spaces and obtain the following diagram: 

$$
\begin{tikzpicture}
\node at (0,4) {$G(f)$};
\node at (4,4) {$\{*\}$};
\node at (0,2) {$B(A)(f)$};
\node at (4,2) {$BK$};
\node at (8,2) {$\{*\}$};
\node at (0,0) {$X$};
\node at (4,0) {$BG$};
\node at (8,0) {$B^2 A$};

\draw [->, thick] (0.55,4) -- (3.6,4) node [midway, above] {$ $};
\draw [->, thick] (0.85,2) -- (3.6,2) node [midway, above] {$ $};
\draw [->, thick] (0.3,0) -- (3.6,0) node [midway, above] {$f$};

\draw [->, thick] (4.5,2) -- (7.6,2) node [midway, above] {$ $};
\draw [->, thick] (4.5,0) -- (7.5,0) node [midway, above] {$B(\alpha)$};

\draw [->, thick] (0,3.7) -- (0,2.3) node [midway, right] {$f_\alpha$}; 
\draw [->, thick] (4,3.7) -- (4,2.3) node [midway, right] {$ $};

\draw [->, thick] (0,1.7) -- (0,0.3) node [midway, right] {$ $}; 
\draw [->, thick] (4,1.7) -- (4,0.3) node [midway, right] {$ $};
\draw [->, thick] (8,1.7) -- (8,0.3) node [midway, right] {$ $};

\node at (0.4,3.5) {$\ulcorner$};
\node at (0.4,1.5) {$\ulcorner$};
\node at (4.4,1.6) {$\ulcorner$};

\end{tikzpicture}
$$ 
\begin{corollary}
If the monodromy $\alpha: G \rightarrow BA$ is a group homomorphism, then there is a principal $B(A)$-bundle $B(A)(f) \rightarrow X$ such that 
the original principal $G$-bundle $G(f) \rightarrow X$ is a principal $K$-bundle $G(f) \rightarrow B(A)(f)$ over it.
\end{corollary}

The main Theorem \cite[Theorem 11.17]{nadler-shende} is that $\cD$ is pullback from $B(A)(f)$. We sketch the proof since a similar discussion will be relevant.
\begin{proposition} \label{descending-sheaf-by-equivariance}
There is a sheaf of categories $\overline{\cD}$ on $B(A)(f)$ such that $\cD = f_\alpha^* \overline{\cD}$.
\end{proposition}

\begin{proof}
We take the argument from \cite[Theorem 11.17]{nadler-shende}, which concerned the special case (of eventual interest to us) of $\alpha = \mathfrak{M}$. Here we formulate it abstractly for clarity. 

Because $G(f) \rightarrow B(A)(f)$ is a principal $K$-bundle, it is enough to show that $\cD$ is $K$-equivariant. 
That is, there is an identification $p^* \cD = a^* \cD$ (with higher coherence data) where $p, a: K \times G(f) \rightarrow G(f)$ are the trivial projection and the action.

To see this, we recall that the action $a: K \times G(f) \rightarrow G(f)$ is restricted from, $\tilde{a}:G \times G(f) \rightarrow G(f)$, the action given by $G$.
The pullback of $\tilde{a}^* \cD$ equals $L \boxtimes_\rho \cD$. By the definition, $K = A(\alpha)$ is the fiber of $G \xrightarrow{\alpha} BA$, which classifies $A(\alpha)$, and thus the pull back $a^* \cD = A_K \boxtimes \cD$ is trivial.
$$
\begin{tikzpicture}
\node at (0,2) {$K \times A$};
\node at (4,2) {$K$};
\node at (8,2) {$\{*\}$};
\node at (0,0) {$K$};
\node at (4,0) {$G$};
\node at (8,0) {$B A$};

\draw [->, thick] (0.85,2) -- (3.6,2) node [midway, above] {$ $};
\draw [->, thick] (0.3,0) -- (3.6,0) node [midway, above] {$f$};

\draw [->, thick] (4.5,2) -- (7.6,2) node [midway, above] {$ $};
\draw [->, thick] (4.5,0) -- (7.5,0) node [midway, above] {$B(\alpha)$};

\draw [->, thick] (0,1.7) -- (0,0.3) node [midway, right] {$ $}; 
\draw [->, thick] (4,1.7) -- (4,0.3) node [midway, right] {$ $};
\draw [->, thick] (8,1.7) -- (8,0.3) node [midway, right] {$ $};

\node at (0.4,1.5) {$\ulcorner$};
\node at (4.4,1.6) {$\ulcorner$};

\draw [->, thick ] (0.3,-0.2) to[out=-20, in=-160] (7.5,-0.2);
\node at (4,-0.7) {$*$};
\node at (2.5,-0.4) {$\circlearrowleft$};

\end{tikzpicture}
$$ 

\end{proof}

We will also consider the situation when there are two such sheaves $\cD_0$ and $\cD_1$ and a morphism $T: \cD_0 \rightarrow \cD_1$ between them so that, on any fiber, the morphism between the local systems 
$$(L_1 \xrightarrow{T|_G} L_2) = (A(\alpha) \times^\rho F_1 \xrightarrow{A(\alpha) \, \times^\rho \, t} A(\alpha) \times^\rho F_2)$$ are induced from a functor $t: F_1 \rightarrow F_2$ between the stalk categories.
The same argument as the above \Cref{descending-sheaf-by-equivariance} implies that we can descend morphisms as well.

\begin{proposition} \label{descending-morphism-by-equivariance}
There is a morphism between sheaves of categories $\overline{T}: \overline{\cD_0} \rightarrow \overline{\cD_1}$ on $B(A)(f)$ such that $T = f_\alpha^* \overline{T}$.
\end{proposition}

Now consider two sections $s_0$, $s_1$ of $G(f) \rightarrow X$ and a homotopy $h: \overline{s_0} = \overline{s_1}$ between their projection to $G(f) \rightarrow B(A)(f)$. Observe that there exists a unique map $g: X \rightarrow G$ such that $s_1 = g \cdot s_0$. 
The homotopy $h: \overline{s_0} = \overline{s_1} = \alpha(g) \cdot s_0$ thus corresponds to a null-homotopy $k: \alpha(g) = *$. 
Because $K = \fib(\alpha)$, this data $k$ further corresponds to a lifting which we will, by abusing of notation, denote by $k: X \rightarrow K$.

\begin{proposition} \label{lifting-the-equivalence}
Let $s_0$, $s_1$ be two sections of $G(f) \rightarrow X$. A homotopy $h: \overline{s_0} = \overline{s_1}$ induces an equivalence $s_0^* \cD = s_1^* \cD$. A similar statement for morphisms holds. 
\end{proposition}

\begin{proof}
The above discussion implies that $s_1: X \rightarrow G(f)$ can be factored as 
$$X \xrightarrow{(k,s_0)} K \times G(f) \xrightarrow{a} G(f).$$
But, as mentioned in the proof of \Cref{descending-sheaf-by-equivariance}, the pullback $a^* \cD = K_A \boxtimes_\rho \cD$ is constant on the first factor and thus $s_1^* \cD = s_0^* \cD$.
\end{proof}

\subsection{From polarization-to-kernel to Maslov-to-kernel} \label{sec: from-p2k-to-m2k}

The goal of this section is to deduce the following from \Cref{thm: polarization-to-kernel}.
\begin{theorem}\label{thm: maslov-to-kernel}
Let $\cU \hookrightarrow S^* M$, $\cV \hookrightarrow S^* N$ be opens and let $\chi: \cU \xrightarrow{\sim} \cV$ be a contactomorphism.
Suppose that $h: \phi_M |_\cU = \chi^* (\phi_N |_\cV)$ is a homotopy \emph{of Maslov data} (i.e.\ of the Maslov data induced by the polarizations). Then under the identification \eqref{equation:microsh-cosphere}, the induced equivalence $\msh_{\cU;\phi_M} = \msh_{\cV: \chi^* \phi_N}$ is given by microkernel convolution defined by \eqref{convolution-microsheaf-kernel} 
$$\cK(h) \circ (-) : \msh_{\cU} = \chi^* \msh_{\cV}$$
where $\cK(h)$ is the image of $1_{\cU}$ under \eqref{induced-microsheaf-kernel}.
\end{theorem}
 
\begin{proof}
Note that the considerations
of Section \ref{remarks principal} 
apply to 
$\alpha = \fM: U/O \rightarrow B Pic(\cC)$, since this is shown to be a group homomorphism by
\cite[Theorem 11.10, Proposition 11.11]{nadler-shende}. 

The discussion before \Cref{lifting-the-equivalence} implies that there is an equivalence $\tilde{h}: \phi_M |_\cU = k \cdot \chi^* (\phi_N |_\cV)$ for some $k \in \Map(\cU, K)$ where $K$ is the fiber
of the classifying map $\fM: U/O \rightarrow B Pic(\cC)$. 
We cannot apply \Cref{thm: polarization-to-kernel} directly but, since the action of $k$ depends only on its composition to $K \rightarrow U/O$, and the thickening $k \cdot \chi^* (\phi_N |_\cV)$ differs from $\chi^* (\phi_N |_\cV)$ by a rotation by some element in $U(l)$, for some $l \gg 0$, so the same proof still applies.
That is, the equivalence $\msh_\cU = \msh_{\cV; k \cdot \chi^* (\phi_N |_\cV)}$ is induced by microsheaf convolution and, by \Cref{lifting-the-equivalence}, the latter equals $\chi^* \msh_\cV$ trivially.
\end{proof}

\begin{corollary} \label{cor: maslov-to-kernel-composition}
Let $\cU_i \hookrightarrow S^* M_i$, $i = 0, 1, 2$ be opens, $\chi_i: \cU_i \xrightarrow{\sim} \cU_{i + 1}$, $i = 0, 1$ be contactomorphisms, $h_i: \phi_{M_i} = \chi_i^* \phi_{M_{i+1}}$, $i = 0, 1$ and $h_2: \phi_{M_0} = (\chi_1 \circ \chi_0)^* \phi_{M_2}$ be homotopy between Maslov data. Then a $2$-homotopy $k: h_2 = h_1 \# h_0$ induces an identification
$$ \cK(h_2) = \cK(h_1) \circ \cK(h_0)$$ 
between microkernels and the is identification is compatible with the identification of the corresponding identifications between microsheaf categories induced from identifications between Maslov data
$$k: \left(\msh_{\cU_0} \overset{h_0}{=} \msh_{\cU_1} \overset{h_1}{=} \msh_{\cU_2} \right) =  \left( \msh_{\cU_0} \overset{h_2}{=} \msh_{\cU_2} \right).$$
\end{corollary}

\begin{proof}
    Similar to the discussion before \ref{lifting-the-equivalence}, if there are sections $s_i$, $i = 0, 1, 2$ and homotopies $h_i: s_i = s_{i+1}$, $i=0, 1$ and $h_2: s_0 = s_2$. Then a $2$-homotopy $h_2 = h_1 \# h_0$ induces an identification between the liftings. The result then follows from \Cref{cor: polarization-to-kernel-composition} and the proof of \Cref{thm: maslov-to-kernel}.
\end{proof}

\begin{corollary}\label{cor: maslov-to-kernel-twisted}
Let $\eta \in Map(M, Pic(\cC))$ and $\xi \in Map(N, Pic(\cC))$ be twistings. In the same setting as \Cref{thm: maslov-to-kernel}, for a homotopy $h: (\eta \cdot \phi_M) |_\cU = \chi^* ( (\xi \cdot \phi_N) |_\cV)$, the induced equivalence $\msh_{\cU;\eta \cdot \phi_M} = \msh_{\cV: \chi^* (\xi \cdot \phi_N)}$ is realized by convolution
$$\cK(h) \circ (-) : \msh^\eta_{\cU} = \chi^* \msh^\xi_{\cV}$$
using Formula \ref{convolution-microsheaf-kernel}, where $\cK(h)$ is given by an equivalence similar to (\ref{induced-microsheaf-kernel}). Furthermore, denoting by $h^\rev$ the reverse isotopy, then
$$\cK(h^\rev) = v^* \mD_{M \times N}(\cK \otimes p_1^* \omega_M)$$
where $v: N \times M = M \times N$ is the coordinates swapping map. 
\end{corollary}

\begin{proof}
The identification $\msh_{\cU;\eta \cdot \phi_M} = \msh^\eta_{\cU}$ follows from \cite[Corollary 4.13]{CKNS} and \Cref{compatibility-twisting-microlocalization}. To apply \Cref{thm: maslov-to-kernel}, one choose a cover $\{U_\alpha\}$ of $M$ and $\{V_\beta\}$ of $N$ with trivializations of $\eta$ and $\xi$. For each $S^* U_\alpha \cap \cU$ and $S^* V_\beta \cap \cV$, over their intersection through $\chi$, there is an identification $$\phi_M = \eta \cdot \phi_M =  \chi^* (\xi \cdot \phi_N) = \chi^* \phi_N$$ and \Cref{thm: maslov-to-kernel} provides a microsheaf kernel $\cK(h)_{\alpha \beta}$ which realizes the equivalence by convolution. But the same gluing description is used to for \Cref{twisted-sheaf-convolution}, which microlocalizes to Formula \ref{convolution-microsheaf-kernel}. The statement for the reverse homotopy follows from \Cref{cor: polarization-to-kernel-reverse-homotopy} by a similar argument.
\end{proof}

\begin{corollary} \label{cor: maslov-to-kernel-twisted-composition}
Let $\eta_i \in \Map(M_1, Pic(\cC))$, for $i = 0, 1, 2$ be twistings.
In the same setting as \Cref{cor: maslov-to-kernel-twisted-composition}. If $h_i: \eta_i \cdot \phi_{M_i} = \chi_i^* (\eta_{i+1} \cdot \phi_{M_{i+1}})$, $i = 0, 1$ and $h_2: \phi_{M_0} = (\chi_1 \circ \chi_0)^* ( \eta_2 \phi_{M_2})$ are homotopy between Malsov data.
Then a $2$-homotopy $k: h_2 = h_1 \# h_0$ induces an identification
$$ \cK(h_2) = \cK(h_1) \circ \cK(h_0)$$ 
between microkernels and the is identification is compatible with the identification of the corresponding identifications between microsheaf categories induced from identifications between Maslov data
$$k: \left(\msh_{\cU_0}^{\eta_0} \overset{h_0}{=} \msh_{\cU_1}^{\eta_1} \overset{h_1}{=} \msh_{\cU_2}^{\eta_2} \right) =  \left( \msh_{\cU_0}^{\eta_0} \overset{h_2}{=} \msh_{\cU_2}^{\eta_2} \right).$$
\end{corollary}

\begin{proof}
    By trivializing the twistings over the corresponding covers as in the proof of \Cref{cor: maslov-to-kernel-twisted}, we obtain the identification $\cK(h_2) = \cK(h_1) \circ \cK(h_0)$ as (untwisted) microkernels locally by \Cref{cor: maslov-to-kernel-composition}. Now, the transition maps needed for gluing them back as twisted microkernels are given by tensoring with local systems, which are compatible with convolutions and thus the identifications glue as well. 
\end{proof}

\begin{remark} \label{rmk: maslov-to-kernel-twisted-higher-compatibility} 
As mentioned in \Cref{rmk: polarization-to-kernel-higher-compatibility} that polarization-to-kernel satisfied higher compatibility.
A similar argument as the above \Cref{cor: maslov-to-kernel-twisted-composition} implies that it descends to the level of Maslov-to-kernel as well.
\end{remark}

\section{Verdier duality on contact manifolds} 

In \eqref{for: mu-verdier} we observed that Verdier duality microlocalizes to cotangent/cosphere bundles. We shall now explain how to globalize Verdier duality to arbitrary (co-oriented) contact manifolds.

Let $V$ be a co-oriented contact manifold and let $\overline{V}$ denote the same contact manifold with the co-orientation reversed.
For a Maslov datum $\tau$, we will explain the existence of a Verdier dualization functor
\begin{equation}\label{for: global-verdier}
\mD_{V; \tau}: \msh_{V; \tau} \rightarrow \msh_{\overline{V}; \tau}^{op}
\end{equation}
which specializes to \eqref{for: mu-verdier} on Darboux charts. In particular, \eqref{for: global-verdier} restricts to an equivalence on constructible microsheaves with perfect microstalks. 

To set the stage, we recall the standard embedding $(\mathbb{R}^{2k+1}, \sum_{i=0}^{k-1} p_i dq_i+ dq_k) \xhookrightarrow{\iota} S^*\mathbb{R}^{k+1}$ taking $(q_0,\dots, q_k, p_0, \dots, p_{k-1}) \mapsto (q_0, \dots, q_k, p_0,\dots, p_{k-1}, -1)$. We will also consider the opposite embedding $\overline{\mathbb{R}}^{2k+1}= (\mathbb{R}^{2k+1}, -(\sum_0^{k-1} p_i dq_i+ dq_k)) \xhookrightarrow{\overline{\iota}:= a \circ \iota}  S^*\mathbb{R}^{k+1}$, where $a$ is the antipodal map. In the forthcoming discussion, it will be important to remember that Gromov's h-principle applies to \emph{co-oriented} contact manifolds; in particular, a high codimensional embedding $V \hookrightarrow \mathbb{R}^{2k+1}$ naturally induces a high codimensional emebdding $\overline{V} \hookrightarrow \overline{\mathbb{R}}^{2k+1}$.

We begin by assuming that our Maslov datum is induced by a polarization. So let $\rho$ be a polarization of $\xi: V \rightarrow B(U/O)$ and $\sigma \coloneqq \rho^{-1}$ a polarization of the stable normal bundle.
Fix a large codimensional (co-oriented) contact embedding $i: V \hookrightarrow \R^{2k + 1} \xhookrightarrow{\iota} S^*(\R^{k+1} )$. So $\sigma$ is realized as a thickening $V^\sigma$ of $V$, which is Lagrangian in the normal direction, and 
the sheaf $\msh_{V; \rho}$ is given by $\msh_{S^*(\R^{k+1} ); V^\sigma} |_V.$
The microlocal Verdier dual $\mD_{\R^{k+1}}: \msh_{S^*(\R^{k+1} )} \rightarrow \msh_{S^*(\R^{k+1} )}^{op}$
then restricts to
$$\mD_{\R^{k+1}}: \msh_{S^*(\R^{k+1} ); V^\sigma} |_{V} \rightarrow a_{\R^{k+1}}^* \msh_{S^*(\R^{k+1} ); (V^\sigma)^a} |_{V^a}.$$ 
But the target is exactly $\overline{V}$ with the same polarization $\rho$, so we obtain a Verdier dual
\begin{equation}\label{equation:general-verdier-polarization}
\mD_{V; \rho}: \msh_{V; \rho} \rightarrow \msh_{\overline{V}; \rho}^{op}.
\end{equation}

On a cosphere bundle equipped with its fiber polarization, we now have two Verdier duality functors. They agree: 
\begin{proposition}
The identification $\msh_{S^* M; \phi_M} = \msh_{S^* M}$ in \eqref{equation:microsh-cosphere} intertwines the Verdier duality functors \eqref{equation:general-verdier-polarization} and \eqref{for: mu-verdier}. 
\end{proposition}

\begin{proof}
As usual, choose a large codimensional embedding of the base manifold $M \hookrightarrow \R^k$ and choose, up to a contractible space of choices,
a splitting of the projection map $T^* \R^k |_M \rightarrow T^* M$ to obtain an inclusion $i: S^* M \hookrightarrow  S^* \R^k$. 
One first sees that the antipodal map $a_{\R^k}$ on $S^* \R^k$ restricts to that on $M$ so $\msh_{\overline{S^* M}; \phi_M} = a^* \msh_{S^* M}$.
The statement, up to microlocalization, then follows from the fact that, if $\iota: Z \hookrightarrow X$ is a closed inclusion of manifolds, then $\iota^* D_X(F) = D_Z( \iota^* F)$
for $F \in \sh(X)$.
\end{proof}

We now extend the construction of a Verdier duality functor to contact manifolds equipped with a Maslov datum which does not necessarily come from a polarization. Tracing through the stabilization argument in \cite[Section 11.1]{nadler-shende}, we see that the above construction defines a morphism $\mD_{U/O(\xi)}: \msh_{U/O(\xi)} \rightarrow \msh_{\overline{U/O(\xi)}}^{op}$.
We elaborate the notation $\overline{U/O(\xi)}$. A Lagrangian submanifold will remain Lagrangian when changing the symplectic form to its minus. Thus both $V$ and $\overline{V}$ shares the same Lagrangian Grassmannian bundle $\ff: LGr(\xi) \rightarrow V = \overline{V}$, viewed as fiber bundles. However, the relative cotangent bundle $T^* \ff$, for $V$, comes with a co-oriented contact structure $\lambda_f$ described in \cite[Lemma 10.5]{nadler-shende}. For the co-orientation reversing contact manifold $\overline{V}$, we will take the reversed co-orientation $- \lambda_f$ on $T^* \ff$
Now, the antipodal map $a: T^* \ff \rightarrow T^* \ff$ turn $\lambda_f$ to $- \lambda_f$ and fixes the zero section $LGr(\xi) \hookrightarrow T^* \ff$. The discussion from the above paragraph thus produces a morphism 
$\mD_{LGr(\xi)}: \msh_{LGr(\xi)} \rightarrow \msh_{\overline{LGr(\xi)}}^{op}$ where we use $\overline{LGr(\xi)}$ to emphasize that it is a set inside the co-orientation reserving contact manifold. Stabilizing the construction as in \cite[Section 11.1]{nadler-shende} and we obtain the morphism $\mD_{U/O(\xi)}: \msh_{U/O(\xi)} \rightarrow \msh_{\overline{U/O(\xi)}}^{op}$.

Recall for a manifold $M$, we have $\mD_M(\cF \otimes L) = \mD_M (F) \otimes L^\vee$, for any microsheaf $\cF$ and local system $L$ on $S^* M$. Similar formulas holds in various situations, e.g. for the high codimensional embedding construction or when considering the relative cotangent bundle $T^* \ff$, we see that we are in the situation to apply \Cref{descending-morphism-by-equivariance}.
That is, we will apply it to the case when the functor between the stalk categories $t: F_1 \rightarrow F_2$ is given by talking the naive dual $(-)^\vee: \cC \rightarrow \cC$ which is given by $x^\vee \coloneqq \underline{\Hom}(x, 1_\cC)$, and this allows us to descend the morphism $\mD_{U/O(\xi)}$

\begin{equation}
\mD_{B Pic(\cC)}: \msh_{B Pic(\cC)(\xi)} \rightarrow  \msh_{\overline{B Pic(\cC)(\xi)}}^{op}
\end{equation}
and thus, for any Maslov data $\tau$ on $V$, we obtain the morphism $\mD_{V; \tau}$ in (\ref{for: global-verdier}).

We end this section by discussing the compatibility of this global Verdier dual with the sheaf kernels obtained from the process of Maslov-to-kernel.
We first recall a lemma regarding the compatibility of the usual Verdier dual with sheaf kernel convolution.

\begin{lemma} \label{lem: convolution-versus-verdier}
Let $M$ and $N$ be manifolds and $K \in sh(M \times N)$ be a sheaf kernel.
Assume $K$ is constructible with perfect stalks, $ss(K) \cap T^* M \times 0_N \subseteq 0_{M \times N}$, and $supp(K) \hookrightarrow M \times N \rightarrow N$ is proper. 
Then for any $F \in sh(M)$, we have
$$D_N ( K \circ (F)) = D_{M \times N}( K \otimes p_1^* \omega_M) \circ D_M(F).$$
\end{lemma}

\begin{proof}
Standard exercise using base change.
\end{proof}

\begin{proposition} \label{prop: untwisted convolution versus verdier}
In the setting of \Cref{thm: polarization-to-kernel}, there is an equivalence of morphisms
$$\mD_N ( \cK(h) \circ (-)) = \mD_{M \times N}( \cK(h) \otimes p_1^* \omega_M) \circ \mD_M(-): \msh_\cU \rightarrow \chi^* a_N^* \msh_\cV^{op}.$$
\end{proposition}

\begin{proof}
Choose a large codimensional contact embedding $- \cU \times \cV \hookrightarrow S^*( \R^k \times \R^k)$ obtained from a large codimensional embedding of the form $M \times N \hookrightarrow \R^k \times \R^k$
as in the proof of \Cref{thm: polarization-to-kernel}.
Now apply \Cref{lem: convolution-versus-verdier} to the GKS sheaf kernel used in \Cref{thm: polarization-to-kernel}. 
Note that we implicitly use the properness assumption when invoking \cite[Proposition 6.3.3]{kashiwara-schapira}.
\end{proof}

\begin{proposition} \label{prop: micro-convolution-versus-verdier}
In the setting of \Cref{cor: maslov-to-kernel-twisted}, there is an equivalence of morphisms
$$\mD_N ( \cK(h) \circ (-)) = \mD_{M \times N}( \cK(h) \otimes p_1^* \omega_M) \circ \mD_M(-): \msh_\cU^\eta \rightarrow \chi^* a_N^* \msh_\cV^{\xi, op}.$$
\end{proposition}

\section{Microsheaves on complex contact manifolds} \label{sec: msh-complex}
From now on, we assume the coefficients $\cC$ to be $R-mod$ for some (discrete) ring $R$.
We write $\det^2: U \rightarrow S^1 = B \Z$ for the colimit of $\det_n^2: U(n) \rightarrow S^1$. We also consider the second Stiefel-Whitney class $w_2: BO \rightarrow B (\Z/2)$, induced from the exact sequence
$$1 \rightarrow \Z/2 \rightarrow \Pin \rightarrow O \rightarrow 1.$$
Consider the ring homomorphism $\Z \rightarrow R$, which restricts to a group homomorphism $\Z/2 = \Z^\times \rightarrow R^\times$. We shall abuse notation and use $w_2$ to also denote the composition
$$w_2: BO \rightarrow B^2 (\Z/2)\rightarrow B^2 R^\times.$$ 
As recalled in \Cref{subsection: maslov-data-and-descent}, for a real contact manifold $V$, there is a universal sheaf $\msh_{B Pic(R)(\xi)}$ on $ B Pic(R)(\xi)$, in the case when $\cC = R-mod$, the obstruction for $\msh_{B Pic(R)(\xi)}$ to descend to $V$
(here termed the {\em Maslov obstruction}) is given by
$$V \xrightarrow{\xi} BU \rightarrow B(U/O) \rightarrow B^2 Pic(R) = B^2 \Z \times B^3 (R^\times)$$
where the last map is given by the map $B\det^2: B(U/O) \rightarrow B^2 \Z$ and the composition $B(U/O) \rightarrow B^2 O \xrightarrow{B w_2} B^3 R^\times$, and a null-homotopy of the first is referred as a {\em grading} and that of the second an {\em orientation}.
Here we notice that $\det^2$ vanishes on $O$ and the map thus descends.
Note that there is always a canonical orientation $o_{can}$ which is given by composing with the canonical null-homotopy from the fiber sequence $BU \rightarrow B(U/O) \rightarrow B^2 O$. 

\begin{example} \label{eg: can=w2-on-cosphere}
Let $M$ be a (real) manifold. The fiber polarization $\phi_M$ provides both a grading $gr_M$ and an orientation $o_M$. 
We've mentioned in \Cref{microsheaves-by-fiber-polarization} that $\msh_{S^* M; gr_M \times o_M} = \msh_{S^* M}$.
By comparing the difference between $o_M$ with the canonical orientation $o_{can}$, as detailed in \Cref{second-stiefel-whitney}, 
one can conclude that, when replacing the fiber orientation by the canonical one, the associated microsheaves $\msh_{S^* M; gr_M \times o_{can}} = \msh_{S^* M}^{w_2(M)}$ is given by the usual microsheaves twisted by the second Stiefel-Whitney class.

Thus, in the situation when there are open subsets $\cU \hookrightarrow S^* M$ and $\cV \hookrightarrow S^* N$ and a contactomorphism $\chi: \cU \xrightarrow{\sim} \cV$,
one needs only a homotopy $h: o_M |_\cU = \chi^* (o_N |_\cV)$ in order to obtain an identification $\msh_{\cU}^{w_2(M)} = \msh_{\cV}^{w_2(N)}.$
\end{example}

Let $(V, \xi)$ be a real contact manifold. The main observation in \cite[Section 3]{CKNS} is that, if, up to stabilization, $\xi$ is the underlying real symplectic vector bundle of a complex symplectic bundle $E$ on $V$, then an identification $E_\R = \xi$ provides $\xi$ a grading.
Now consider the setting when $(V, \xi)$ is a complex contact manifold and we use $\pi: \tilde{V} \xrightarrow{\C^\times: 1} V$ to denote its complex symplectization, which is a $\C^\times$-bundle over it. 
In this case, the real contact manifold $p: V_0 \coloneqq \tilde{V}/ \Rp$, as discussed in \cite[Section 2]{CKNS}, admits a grading coming from the complex symplectic vector bundle $\xi$.

\begin{definition}[{\cite[Definition 5.1, Remark 5.4.]{CKNS}}]  \label{def: pmsh}
Let $o$ be an orientation of $V_0$. We denote by $\msh_{V_0, o}$ the microsheaves induced by the grading from the complex symplectic structure of $\xi$ and the orientation $o$. When $o = o_{can}$ is the canonical orientation, we write simply $\msh_{V_0}$.
The sheaf $\P \msh_V$, defined in \cite{CKNS}, is a sheaf of full subcategories of $p_* \msh_{V_0}$. 
\end{definition}

As observed in \cite[Lemma 5.6]{CKNS}, a complex contactomorphism $\chi: V \rightarrow U$ between two complex contact manifold is, by definition, a biholomorphic map such that $\chi_* (\xi_V) = \xi_U$, so such a map automatically preserves the complex grading, even before stabilization, and induces canonically an identification $\chi^* \msh_U = \msh_V$. 

In fact, as explained in \cite[Remark 4.20]{CKNS}, for a complex contact manifold $V$, the classifying map
$V_0 \rightarrow B(U/O) \rightarrow B^2 Pic(\cC)$ always factorizes to $V_0 \rightarrow B^2 Pic(\cC)_0$ where we denote by $Pic(\cC)_0$ the connected component of $Pic(\cC)$.
The sheaf $\msh_{V_0}$ is thus a pullback of a universal sheaf $\msh_{B Pic(\cC)_0} (\xi)$ on $B Pic(\cC)_0(\xi)$ whose monodromy on fibers are given by $Pic(\cC)_0$ instead of general objects in $Pic(\cC)$. 

\begin{example} \label{eg: canonical-microsheaf-on-coprojective}
Let $X$ be a complex manifold. Its coprojective bundle $\P^* X$ is a complex contact manifold, whose complex symplectization is the complex cotangent bundle $T^* X$. 
The main observation in \cite[Lemma 3.2]{CKNS} is that its canonical grading is the same as the fiber polarization grading $gr_{X_\R}$, when viewing $X$ as a real manifold $X_\R$. 
Concretely, both the complex symplectic vector bundle $T (T^*X)$ and the real vector bundle $T X_\R$ are obtained from the complex vector bundle $T X$, one by forming the quaternion bundle $(-) \otimes_\C \H$ and the other by taking the underlying real bundle $(-)_\R$,
and both gives rise to the real symplectic vector bundle 
$$T (T^* X_\R) = (T (T^*X))_\R = (T X_\R) \otimes_\R \C.$$

In summary, on the real contact manifold $S^* X_\R$, we have that the canonical microsheaves $\msh_{S^* X_\R} = \msh_{S^* X_\R}^{w_2(X)}$ are given by $w_2$ twisted microsheaves by \Cref{eg: can=w2-on-cosphere}, which we can push to $\P^* X$ to obtain $\P \msh_{P^* X} = \P \msh_{P^* X}^{w_2(X)}$.
Furthermore, if there is a complex contactomorphism $\chi: \cU \rightarrow \cV$ between $\cU \subseteq \P^* X$ and $\cV \subseteq \P^* Y$, there is a canonical identification $\msh_{\cU_0}^{w_2(X)} = \chi^* \msh_{\cV_0}^{w_2(Y)}$, which by \Cref{cor: maslov-to-kernel-twisted}, is induced by a twisted microsheaf kernel
\begin{equation} \label{Kchi} \cK(\chi) \in \P \msh_{\Gamma_\chi^a}^{w_2(X \times Y)}(\cU \times \cV ).
\end{equation}

The next lemma shows that $\cK(h)$ admits extra symmetry. 

\begin{lemma} \label{swapping-coordinates-equal-inverse}
We have a canonical identification $\cK(\chi^{-1}) = v^* \cK(\chi)$ where $v: X \times Y = Y \times X$ is the coordinates swapping map. 
In particular, $\cK(\chi)[n]$ is Verdier self-dual (where $n = \dim_\C X$).
\end{lemma}

\begin{remark}
The case where $\cU$ is a ball is used implicitly in the proof of \cite[Thm.\ 11.4.9]{kashiwara-schapira}.  
\end{remark}

\begin{proof}
The second statement is a direct consequence of the first statement and \Cref{cor: maslov-to-kernel-twisted}. To see the first statement, we apply \cite[Lem.\ 5.6]{CKNS} to the induced map
$$(dv^\vee): P^*(Y \times X) \rightarrow P^*(X \times Y)$$
on the coprojective bundles. This implies that $\P \msh_{P^*(X \times Y)} = (d v^\vee)^* \P \msh_{P^*(Y \times X)}$ canonically and it restricts to the subsheaves
$$\P \msh_{P^*(X \times Y); \Gamma_{\chi}^a} = (d v^\vee)^* \P \msh_{P^*(Y \times X); \Gamma_{\chi^{-1}}^a},$$
and $\cK(\chi)$ is sent to $\cK(\chi^{-1})$ under this identification.
However, realizing the equivalence under the identification $\P \msh_{P^* X} = \P \msh_{P^* X}^{w_2(X)}$ implies that the equivalence is given by
\begin{equation} \label{swapping-microsheaves}
v^*: \P \msh_{P^*(X \times Y)}^{w_2(X \times Y)} =\P \msh_{P^*(Y \times X)}^{w_2(Y \times X)},
\end{equation}
microlocalized from $v^*: sh(X \times Y) = sh(Y \times X)$.
\end{proof}

When there is another $\chi^\prime: \cV \rightarrow \cW$ for some $\cW \hookrightarrow \P^* Z$, the strict identifications of complex symplectic vector bundles $(\chi^\prime \circ \chi)_* \xi_\cU = \chi^\prime_* (\chi_* \xi_\cU) =  \chi^\prime_* \xi_\cV$ implies that there is a canonical identification
\begin{equation} \label{canonical-kernel-composition}
c_{\chi^\prime, \chi}: \cK(\chi^\prime) \circ \cK(\chi) = \cK( \chi^\prime \circ \chi) \in \P \msh_{\Gamma_{\chi^\prime \circ \chi}^a}^{w_2(X \times Z)}( \cU  \times \cW ).
\end{equation}
A similar statement holds for $2$-morphisms between such identification, in that the \v{C}ech $3$-cocycle condition exhibited in \Cref{fig: cech-three-cocycle} holds automatically. 
That is, assume there is a further contactomorphism $\chi^{\prime \prime}: \cW \rightarrow \cQ$, we will need to show that the two different compositions of identifications between microkernels,
$$\cK(\chi^{\prime \prime} ) \circ  \cK(\chi^\prime) \circ \cK(\chi) \xrightarrow{\cK(\chi^{\prime \prime}) \, \bigcirc \, c_{\chi^\prime, \chi}} \cK(\chi^{\prime \prime} ) \circ \cK( \chi^\prime \circ \chi) \xrightarrow{c_{\chi^{\prime \prime}, \chi^\prime \circ \chi}} \cK(\chi^{\prime \prime} \circ \chi^\prime \circ \chi),$$ 
and 
$$\cK(\chi^{\prime \prime} ) \circ  \cK(\chi^\prime) \circ \cK(\chi) \xrightarrow{ c_{\chi^{\prime \prime}, \chi^\prime} \, \bigcirc \, \cK(\chi) } \cK(\chi^{\prime \prime} \circ \chi^\prime)  \circ  \cK( \chi) \xrightarrow{c_{\chi^{\prime \prime} \circ \chi^\prime, \chi}} \cK(\chi^{\prime \prime} \circ \chi^\prime \circ \chi)$$ 
are the same. Per \Cref{rmk: maslov-to-kernel-twisted-higher-compatibility}, it will be implied by the same equality between Maslov data.
But, in this case, it is just the associativity for the identification $\chi^{\prime \prime}_* \chi^\prime_* \chi_* \xi_\cU = \xi_\cQ$.
\end{example}

We end this section with a gluing description for $\P \msh_V$. First, take a cover of $V$ by complex Darboux charts $\{\cU_\alpha\}$ with complex contact embeddings $f_\alpha: \cU_\alpha \hookrightarrow \P^* X_\alpha$ for some complex manifold $X_\alpha$.
The above \Cref{eg: canonical-microsheaf-on-coprojective} implies that $\P \msh_V  |_{\cU_\alpha} = f_\alpha^* \msh_{P^* X_\alpha}^{w_2(\alpha)}$ where we use the notation $w_2(\alpha) \coloneqq w_2(X_\alpha)$, and we will further simplify the notation by $\msh_\alpha^{w_2} \coloneqq f_\alpha^* \msh_{P^* X_\alpha}^{w_2(\alpha)}$. 
The complex contactomorphism, $\chi_{\beta \alpha}: f_\alpha( \cU_{\alpha \beta}) \xrightarrow{\sim} f_\beta( \cU_{\alpha \beta})$ produces a $w_2$-microkernel and we denote its pullback on $V$ by $\cK_{\beta \alpha}$.

On triple overlaps, by \Cref{cor: maslov-to-kernel-twisted-composition}, there is a canonical isomorphism between microkernels  
$$ c_{\gamma \beta \alpha}: \cK_{\gamma \beta} \circ \cK_{\beta \alpha} = \cK_{\gamma \alpha}.$$

\begin{theorem} \label{thm: gluing-description-microsheaf}
Once fixed a Darboux chart $\{(\cU_\alpha, f_\alpha)\}$, the sheaf $\P \msh_V$ can be obtained by the gluing data
\begin{equation} \label{microsheaf-gluing-data}
\left( \{ \P \msh_\alpha\}, \{\cK_{\beta \alpha}\}, \{ c_{\gamma \beta \alpha} \} \right)
\end{equation}
where $\cK_{\beta \alpha}$ are $w_2$-microsheaf kernel satisfying $\mD_{\beta \alpha} (\cK_{\beta \alpha}[2n]) = \cK_{\beta \alpha}$ compatible with the identifications $c_{\gamma \beta \alpha}$.
where $\dim_\C V = 2n -1$ and $\mD_{\beta \alpha}$ is the corresponding Verdier dual defined in (\ref{for: global-verdier}).
\end{theorem}

\begin{remark}
A priori, there is an infinite layer of compatible identifications when gluing a sheaf of categories. 
However, our identifications comes from the homotopical structure of $Pic(R)$, which has $\pi_k (Pic(R)) = 0$ for $k > 1$, so all $k$-morphisms for $k > 2$ are automatically trivial. 
\end{remark}

\begin{proof}
By \Cref{prop: micro-convolution-versus-verdier}, the sheaf $\msh_{\overline{V_0} }$ is glued by $\mD_{\beta \alpha} (\cK_{\beta \alpha}[2n])$. 
However, reversing the co-orientation is invisible on the complex contact level and so $\P \msh_{\overline{V} } = \P \msh_V$ as well.
\end{proof}

\begin{remark}
Because the null-homotopy on the $R^\times$-component already happen before composing with $\Z^\times \rightarrow R^\times$. The above data can in fact be obtained from tensoring the same data over $\Z$ with $(-) \otimes_\Z R$.
\end{remark}

\section{Modules over twisted differential operators} \label{section: modules-over-twisted}

In \Cref{section:review-sheaf-theory}, we discussed twisting sheaves of categories. 
In the special case when the sheaf is formed by taking modules of an algebroid, we can instead twist the algebroid and consider ordinary modules. We will not try to pursue the general theory in this section but restrict ourselves to the case of the ring of differential operators $\cD_X$.

\subsection{Algebroids} \label{Algebroids} Given a sheaf of categories $\mathcal{A}$ on a space $X$, we denote by $\pi_0(\mathcal{A})$ the sheafification of the presheaf of sets $X \supseteq U \mapsto \pi_0(\mathcal{A}(U))$. (Recall that $\pi_0(-)$ of a category is the set of isomorphism classes of objects).  

Recall that a $\C-$algebroid $\cA$ is a sheaf of $\C-$linear 1-categories such that locally $\pi_0(\cA) = \{*\}$, i.e., $\cA$ is locally non-empty and all sections are isomorphic to each other. If $A$ is a sheaf of $\C$-algebras, then one can consider the presheaf of categories $A^{+,pre}$ which assigns an open set $U$ to the category with one object $\{*_U\}$ whose endomorphism is $A(U)$. Its sheafification $A^+$ is a $\C$-algebroid. Although $A^+$ does not contain more data than $A$, there can be more algebroid morphisms than algebra morphisms.

A leisurely account of algebroids and their use in D-module theory can be found in \cite[Sec. 2.1]{Kashiwara-Schapira-DQ}.

\subsection{Line bundles}\label{subsection:line-bundles}

Let $X$ be a complex manifold. We let $Pic(\cO_X)$ be the algebroid whose objects are the (holomorphic) line bundles on $X$, whose morphisms are isomorphisms of line bundles, and whose multiplication is tensor product. More generally, we consider the ``Picard stack'' of $\C$-algebroids $Pic_X(-)$, which assigns to $U \subseteq X$ the algebroid $Pic(\cO_U)$.

\begin{lemma}\label{lemma:pic-O}
$Pic_X(-)= \mathcal{O}_X^+$. Furthermore, under this identification, the canonical ring anti-isomorphism 
$\cO_X \xrightarrow{\sim} \cO_X^{op}$ is identified with taking the dual line bundle
\begin{align*}
Pic_X &\xrightarrow{(-)^\vee} Pic_X \\
L &\mapsto L^{\vee} \left(= L^{-1} \right).
\end{align*}
\end{lemma}
\begin{proof}
For $U \subseteq X$, the category $\cO_X^+(U)$ is given by $Pic(\cO_U)$, the category of (holomorphic) line bundles on $U$. Indeed there is a morphism $\cO_X^{+,pre}(U) \rightarrow Pic(\cO_U)$ which send the point $\{*_U\}$ to the trivial line bundle $\cO_U$ with its effect on morphisms is the identity on $\cO(U)$. This morphism is an equivalence by $\bar{\partial}$-Poincar\'e lemma.
For the second statement, we notice that when reversing the order of multiplication, a transition map $g_{\beta \alpha} \in \cO^\times$ of $L$ will be seen as going the reverse direction. To have the correct direction, one turns $g_{\beta \alpha}$ to $g_{\beta \alpha}^{-1}$ but the later is a cocycle of $L^{-1}$.
\end{proof}

\begin{remark}\label{remark:gluing-data-lb}
Let us recall the \u Cech description of  line bundles, i.e. elements of $\pi_0(Pic(\cO_X)) \simeq H^1(X, \mathcal{O}_X^\times)$. 
It will be useful later to spell this out explicitly. 
Let $U_\alpha$ be a cover of $X$ such that $L |_{U_\alpha}$ admits a section $s_\alpha$ and, on the overlap $U_{\beta \alpha}$, there exists a holomorphic function $g_{\beta \alpha}: U_{\beta \alpha} \rightarrow \C^\times$ such that $s_\beta = g_{\beta \alpha} s_\alpha$. Note that $L$ is a (holomorphic) line bundle implies that $g_{\gamma \beta} g_{\beta \alpha} = g_{\gamma \alpha}$. In addition, the natural equivalence $L \otimes_\cO L^{-1} = \cO_X$ picks out sections $s_{\alpha}^{*}$ such that $\langle s_\alpha^{*}, s_\alpha \rangle = 1$. This later equality implies that $s_\beta^* = g_{\beta \alpha}^{-1} s_\alpha$. In more concise term, the there is an isomorphism of groups
\begin{align*}
\pi_0(Pic(X)) &\xrightarrow{\sim} H^1(X; \cO_X^\times) \\
L& \mapsto [(U_\alpha, g_{\beta \alpha})]
\end{align*}
between holomorphic line bundles and equivalence classes of \v{C}ech cocycles.
\end{remark}

We now discuss ``fractional'' line bundles. Let $G$ be a closed subgroup of $\C^\times$ and fix $\eta: X \to B^2G$. Then we can consider the stack
$$Pic_X^\eta(-):= Pic_X(-) \otimes_{BG} P_\eta.$$
Here $BG$ acts on $Pic(U)$ for any $U \subset X$ by 
$G \to \Omega^* BG \to \Omega^* Aut(Pic(U))= End_{\mathcal{O}_U-mod}(\mathcal{O}_U)= \mathcal{O}_X(U),$ which is the inclusion of $G$ as constant functions into the ring of holomorphic function on $U$.

A \emph{($\eta$-)twisted line bundle} is a global section of $Pic_X^\eta(-)$. Note that we have maps $Pic_X^\eta(-) \otimes Pic_X^\mu(-) \to Pic_X^{\eta+\mu}(-)$. In particular, if $G=\mathbb{Z}/k$, we have $Pic_X^\eta(X)^{\otimes k} \to Pic_X(X)$. 

\begin{example}
If $L$ is a line bundle on $X$, a \emph{$k$-th root} of $L$ is a line bundle $L^{1/k}$ along with an isomorphism $(L^{1/k})^{\otimes k} \simeq L$. Such a $k$-th root need not exist in general as a line bundle. However, we can always construct a $k$-th root as a \emph{twisted} line bundle. 

To do this, consider the exact sequence $0 \to \underline{\mathbb{Z}/k} \to \mathcal{O}_X^* \xrightarrow{z \mapsto z^k} \mathcal{O}_X^* \to 1$. Let $\eta: X \to B^2(\mathbb{Z}/k)$ be image of $L \in H^1(X, \mathcal{O}_X^*) \to H^2(X, \mathbb{Z}/k)$ under the connecting map. 

Then we can define a twisted sheaf $L^{1/k} \in Pic^{\eta}(X)$ concretely as follows: choose a cover $\{U_\alpha\}$ and choose a $k$-th root of the transition functions $f_{\alpha \beta}^{1/k}$. On triple overlaps, the cocycle condition is satisfied up to a $\mathbb{Z}/k$ ambiguity, which (tracing through the definitions) is exactly what is needed to define an object $L^{1/k} \in Pic^{\eta}(X)$.  By construction, $(L^{1/k})^{\otimes k} \simeq L$. 
\end{example}

Note that a sufficient condition for $L$ to admit a $k$-th root as a genuine line bundle is for the image of $L \in H^1(X, \mathcal{O}_X^*) \to H^2(X, \mathbb{Z}/k)$ to vanish. This condition is evidently also necessary since $H^2(X, \mathbb{Z}/k)$ is $k$-torsion. 

\begin{example}
The case $k=2$ shows that a line bundle $L$ admits a square root iff $w_2(L)=0$. In particular, the canonical bundle $\Omega_X$ admits a square root iff $w_2(\Omega_X)=0$ iff $w_2(TX)=0$ iff $X$ is spin. 
\end{example}

\subsection{Twisted differential operators}
Let $X$ be a complex manifold and let $\cD_X$ be the sheaf of differential operators of $X$. This is a sheaf of non-commutative unital rings, and we have inclusions
$$\mathcal{O}_X \subset \cD_X \subset  \cE nd_{\C} (\mathcal{O}_X).$$

The sheaf $\cD_X$ can be defined in multiple equivalent ways. Most concretely a section $D \in \cD_X(U)$ is defined to be a section $D \in Hom_{\C}(\mathcal{O}_X(U), \mathcal{O}_X(U))$ which is locally of the form
$$D(x)= \sum_{\alpha} a_\alpha(x) \partial^\alpha.$$
Here $\alpha=(\alpha_1,\dots, \alpha_k) \in \mathbb{N}^k$ ranges over all multi-indices of length $k \geq 0$, and we write $\partial^\alpha= \partial^{\alpha_1}\dots \partial^{\alpha_n}$. 
We also mention the coordinate-free description, which goes back to Sato, as $\Gamma_{[\Delta]} (\cO_X \boxtimes_\cO \Omega_X)[-n]$ where for complex variety $V \subseteq X$, the functor $\Gamma_{[V]}$ is the temperate support in $V$ \cite[II.5]{bjork-analytic-D}.

Let $L$ be a holomorphic line bundle on $X$. We denote by $\cD_X^L \coloneqq L \otimes_\cO \cD_X \otimes_\cO \otimes L^{-1}$ the $\cO_X$-algebra of differential operators twisted by $L$, whose multiplication is given by,
$$ (L \otimes_\cO \cD_X \otimes_\cO \otimes L^{-1}) \otimes_\cO (L \otimes_\cO \cD_X \otimes_\cO \otimes L^{-1}) = L \otimes_\cO \cD_X \otimes_\cO \cD_X \otimes_\cO \otimes L^{-1} \rightarrow L \otimes_\cO \cD_X \otimes_\cO \otimes L^{-1}$$ where we use $L \otimes_\cO L^{-1} = \cO_X$ for the first equality.

Equivalently and more concretely, $\cD_X^L$ can be described by gluing data. Following the notation of \Cref{remark:gluing-data-lb}, let $U_\alpha$ be a cover of $X$ such that $L |_{U_\alpha}$ admits a section $s_\alpha$ and, on the overlap $U_{\beta \alpha}$, there exists a holomorphic function $g_{\beta \alpha}: U_{\beta \alpha} \rightarrow \C^\times$ such that $s_\beta = g_{\beta \alpha} s_\alpha$. Then on $U_\beta$, a section of $\cD_X^L$ has the form $s_\beta \otimes P \otimes s_\beta^*$, which when restricting to $U_{\beta \alpha}$, becomes the same as 
$$g_{\beta \alpha} s_\alpha \otimes P \otimes g_{\beta \alpha}^{-1} s_\alpha^* 
= s_\alpha \otimes g_{\beta \alpha} P g_{\beta \alpha}^{-1} \otimes s_\alpha^* 
= s_\alpha \otimes \Ad( g_{\beta \alpha}) (P)  \otimes s_\alpha^*.$$ 
In short, $\cD_X^L$ is glued from $\left(\cD_{U_\alpha}, \Ad(g_{\beta \alpha}) \right)$ as an algebra.

We can also twist the sheaf of rings $\cD_X$ by a \emph{fractional} line bundle.

\begin{definition}
We denote by $\cD_X^{\sqrt{L}}$ the sheaf of $\cO$-algebras $L^{1/2} \otimes_\cO \cD_X  \otimes_\cO L^{-1/2}$.
\end{definition}

The expression $L^{1/2} \otimes_\cO \cD_X \otimes_\cO L^{-1/2}$ is a sheaf since $L^{1/2}$ and $L^{-1/2}$ have twistings inverse to each other. More concretely, by shrinking the open cover, we may assume that $g_{\beta \alpha}$ admits a square root $\sqrt{g_{\beta \alpha}}: U_\alpha \rightarrow \C^\times$. Then the sheaf $\cD_X^{\sqrt{L}}$ is glued by the data $\left(\cD_{U_\alpha}, \Ad(\sqrt{g_{\beta \alpha}}) \right)$. We note that since $\Ad(\pm1) = 1$, the expression does not depend on the choice of the square root $\sqrt{g_{\beta \alpha}}$.

\subsection{Twisted $\cD$-modules}
A left/right \emph{$\cD_X$-module} is simply a left/right module over the sheaf of non-commutative rings $\cD_X$.  Similarly one defines left/right $\cD_X^L$-module for any possibly twisted line bundle $L$. When we refer to $\cD_X$-modules without any further adjective, we always mean left $\cD_X$-modules.

\begin{example}[{\cite[2.5.17]{bjork-analytic-D}}] 
Tautologically, the structure ring $\cO_X$ is a left $\cD_X$-module.
Let $T \subseteq X$ be an analytic hypersurface. Then $\cO( * T)$, the sheaf of meromorphic functions with poles contained in $T$, is a left $D_X$-module which contains $\cO_X$ as a submodule.
\end{example}

\begin{example}\label{example: right-canonical}
     Recall that a vector field $v$ acts on top forms $\Omega_X$ by the Lie derivative $L_v$. By Cartan's formula, for $v \in \Theta$ and $\omega \in \Omega_X$, it is simply given by
$$ L_v(\omega) = (d \iota_v + \iota_v d)\omega$$
where $\iota_v$ is the natural contraction of forms by vector fields (at the first component). One can check that 
$$ \omega v \coloneqq -L_v \omega$$
equips $\Omega_X$ with the structure of a right $\cD_X$-module. See for example  \cite[Thm.\ 1.2.14]{bjork-analytic-D} or \cite[Lem.\ 1.8]{kashiwara-microlocal-D}. 
\end{example}

In general, the ring $\cD_X$ is different from $\cD_X^L$. But one categorical level up, the twisting is trivial.
\begin{lemma}[{\cite[Prop.\ 1.9]{kashiwara-microlocal-D} }]\label{lem: equivalence-modules}
There is an equivalence of sheaves $\cD_X^L -mod = \cD_X -mod$.
\end{lemma}

\begin{proof}
The equivalence is simply given by 
\begin{align*}
\cD_X -mod &\xrightarrow{\sim} \cD_X^L -mod \\
\cM &\mapsto L \otimes_\cO \cM
\end{align*}

However, for the purpose of Proposition \ref{twisting-algebras-versus-twisting-modules}, we give a more involved proof from the point of view of gluing. We recall that, for sheaves of sets $G$, $F$, the limit
$$\Hom(G,F) = \lim \left( \prod_{\alpha \in I} \Hom(G_\alpha,F_\alpha) \rightrightarrows \prod_{\alpha, \beta \in I} \Hom(G_{\beta \alpha}, F_{\beta \alpha}) \right),
$$
where the subscript indicates restrictions, states that, in order to define a morphism from $G$ to $F$, it suffices to define them on each $U_\alpha$ and check that they agree on the overlap $U_{\beta \alpha}$. When $C$ and $D$ are sheaves of abelian categories, the existence of non-trivial $2$-morphisms increase the length of the limit by one to 

$$\Hom(C,D) = \lim \left( \prod_{\alpha \in I} \Hom(C_\alpha, D_\alpha) \rightrightarrows \prod_{\alpha, \beta \in I} \Hom(C_{\beta \alpha}, D_{\beta \alpha})
 \mathrel{\substack{\textstyle\rightarrow\\[-0.6ex]
                      \textstyle\rightarrow \\[-0.6ex]
                      \textstyle\rightarrow}}
\prod_{\alpha, \beta, \gamma \in I} \Hom(C_{\gamma \beta \alpha}, D_{\gamma \beta \alpha}) \right).
$$
That is, equality between two functors are now a structure, a natural equivalence, so we have to check whether they are compatible on triple overlap. (Since there is no nontrivial $3$-morphism, equality on that level is again a property.)

By Corollary \ref{equivalence-to-identity}, we have $L_{g_{\beta \alpha} }: \id_{\cD -mod} \xrightarrow{\sim} \forcoe{\Ad(g_{\beta \alpha})}{(-)}$, and we can define a isomorphism from $\cD^L -mod$ to $\cD -mod$ by the following diagram:

$$
\begin{tikzpicture}

\node at (0,2) {$\cD -mod$};
\node at (6,2) {$\cD -mod$};
\node at (0,0) {$\cD -mod$};
\node at (6,0) {$\cD -mod$};

\draw [double equal sign distance, thick] (1,2) -- (5,2) node [midway, above] {$ id$};
\draw [double equal sign distance, thick] (1,0) -- (5,0) node [midway, above] {$ id$};

\draw [double equal sign distance, thick] (0,1.7) -- (0,0.3) node [midway, left] {$\forcoe{\Ad(g_{\beta \alpha})}{(-)}$}; 
\draw [double equal sign distance, thick] (6,1.7) -- (6,0.3) node [midway, left] {$id$};

\draw [double equal sign distance, thick] (3.6,1.2) -- (2.7,0.8) node [midway, right] {$ $};
\draw [->, thick] (2.7,0.8) -- (2.64,0.78);

\node at (2.7,1.45) {$L_{g_{\beta \alpha}}$};

\end{tikzpicture}  
$$
To check that this functor is well-defined, we need to show that $L_{g_{\gamma \beta}} \circ L_{g_{\beta \alpha}} = L_{g_{\gamma \alpha}}$, which follows from the fact that $\{g_{\beta \alpha} \}$ is a \v{C}ech cocycle.

\end{proof}

\begin{proposition} \label{twisting-algebras-versus-twisting-modules}
There is an equivalence of sheaves $\cD_X^{\sqrt{L}} -mod = (\cD_X -mod)^{w_2(L)}$
\end{proposition}

\begin{proof}
Similarly to the proof of Lemma \ref{lem: equivalence-modules}, we would like to untwist $\cD^{\sqrt{L}}-mod$ by the following diagram:
$$
\begin{tikzpicture}

\node at (0,2) {$\cD -mod$};
\node at (6,2) {$\cD -mod$};
\node at (0,0) {$\cD -mod$};
\node at (6,0) {$\cD -mod$};

\draw [double equal sign distance, thick] (1,2) -- (5,2) node [midway, above] {$ id$};
\draw [double equal sign distance, thick] (1,0) -- (5,0) node [midway, above] {$ id$};

\draw [double equal sign distance, thick] (0,1.7) -- (0,0.3) node [midway, left] {$\forcoe{\Ad(\sqrt{g_{\beta \alpha}})}{(-)}$}; 
\draw [double equal sign distance, thick] (6,1.7) -- (6,0.3) node [midway, left] {$id$};

\draw [double equal sign distance, thick] (3.6,1.2) -- (2.7,0.8) node [midway, right] {$ $};
\draw [->, thick] (2.7,0.8) -- (2.64,0.78);

\node at (2.7,1.45) {$L_{\sqrt{g_{\beta \alpha}}}$};

\end{tikzpicture}
$$
On the left side, the natural transformations over the triple overlaps are given by $id$, since $\Ad( \pm 1) = 1$. To absorb the twisting coming from the $L_{ \sqrt{g_{\beta \alpha}} }$'s, those for the right hand side have to be the $c_{\gamma \beta \alpha}$'s. 
\end{proof}

\subsection{The star anti-involution} \label{star-involution-D-modules}
The right module structure of $\Omega_X$ discussed in \Cref{example: right-canonical} provides an equivalence 
\begin{align*}
\Omega_X \otimes_{\cO_{\dT^* X}} (-): \cD_X -\sMod &\xrightarrow{\sim} \cD_X^{op} - \sMod \\
\cM &\mapsto \Omega_X \otimes_{\cO_X} \cM.
\end{align*} 
This follows from the general fact that the tensor product of a left and a right $D$-module over $\cO$ has a right module structure as explained in \cite[(4), 1.3.1 Theorem]{bjork-analytic-D}.
In our case, for $\cM \in \cD_X -mod$, the right module structure on $\Omega_X \otimes_{\cO_X} \cM$ is given by $$ (\omega \otimes m) v \coloneqq (-L_v \omega) \otimes m + \omega \otimes (-v m)$$
for $v \in \Theta$, $\omega \in \Omega$, and $m \in \cM$. One can also view the equivalence as giving by the tautological left $\cD_X^{op}$ and right $\cD_X$-module $\Omega_X \otimes_{\cO_X} \cD_X$ where the `interesting' left $\cD_X^{op}$-module structure is the one described above and the right $\cD_X$-module structure coming from simply multiplying on the right. General Morita theory discussed in \Cref{lem: bimodule-ring-morhpism}, or rather its sheaf version studied in \cite[Section 3]{dagnolo-polesello}, thus produces 
an $\cO_X$-algebra isomorphism 
\begin{equation}\label{equation:star-map}
*: \cD_X \xrightarrow{\sim} \cD_X^{\Omega_X,op}.
\end{equation}

\begin{proposition}
The anti-involution $*: \cD_X  \xrightarrow{\sim} \cD_X^{\Omega_X,op}$ is the usual \emph{star anti-involution} (or \emph{formal adjoint}).
That is, in local coordinates $(x_1,\dots, x_n)$, we have a natural section $dx_1 \wedge \dots \wedge dx_n$ of $\Omega_X$, and \eqref{equation:star-map} is then just the map $P=\sum_\alpha a^\alpha \partial_\alpha \mapsto *(P):= P^*:= \sum_\alpha(-1)^{|\alpha|} \partial_\alpha a^\alpha$.

\end{proposition}

\begin{proof}
Set $d x = dx_1 \wedge \dots \wedge dx_n$. Then $dx \otimes 1 \in \Omega_X \otimes_{\cO_X} \cD_X$ forms a common generator of the two module structures. 
To see its effect, according to \Cref{lem: bimodule-ring-morhpism}, we have to apply $P$ on $dx \otimes 1$ by the interesting $\cD_X^{op}$ action.
But for any $k$, $\partial_k (dx \otimes 1) \coloneqq (-L_{\partial_k}dx \otimes 1) - (dx \otimes 1 \cdot  \partial_k) = - dx \otimes \partial_k$.
Thus, for any monomial $a(x) \partial^\alpha$, we have 
\begin{align*}
(a(x) \partial^\alpha) \cdot (d x \otimes 1) &\coloneqq (a(x) \cdot( \partial^\alpha \cdot (dx \otimes 1) ) 
= a(x) \cdot (dx \otimes ( (-1)^{|\alpha|} \partial^\alpha) ) \\
&=  dx \otimes ( (-1)^{|\alpha|} \partial^\alpha) \cdot a(x) ). 
\end{align*}
\end{proof}

\begin{lemma} \label{lem: compatibility-with-star}
There is an identification $(\cD_X -\sMod)^{w_2(X)} = \cD_X^{\sqrt{\Omega_X} } -\sMod$ where $\cD_X^{\sqrt{\Omega_X} }$ is the ring of differential operators twisted by any square root of $\Omega_X$.
Under this identification, the equivalence $\Omega_X \otimes_{\cO_X} (-): (\cD_X -\sMod)^{w_2(X)} = \cD_X^{op} -\sMod^{w_2(X)}$ is identified with 
\begin{align*}
\cD_X^{\sqrt{\Omega_X} } -\sMod &= \cD_X^{\sqrt{\Omega_X}, op} -\sMod \\
\cM &\mapsto \cM.
\end{align*}
Here we use the fact that any left $\cD_X^{\sqrt{\Omega_X} }$-module automatically admits right $\cD_X^{\sqrt{\Omega_X} }$-module structure to define the map.
\end{lemma} 

\begin{proof}
By \Cref{twisting-algebras-versus-twisting-modules}, the equivalence $\cD_X^{\sqrt{\Omega_X} } -mod = (\cD_X -mod)^{w_2(L)}$ is given by $\Omega_X^{-1/2} \otimes_\cO (-)$. 
Since $(\cD_X^L)^{op} = (\cD_X^{op})^{L^{-1}}$, the same equivalence for the `op' version is given by  $\Omega_X^{1/2} \otimes_\cO (-)$.
But then the composition 
$$\cD_X^{\sqrt{\Omega_X} } -mod = (\cD_X -mod)^{w_2(L)} = (\cD_X^{op} -mod)^{w_2(L)} = \cD_X^{\sqrt{\Omega_X, op} } -mod$$
is given by $\Omega_X^{-1/2} \otimes_{\cO_X} \Omega_X \otimes_{\cO_X} \Omega^{-1/2} \otimes_{\cO_X} (-)  = \cO_X \otimes_{\cO_X} (-)$, which does nothing on the underlying $\cO_X$-module.
Said differently, the bimodule which induces the equivalence $\cD_X^{\sqrt{\Omega_X} } -mod = \cD_X^{\sqrt{\Omega_X, op} } -mod$ is obtained from the bimodule $\Omega_X \otimes_{\cO_X} \cD_X$ by averaging the twisting
$$\Omega_X^{-1/2} \otimes_{\cO_X} \Omega_X \otimes_{\cO_X} \cD_X \otimes_{\cO_X} \Omega^{1/2} = \cD_X^{\sqrt{\Omega_X}}$$
and thus is invisible on the O-module level.
\end{proof}

\section{The canonical sheaf of microdifferential operators}

\subsection{Microlocal operators}

Let $X$ be a complex manifold. There is a sheaf of algebras $\cE_X^\R$ on $T^*X$ whose sections are sometimes called \emph{microlocal operators}. It was originally introduced in the foundational paper of Sato, Kawai, and Kashiwara \cite[II.1.1, II.1.2]{sato-kashiwara-kawai}. 

In this subsection, we summarize the construction and some key properties of $\cE_X^\R$ following \cite[Sec.\ 11.4]{kashiwara-schapira}. We will ultimately be interested in a ``tempered'' variant $\cE_X^{\R, f}$ which will be introduced in the next subsection. However, since both versions enjoy the same formal properties, it seems pedagogically preferable to begin our discussion with $\cE_X^{\mathbb{R}}$.

Let $X_1, X_2, X_3$ be manifolds. Following the notation of \cite[Section 7.1]{kashiwara-schapira}, we let $q_{ij}: X_1 \times X_2 \times X_j \to X_i \times X_j$ be the projection. Similarly we let $p_{ij}:= T^*X_1 \times T^*X_2 \times T^*X_j \to T^*X_i \times T^*X_j$ be the projection on cotangent bundles, and write $p_{ij}^a$ for the composition of $p_{ij}$ with the antipodal map on the $j$-th component. For sheaf kernels $K_1 \in sh(X_1 \times X_2)$, and $K_2 \in sh(X_2 \times X_3)$ we have the convolution
\begin{equation}
    K_1 \circ K_2:= {q_{13}}_! (q_{23}^{*} K_1 \otimes q_{12}^{*} K_2)
\end{equation}
The functor $\mhom(-,-)$ intertwines convolution with tensor product \cite[Proposition 4.4.11]{kashiwara-schapira}. More precisely, given sheaf kernels $K_1, F_1 \in sh(X_1 \times X_2)$ and $K_2, F_2 \in sh(X_ \times X_3)$, we have
\begin{equation}
    {p^a_{13}}_!({p^a_{23}}^{*} \mhom(K_1, F_1) \otimes {p^a_{12}}^{*}\mhom(K_2, F_2)) \to \mhom(K_1 \circ K_2, F_1 \circ F_2).
\end{equation}

Let us now identify $N^*_{\Delta_X}(X \times X) = T^* X$ by projecting onto the first component. We shall also write $\Omega_X \boxtimes_{\cO} \cO_Y:= q_1^* \Omega_X  \otimes_{q_1^* \cO_X} \cO_{X \times Y}$ 

We now consider the sheaf of abelian groups 
\begin{equation} 
\cE_X^\R \coloneqq \mhom(\C_{\Delta_X}[-n], \Omega_X \boxtimes_{\cO} \cO_X) \in  sh(T^* X)^\heartsuit
\end{equation}
It can be shown \cite[Prop.\ 11.4.1, 11.4.4]{kashiwara-schapira} that $\cE_X^\R$ has a ring structure. This is ultimately because of the fact that $\C_{\Delta_X} \circ \C_{\Delta_X} = \C_{\Delta_X}$. The `multiplication' comes from the `integration' morphism \cite[Theorem 11.1.4]{kashiwara-schapira}
\begin{equation} \label{for: integration-map} 
f_! \Omega_Y [\dim_\C Y] \rightarrow \Omega_Z [\dim_\C Z],
\end{equation}
which exists for any holomorphic map $f: Y \rightarrow Z$, specializing to the case of $q_2: X \times X  \rightarrow X$. 
Similarly, for any microsheaf $F$, the sheaf $H^k \mhom(F,\cO_X)$ (resp. $H^k \mhom(F,\Omega_X)$) has left (resp. right) $\cE_X^\R$-module structure. In particular, when taking $F = 1_X$, we see that $\cO_X$ (resp. $\Omega_X$) is a left (resp.) $\cE_X^\R$-module. See also \cite[Proposition 9.3.1]{kashiwara-schapira-microlocal}.

\subsection{From microsheaf kernels to $\mathcal{E}$ kernels} \label{sec: mRH-for-kernel}

Let $\cU \subseteq \dT^* X$, $\cV \subseteq \dT^* Y$ and let $\chi: \cU \xrightarrow{\sim} \cV$ be a complex homogeneous symplectomorphism. 
As discussed in \Cref{eg: canonical-microsheaf-on-coprojective}, there is a canonical (twisted)-microkernel 
$$\cK(\chi) \in \msh^{w_2(X \times Y)}( -\cU \times \cV)$$ such that convolution with $\cK(\chi)$ induces an equivalence
$$\cK(\chi): \msh^{w_2(X)}(\cU) \xrightarrow{\sim} \msh^{w_2(Y)}(\cV)$$
between microsheaves. 
Define 
\begin{equation}\label{equation:hchi-big}
    \cH(\chi) \coloneqq \mhom( \cK(\chi)[-n], \Omega_X \boxtimes_\cO \cO_Y) \in sh^{w_2(X \times Y)}(-\cU \times \cV)^\heartsuit
\end{equation}

For the next proposition, we will use the notation $\cE_X^\R -\sMod$ to denote the sheaf of categories whose sections on $\cU$ is 
$\cE_X^\R -\sMod(\cU) \coloneqq \cE_X^\R|_\cU -mod$. 

\begin{proposition} \label{prop: qct-by-microkernel-ER}
$\cH(\chi)$ is a right $r_1^* \cE_X^\R$ and a left $r_2^* \cE_Y^\R$-module where $r_1: \Gamma_\chi^a \xrightarrow{\sim} \cU$ is the projection and similarly for $r_2$.
Furthermore, tensoring with the bimodule $\cH(h)$ induces an equivalence
$$
\cH(\chi) \otimes_{\cE_X^\R} (-): (\cE_X^\R -\sMod)^{w_2(X)}|_\cU \xrightarrow{\sim}  (\cE_Y^\R - \sMod)^{w_2(Y)} |_\cV,
$$ 
where $\cH(\chi) \otimes_{\cE_X^\R} \cM \coloneqq   {r_2}_* (\cH(\chi) \otimes_{r_1^* \cE_X^\R |_\cU} r_1^* \cM)$, between the (twisted) module categories. 
\end{proposition}

\begin{proof}
The bimodule structure is explained in \cite[Lemma 11.4.3]{kashiwara-schapira}, and is the same structure used to define the ring structure of $\cE_X^\R$. 
To show the equivalence, it is enough by Morita theory to find an inverse $(\cE_X^\R, \cE_Y^\R)$-bimodule. We claim that this bimodule is
$$\cH(\chi^{-1}) \coloneqq \mhom(\cK(\chi^{-1})[n], \Omega_Y \boxtimes_\cO \cO_X).$$
Indeed, the tensor product 
\begin{align*}
\cH(\chi^{-1}) \otimes_{\cE_Y^\R} \cH(\chi) 
&\rightarrow \mhom\left(\cK(\chi^{-1}) [-n] \circ \cK(\chi) [-n],   (\Omega_Y \boxtimes_\cO \cO_X) \circ (\Omega_X \boxtimes_\cO \cO_Y) \right) \\
&\rightarrow \mhom\left(1_{\Delta_X}[-n],  \Omega_X \boxtimes_\cO \cO_X \right) = \cE_X^\R.
\end{align*}
admits a canonical map to $\cE_X^\R$ and similarly to $\cE_Y^\R$.
To define the second arrow, we are using (i) the fact that $\cK(\chi^{-1}) \circ \cK(\chi) = 1_{\Delta_X}$ canonically and (ii) the integration map $(\Omega_Y \boxtimes_\cO \cO_X) \circ (\Omega_X \boxtimes_\cO \cO_Y) \to \Omega_X \boxtimes_\cO \cO_X[-n]$; see \cite[Lem.\ 11.4.3]{kashiwara-schapira}. 
Thus it is sufficient to check on stalks that this map is an equivalence. But this is \cite[Theorem 11.4.9]{kashiwara-schapira} once we restrict to smaller open sets. 
\end{proof}

Similarly to the discussion in \Cref{star-involution-D-modules}, there is an anti-involution on $\cE_X^\R$. The slight complication, as explained in \cite[5.2]{schapira-microdifferential} or \cite[Theorem 7.7]{kashiwara-microlocal-D}, is that the map $*: \cE_X^\R \rightarrow a^* \cE_X^{\R, \Omega_X, op}$ now reverses the cotangent direction, due to the fact that $\xi$ is sent to $-\xi$. Similarly to \Cref{lem: equivalence-modules}, the twisting is invisible at the module level and is given by tensoring with the line bundle $\Omega_X$.

\begin{lemma}
There is an equivalence of sheaves of categories
$$\Omega_X \otimes_{\cO_{T^* X}} (-): \cE_X^\R -\sMod \xrightarrow{\sim} a^* \cE_X^{\R,op} - \sMod$$
whose inverse is given by $\Omega_X^{-1} \otimes _{\cO_{T^* X}} (-)$.
\qed
\end{lemma}

\begin{proposition} \label{prop: qct-by-microkernel-ER-versus-star}
We have the following commuting diagram:

$$
\begin{tikzpicture}
\node at (0,2.5) {$ (\cE_X^\R -\sMod)^{w_2(X)}|_\cU$};
\node at (8,2.5) {$(\cE_Y^\R - \sMod)^{w_2(Y)} |_\cV$};
\node at (-0.2,0) {$ (a^*_X \cE_X^{\R,op} -\sMod)^{w_2(X)}|_{\cU}$};
\node at (8,0) {$ (a^*_Y \cE_Y^{\R,op} - \sMod)^{w_2(Y)} |_{\cV}$};

\draw [->, thick] (1.8,2.5) -- (6.2,2.5) node [midway, above] {$\cH(\chi) \otimes_{\cE_X^\R} (-) $};
\draw [->, thick] (2,0) -- (5.8,0) node [midway, above] {$\cH(\chi^{-1})^{op} \otimes_{\cE_X^{\R,op} } (-)$};

\draw [->, thick] (0,2.2) -- (0,0.3) node [midway, right] {$\Omega_X \otimes_{\cO_{T^* X}} (-)$}; 
\draw [->, thick] (8,2.2) -- (8,0.3) node [midway, right] {$\Omega_Y \otimes_{\cO_{T^* Y}} (-)$};

\node[scale=1.5] at (5,1.4) {$\circlearrowleft$};

\end{tikzpicture}
$$
Here, for the bottom row, we identify $(a^*_X \cE_X^{\R,op} -\sMod)^{w_2(X)}|_{\cU}$ with $( \cE_X^{\R,op} -\sMod)^{w_2(X)}|_{\cU^a}$ and similarly for $Y$.
\end{proposition}

\begin{proof}
We see in the proof of \Cref{prop: qct-by-microkernel-ER} that $\cH(\chi^{-1}) \in \sh( \cV \times \cU)^\heartsuit$ is a sheaf viewed as an $(\cE_X^\R, \cE_Y^\R)$-bimodule using $\chi$.
Since the bimodule structure is given by \cite[Lemma 11.4.3]{kashiwara-schapira}, a microlocal version of convolution, in order to consider its opposite bimodule, we have to swap its coordinate in a way similar to \Cref{lem: opposite-sheaf-bimodule}.
That is, 
$$\cH(\chi^{-1})^{op} = v^* \mhom( \cK(\chi^{-1})[-n], \Omega_Y \boxtimes_\cO \cO_X) \in \sh( \cU \times \cV)^\heartsuit.$$

We thus begin in the bottom left corner of the diagram; we shall argue that traveling around the left, top, and right arrows is the same as traveling along the bottom arrow and we would like to see $\cH(\chi^{-1})^{op}$.

The first itinerary is realized by the composition of bimodules
\begin{align}
\Omega_Y \otimes \cH(\chi) \otimes_{\cE_X^\R} \Omega_X^{-1} &\simeq \mhom(1_Y, \Omega_Y) \otimes \cH(\chi) \otimes_{\cE_X^\R} \mhom(1_X, \Omega_X^{-1}) \nonumber \\
&\simeq \mhom(\cK(\chi)[-n], \cO_X \boxtimes_{\cO} \Omega_Y) \in sh( \cU \times \cV)^\heartsuit \label{equation:intermediate-up-right-left}
\end{align}
 
Applying $v^*$ (see \ref{swapping-microsheaves}) transforms \eqref{equation:intermediate-up-right-left} into
\begin{equation}\label{equation:up-right-down}
    \mhom( v^* \cK(\chi)[-n], \Omega_Y \boxtimes_\cO \cO_X) \in \sh( \cV \times \cU)^\heartsuit.
\end{equation}

But by \Cref{swapping-coordinates-equal-inverse}, \eqref{equation:up-right-down} may be rewritten as
\begin{equation}\label{equation:up-right-down-verdier}
    \mhom( \cK(\chi^{-1})[-n], \Omega_Y \boxtimes_\cO \cO_X) \in \sh( \cV \times \cU)^\heartsuit.
\end{equation}
As explained in the proof of \Cref{prop: qct-by-microkernel-ER}, \eqref{equation:up-right-down-verdier} is the bimodule 
 
$$\cH(\chi^{-1}) \otimes_{\cE_Y^\R}(-): (\cE_Y^\R - \sMod)^{w_2(Y)} |_\cV \xrightarrow{\sim} (\cE_X^\R -\sMod)^{w_2(X)}|_\cU.$$
Finally, we apply $v^*$ (undoing our previous application of $v^*$) to obtain
\begin{equation}
    \cH(\chi^{-1})^{op} \otimes_{\cE_X^{\R,op} } (-): ( \cE_X^{\R,op} -\sMod)^{w_2(X)}|_{\cU^a}\to (\cE_Y^{\R,op} - \sMod)^{w_2(Y)} |_{\cV^a}
\end{equation}
as desired.
\end{proof}

\begin{remark}
    The only non-formal step in the proof of \Cref{prop: qct-by-microkernel-ER-versus-star} is the appeal to \Cref{swapping-coordinates-equal-inverse}. This relies on our assumption that $\chi$ is a \emph{complex} homogeneous symplectomorphism, and can be understood as a manifestation of Verdier duality. \Cref{prop: qct-by-microkernel-ER-versus-star} would be false if we only assumed that $\chi$ was a real homogeneous symplectomorphism.
\end{remark}

\begin{corollary} \label{cor: composition-qct-by-microkernel-ER-versus-star}
In the situation of \Cref{prop: qct-by-microkernel-ER-versus-star},
if there is another $\chi^\prime: \cV \rightarrow \cW$ for some $\cW \hookrightarrow \P^* Z$, then the commuting diagrams for $\chi$ and $\chi^\prime$ composes to the same diagram for $\chi^\prime \circ \chi$.
\end{corollary}
\begin{proof}
Because $\Omega_Y^{-1} \otimes_{\Omega_{T^* Y}} \Omega_Y = \cO_{T^* Y}$ canonically, the composition of the diagram of $\chi$ and $\chi^\prime$ can thus be identified as the diagram

$$
\begin{tikzpicture}
\node at (0,2.5) {$ (\cE_X^\R -\sMod)^{w_2(X)}|_\cU$};
\node at (6,2.5) {$(\cE_Y^\R - \sMod)^{w_2(Y)} |_\cV$};
\node at (12,2.5) {$(\cE_Z^\R - \sMod)^{w_2(Z)} |_\cW$};
\node at (-0.2,0) {$ (a^*_X \cE_X^{\R,op} -\sMod)^{w_2(X)}|_{\cU}$};
\node at (6,0) {$ (a^*_Y \cE_Y^{\R,op} - \sMod)^{w_2(Y)} |_{\cV}$};
\node at (12,0) {$ (a^*_Z \cE_Z^{\R,op} - \sMod)^{w_2(Z)} |_{\cW}$};

\draw [->, thick] (1.8,2.5) -- (4.2,2.5) node [midway, above] {$\cH(\chi) \otimes_{\cE_X^\R} (-) $};
\draw [->, thick] (7.8,2.5) -- (10.2,2.5) node [midway, above] {$\cH(\chi^\prime) \otimes_{\cE_Y^\R} (-) $};
\draw [->, thick] (2,0) -- (3.8,0) node [midway, above] {$ $};
\draw [->, thick] (8.2,0) -- (9.8,0) node [midway, above] {$ $};

\draw [->, thick] (0,2.2) -- (0,0.3) node [midway, right] {$\Omega_X \otimes_{\cO_{T^* X}} (-)$}; 
\draw [->, thick] (12,2.2) -- (12,0.3) node [midway, right] {$\Omega_Z \otimes_{\cO_{T^* Y}} (-)$};

\node[scale=1.5] at (6,1.4) {$\circlearrowleft$};

\end{tikzpicture}
$$
where the bottom arrows are similarly given by the `op' version of the top arrows. But similar to the proof of \Cref{prop: qct-by-microkernel-ER}, there is a morphism
$
\cH(\chi^\prime) \otimes_{\cE_Y^\R} \cH(\chi) \rightarrow \cH(\chi^\prime \circ \chi) 
$
given by the morphism
\begin{align*}
&{p^a_{13}}_!({p^a_{23}}^{*} \mhom(\cK(\chi^\prime), \Omega_Y \boxtimes_\cO \cO_Z) \otimes {p^a_{12}}^{*}\mhom(\cK(h), \Omega_X \boxtimes_\cO \cO_Y))  \\
&\rightarrow \mhom\left(\cK(\chi^\prime) \circ \cK(\chi), (\Omega_Y \boxtimes_\cO \cO_Z) \circ (\Omega_X \boxtimes_\cO \cO_Y) \right) \\
&\rightarrow \mhom\left(\cK(\chi^\prime) \circ \cK(\chi), \Omega_X \boxtimes_\cO \cO_Z \right) = \mhom\left(\cK(\chi^\prime \circ \chi), \Omega_X \boxtimes_\cO \cO_Z \right)
\end{align*}
where we use the \eqref{canonical-kernel-composition} for the last equality. But \cite[Proposition 11.4.7]{kashiwara-schapira} then implies that this map is an equivalences. A similar computation holds for the `op' version as well.
\end{proof}

\subsection{$\mathcal{E}$-modules and microlocal Riemann-Hilbert} \label{sec: micro-diff}

For the purpose of microlocal Riemann-Hilbert, the sheaf $\cE_X^\R$ is too large.  In particular, this sheaf is only real conic, i.e., conic with respect to $\R_+$, but not complex conic, i.e., conic with respect to $\C^\times$. This is reflected in the fact that $\cE_X^\R$ can contain functions of the form $\log(\xi)$.

A naive way to  enforce respect for $\C^\times$ 
is the following. 
Let 
$\gamma: \dT^* X \rightarrow \P^* X$ be the quotient map by $\C^\times$, and consider 
$$\cE_X^\infty := \gamma^* \gamma_* \cE_X^{\R}.$$
Sections of this were called ``microdifferential operators of infinite order'' in \cite[8.2.15]{bjork-analytic-D}, and ``pseudo-differential operators'' in \cite{sato-kashiwara-kawai}. 

\begin{remark}
As mentioned at the beginning of \cite[2.1]{sato-kashiwara-kawai} or \cite[Remark 11.4.5]{kashiwara-schapira}, the sheaf of (infinite order) differential operator can be obtained by
$$ \cD_X^\infty = \Gamma_{\Delta_X}( \Omega_X \boxtimes_{\cO} \cO_X)[-n].$$
\end{remark}

In fact, certain important `finite degree' properties are lost in $\cE_X^\infty$.  This is a microlocal version of the following issue in the theory of $\cD$-modules: 
 for a closed analytic set $Z \subseteq X$, there is an important difference between sections with support in $Z$ and sections with `temperate' support in $Z$  (see e.g. \cite[II.5]{bjork-analytic-D}).

There is a better version, denoted by $\cE_X$ \footnote{Originally denoted by $\cP_X^f$ in \cite[Definition1.5.6]{sato-kashiwara-kawai}}, whose sections are called (finite-order) microdifferential operators. In  \cite{kashiwara-microlocal-D, bjork-analytic-D,schapira-microdifferential}, this sheaf of rings is introduced 
by taking formal series of differential operators, adjoining `$\partial_x^{-1}$', to obtain what is commonly denoted as $\hat{\cE}_X$\footnote{This is the ring of `formal' microdifferential operators, which \cite{kashiwara-microlocal-D} mostly works with.}, and imposing convergence conditions. More convenient to us is a variant of the original approach to $\cE_X$ in \cite{sato-kashiwara-kawai},
developed further by Andronikof \cite{andronikof-microlocalisation}, 
where the ring structure is functorially extracted from $\mhom$.
Some discussions from a modern point of view can be found in \cite{guillermou-dg-microlocalization, prelli}.

The key ingredient is the tempered microlocal hom \cite[Def.\ 2.3.1]{andronikof-microlocalisation}:  
$$\tmhom(-,\cO_X) : \sh_{\mathbb{R},c}(X) \to \sh(T^*X),$$
where $sh_{\mathbb{R},c}(X)$ is the category of $\R$-constructible sheaves.

The construction of $\cE_X$ is completely parallel to that of $\cE^\R$ in the previous subsection, except that all occurrences of $\mhom$ are replaced with $\tmhom$. Namely, one first defines
$$ \cE_X^{\R,f} \coloneqq \tmhom(\C_{\Delta_X}[-1], \Omega_X \boxtimes_{\cO} \cO_X) \in  \sh(T^* X)^\heartsuit.$$
Then one sets
$$\cE_X := \gamma^* \gamma_* \cE_X^{\R,f} \hookrightarrow \cE_X^{\R,f},$$
where we recall that $\gamma: \dT^* X \rightarrow \P^* X$ is the quotient map. 
Andronikof checks that this agrees with other definitions \cite[Cor.\ 5.5.2, Cor.\ 5.6.1]{andronikof-microlocalisation}.

Since $\tmhom$ and $\mhom$ enjoy the same formal properties, the algebras $\cE_X^{\R,f}$ and $\cE_X$ correspondingly satisfy all the properties enjoyed by $\cE^\R$ which were discussed in the previous subsection (\cite[Cor.\ 5.5.2, Cor.\ 5.6.1]{andronikof-microlocalisation}). For example, 
$H^j \tmhom(F,\cO_X)$ (resp. $H^j \gamma^* \gamma_* \tmhom(F,\cO_X)$ ) is a left $\cE_X^{\R,f}$ (resp. $\cE_X$) module.
Similarly, quantized contact transforms exist both for $\cE^{\R,f}$-modules and $\cE$-modules.

The  works of Andronikov and Waschkies \cite{andronikof-microlocal-RH, waschkies-microlocal-RH} microlocalized the usual Riemann-Hilbert correspondence to projectivized cotangent bundles:
\begin{theorem}[Local microlocal Riemann--Hilbert; Andronikov \cite{andronikof-microlocal-RH} Waschkies \cite{waschkies-microlocal-RH}] \label{thm: local mrh}
Let $X$ be a complex manifold. There is an equivalence 
$$\mRH_X: \P erv_{\P^* X} \xrightarrow{\sim} \cE_X -\sMod_{rh}$$
between perverse microsheaves on $\P^* X$ and regular holonomic $\cE_X$-modules. 
Furthermore, the equivalence respects contact transform \cite[Corollary 3.4.3]{waschkies-microlocal-RH}. That is, in the setting of \Cref{sec: mRH-for-kernel}, we have the following commuting diagram:
$$
\begin{tikzpicture}
\node at (0,2.5) {$ \P erv_{\cU}^{w_2(X)} $};
\node at (8,2.5) {$ \P erv_{\cV}^{w_2(Y)} $};
\node at (-0.2,0) {$ (\cE_X^\R -\sMod_{rh})^{w_2(X)}|_\cU$};
\node at (8,0) {$ (\cE_Y^\R - \sMod_{rh})^{w_2(Y)} |_\cV$};

\draw [->, thick] (1, 2.5) -- (7, 2.5) node [midway, above] {$\cK(\chi) \circ (-)$};
\draw [->, thick] (1.8,0) -- (6.1,0) node [midway, above] {$\cH(\chi) \otimes_{\cE_X^\R} (-)$};

\draw [->, thick] (0,2.1) -- (0,0.3) node [midway, right] {$\mRH_X |_\cU$}; 
\draw [->, thick] (8,2.1) -- (8,0.3) node [midway, right] {$\mRH_Y |_\cV$};

\node[scale=1.5] at (4,1.4) {$\circlearrowleft$};

\end{tikzpicture}
$$
\end{theorem}

Lastly, we mention that the algebra $\cE_X$ has the nice Morita-theoretical property of being Picard good. See \Cref{def: picard-good} for the definition.

\begin{theorem}[{\cite[Theorem 4.3.6]{dagnolo-polesello}}] \label{thm: E-is-picard-good}
The $\C$-algebra $\cE_X$ is Picard good.
\end{theorem}

\subsection{The microlocal Riemann--Hilbert correspondence for complex contact manifolds}

We recall Kashiwara's quantization of contact manifolds and Polsello's uniqueness criterion. 

\begin{theorem}[{\cite[Theorem 2]{kashiwara-quantization-contact}, \cite[Theorem 3.3]{polesello-uniqueness}}] \label{uniqueness-of-kashiwara's}
Let $(V,\cH)$ be a complex contact manifold. There exists a unique $\C$-algebroid $\cE_V$ admitting the following structures
\begin{itemize}
	\item[(i)] $\cE_V$ is filtered. (This means that for any $U \subset V$ and an objects $K, L \in \cE_V(U)$, the $\C$-module $Hom_{\cE_V(U)}(K,L)$ is filtered, the filtrations are compatible with composition, etc. see \cite[Sec.\ 2]{polesello-uniqueness}). 
	\item[(ii)] there is an isomorphism of commutative algebroid stacks $$\sigma: Gr(\cE_V) \simeq \left(\bigoplus_{m \in \Z} \cL^{m} \right)^+$$
where $\cL$ is the line bundle $(TV/\cH)^\vee$, and its zeroth degree associated graded is given by taking dual line bundle $(-)^\vee$, and on each local Darboux chart $(\cE_V,*)$ is equivalent to the canonical filtered algebroid $\cE_Y$ on $(\P^* Y,*)$.
	\item[(iii)] $\cE_V$ is equipped with an anti-involution $*$ (by definition: an equivalence of stacks $*: \cE_V \to \cE_V^{op}$ along with an invertible natural transformation $\epsilon: *^2 \Rightarrow \operatorname{id}$ such that $\epsilon \circ id_*: *^3 \Rightarrow *$ and $id_* \circ \epsilon: * \Rightarrow \*^3$ are inverse). 
\end{itemize}
These structures must be compatible; this means that $*$ is a filtered functor, $\epsilon$ is a filtered natural transformation, and there exists an invertible natural transformation 
\begin{equation}\label{eq: compatibility-1} 
\delta_0: \sigma_0 \circ Gr_0(*) \Rightarrow D \circ \sigma_0
\end{equation}
such that the following diagram commutes:
\begin{equation}\label{eq: compatibility-2} 
\begin{tikzcd}
        \sigma_0 \circ Gr_0(*^2) \ar[rr, Rightarrow, "\delta_0 \circ id_{Gr_0(*)}"] \ar[d, Rightarrow, "id_{\sigma_0} \circ Gr_0(\epsilon)"]  && D \circ \sigma_0 \circ Gr_0(*) \\
    \sigma_0 \ar[rr, Rightarrow, "\simeq"] && D^2 \circ \sigma_0^2 \ar[u, Rightarrow, "id_D \circ \delta_0"]
\end{tikzcd}
\end{equation}
Here $D: Pic_V \to Pic_V^{op}$ is the functor sending a line bundle to its dual.
\end{theorem}

We note that the notion of holonomic and regular holonomic is well-defined microlocally 
(\cite[Chapter 8]{kashiwara-microlocal-D}, \cite[Chapter VIII]{bjork-analytic-D}) 
and the work of Andronikov \cite[Chapter 5]{andronikof-microlocalisation} and Waschkies \cite[Section 3]{waschkies-microlocal-RH} shows that they are invariant under quantized contact transform.
Similarly on the microsheaf side, we show in our previous paper \cite[Section 6]{CKNS}, the notion of perverse $t$-structure exists on the canonical microsheaf category $\msh_V$ and, on Darboux charts, it coincides with the usual notion. 

\begin{definition}
We denote by $\P\msh_{V;\C-c}^b$ the subsheaf of $\P\mu sh_V$ consisting of complex constructible microsheaves with perfect microstalks, 
$\P erv_V  \subseteq \P\msh_{V;\C-c}^b$ the subsheaf (of abelian categories) given by the heart of the perverse $t$-structure. 
Similarly, we denote by $\cE_V -\sMod_h$ (resp. $\cE_V -\sMod_{rh}$) the subsheaf of $\cE_V -\sMod$ consisting of holonomic (resp. regular holonomic) $\cE_V$-modules.
\end{definition}

We turn now to our main task, of globalizing the Andronikov-Waschkies results to a comparison of $\P erv_V$ and $\cE_V -\sMod_{rh}$ on an arbitrary complex contact manifold $V$. 

Let $\cU \hookrightarrow \P^* X$, $\cV \hookrightarrow \P^* Y$ be open subsets and $\chi: \cU \to \cV$ be a complex contactomorphism. Then we define the following ``tempered'' analog of \eqref{equation:hchi-big}: \begin{equation}\label{equation:hchi-small}
\cH(\chi) \coloneqq \gamma^* \gamma_* \tmhom( \cK(\tilde{\chi})[-n], \Omega_X \boxtimes_\cO \cO_Y) \in \sh^{w_2(X \times Y)}(\cU \times \cV)^\heartsuit.
\end{equation}

\begin{proposition} \label{prop: qct-by-microkerne-E}
Tensoring with the bimodule $\cH(\chi)$ (resp. $\cH(\chi^{-1})^{op}$ ) induces an equivalence
\begin{align*}
\cH(\chi) \otimes_{\cE_X } (-)&: (\cE_X -\sMod)^{w_2(X)}|_\cU \xrightarrow{\sim}  (\cE_Y - \sMod)^{w_2(Y)} |_\cV \\
(resp. \, \cH(\chi^{-1})^{op} \otimes_{\cE_X^{op} } (-)&: (\cE_X^{op} -\sMod)^{w_2(X)}|_\cU \xrightarrow{\sim}  (\cE_Y^{op} - \sMod)^{w_2(Y)} |_\cV. )    
\end{align*}
 
Furthermore, the obvious tempered analog of \Cref{prop: qct-by-microkernel-ER-versus-star} and \Cref{cor: composition-qct-by-microkernel-ER-versus-star} hold.
\end{proposition}

\begin{proof}
As in \Cref{prop: qct-by-microkernel-ER} and \Cref{prop: qct-by-microkernel-ER-versus-star}, replacing everywhere when appropriate $\mhom$ by $\tmhom$.
\end{proof}

This proposition gives us (what we will eventually show to be) an alternative construction of Kashiwara's algebroid $\cE_V$:
Choose a Darboux cover $\{(\cU_\alpha, f_\alpha)\}$ as in  \Cref{thm: gluing-description-microsheaf}; this provides gluing data for $\P erv_V$.  Apply the construction of \eqref{equation:hchi-small} to said gluing data to obtain the $\mathcal{H}$ and $c$ of 
\begin{equation} \label{for: gluing-description-E-module}
\left( \{ (\cE_\alpha - \sMod)^{w_2(\alpha)} \}, \{\cH_{\beta \alpha}\}, \{ c_{\gamma \beta \alpha} \} \right)
\end{equation}
We have applied a functor to gluing data, hence obtain gluing data for some sheaf of categories which is locally isomorphic to $\cE_V - \sMod_{rh}$.  We denote the resulting sheaf of categories by $\cC_V$.  

\begin{definition}[{\cite[Definition 5.2.1]{dagnolo-polesello}}]
Let $V$ be a contact manifold. 
\begin{enumerate}
\item
An \emph{$\cE$-algebroid} $\cA$ on $V$ is a $\C$-algebroid with the property that for any Darboux ball $V \supset \cU \xhookrightarrow{f} \P^* X$, there is an equivalence $\cA |_\cU \cong (f^* \cE_X)^+$.

\item A sheaf of $\C$-linear categories $\cC$ on $V$ is called a \emph{sheaf of twisted $\cE$-modules} if for any Darboux ball $V \supset \cU \xhookrightarrow{f} \P^* X$, there is an equivalence $\cC |_\cU \cong f^*  (\cE_X -mod)$.
\end{enumerate}
\end{definition}

The Picard good property of $\cE_X$, recalled in \Cref{thm: E-is-picard-good}, implies that there is a simple description of sheaves of twisted $\cE$-modules.
\begin{theorem}[{\cite[Theorem 5.2.3]{dagnolo-polesello}}] \label{thm: classification-twisted-E-modules}
\
\begin{enumerate}
\item Any sheaf of twisted $\cE$-modules $\cC$ is equivalent to $\cA -mod$ for some $\cE$-algebra.
\item Two $\cE$-algebroids $\cA$ and $\cB$ are equivalent if and only if the associated sheaves of twisted $\cE$-modules $\cA -mod$ and $\cB-mod$ are equivalent.
\end{enumerate}
\end{theorem}

\begin{proof}
This is a direct consequence of \Cref{thm: E-is-picard-good} and \Cref{prop: algebroid-eq-=-module-eq}.
\end{proof}

The above \Cref{thm: classification-twisted-E-modules} implies that $\cC_V$ is of the form $\cE_V^\prime -mod$ for some $\cE$-algebroid $\cE_V^\prime$.

\begin{theorem} \label{thm: altenative-construction-of-kashawara-E}
The algebroid $\cE_V^\prime$ is equivalent to the canonical algebroid $\cE_V$. 
\end{theorem}

\begin{proof}
We will eventually use Polsello's criterion, for which we will need a  $*: \cE_V^\prime \xrightarrow{\sim} \cE_V^{\prime, op}.$  Let us construct it.

First note that appealing to the statement about right modules in \Cref{prop: qct-by-microkerne-E}, the same choice of Darboux charts $\{ (\cU_\alpha, f_\alpha)\}$, gives us the gluing data
$$
\left( \{ (a_\alpha^* \cE_\alpha^{op} - \sMod)^{w_2(\alpha)} \}, \{\cH_{\alpha \beta}^{op}\}, \{ c_{\gamma \beta \alpha} \} \right)
$$
where $\cH_{\alpha \beta}^{op} = \cH(f_\alpha \circ f_\beta^{-1})^{op}$ and, by the proof of \Cref{thm: classification-twisted-E-modules}, it glues to $\cE_V^{\prime, op} -mod$.
Furthermore, the tempered version of \Cref{prop: qct-by-microkernel-ER-versus-star} and \Cref{cor: composition-qct-by-microkernel-ER-versus-star} implies that there exists an equivalence  
$$*: \cE_V^\prime \xrightarrow{\sim} \cE_V^{\prime, op},$$
which can be recovered from its corresponding twisted $\cE$-module equivalence, glued from 
$$\Omega_X \otimes_{\cO_{\dT^* X}} (-): \cE_X -\sMod^{w_2(X)} \xrightarrow{\sim} a^* \cE_X^{op} - \sMod^{w_2(X)}$$ 
on Darboux charts. Indeed, in order to glue, we have to check that the diagram of \Cref{prop: qct-by-microkernel-ER-versus-star}  is functorial with respect to horizontal compositions; 
since we are gluing morphisms between $(2,1)$-sheaves, we need to check only that, for triple intersections indexed by $\alpha$, $\beta$, and $\gamma$, the commuting diagram compose in the way illustrated below:
$$
\begin{tikzpicture}
\node at (0,2) {$\alpha$};
\node at (2.5,2) {$\beta$};
\node at (0,0) {$\alpha$};
\node at (2.5,0) {$\beta$};
\node[scale=1.5] at (1.25,1.1) {$\circlearrowleft$};

\draw [->, thick] (0.4,2) -- (2.2,2) node [midway, above] {$\sim$};
\draw [->, thick] (0.4,0) -- (2.2,0) node [midway, above] {$\sim$};

\draw [->, thick] (0,1.7) -- (0,0.3) node [midway, right] {$ $}; 
\draw [->, thick] (2.5,1.7) -- (2.5,0.3) node [midway, right] {$ $};

\node at (3,1) {$\bigcirc$};

\node at (3.5,2) {$\beta$};
\node at (6,2) {$\gamma$};
\node at (3.5,0) {$\beta$};
\node at (6,0) {$\gamma$};
\node[scale=1.5] at (4.75,1.1) {$\circlearrowleft$};

\draw [->, thick] (3.9,2) -- (5.7,2) node [midway, above] {$\sim$};
\draw [->, thick] (3.9,0) -- (5.7,0) node [midway, above] {$\sim$};

\draw [->, thick] (3.5,1.7) -- (3.5,0.3) node [midway, right] {$ $}; 
\draw [->, thick] (6,1.7) -- (6,0.3) node [midway, right] {$ $};

\node at (6.5,1) {$=$};

\node at (7,2) {$\alpha$};
\node at (9.5,2) {$\gamma$};
\node at (7,0) {$\alpha$};
\node at (9.5,0) {$\gamma$};
\node[scale=1.5] at (8.25,1.1) {$\circlearrowleft$};

\draw [->, thick] (7.4,2) -- (9.2,2) node [midway, above] {$\sim$};
\draw [->, thick] (7.4,0) -- (9.2,0) node [midway, above] {$\sim$};

\draw [->, thick] (7,1.7) -- (7,0.3) node [midway, right] {$ $}; 
\draw [->, thick] (9.5,1.7) -- (9.5,0.3) node [midway, right] {$ $};

\end{tikzpicture}
$$ 
But this is \Cref{cor: composition-qct-by-microkernel-ER-versus-star}.  

Up until now, our strategy is to define the relevant functors by first defining them as functors between module categories using bimodules, and turn them into algebroid morphisms by appealing to (the proof of) \Cref{thm: classification-twisted-E-modules}.
However, as explained in \Cref{prop: algebroid-eq-=-module-eq},  on small open sets where the algebroid admits a section so there exists a common generator for the invertible bimodules, the passage from bimodules to ring homomorphisms is given by \Cref{lem: bimodule-ring-morhpism}.
In fact, this process is exactly how the quantized contact transform, as ring isomorphisms, are classically obtained in for example \cite[Theorem 11.4.9]{kashiwara-schapira}. 

Thus, refining the Darboux cover if needed, we can assume the transition maps comes from genuine ring isomorphisms $\Phi(\chi): \cE_X |_\cU \xrightarrow{\sim} \cE_Y^{op} |_\cV$.
By \cite[Theorem 7.2.2]{kashiwara-introduction}, this isomorphism is filtered and compatible with the symbols, i.e., there is a commuting diagram

$$
\begin{tikzpicture}
\node at (0,2) {$\cE_X(m)$};
\node at (5,2) {$\chi^* \cE_Y^{op}(m)$};
\node at (0,0) {$\cO_{T^* X}(m)$};
\node at (5,0) {$\chi^* \cO_{T^* Y}^{op}(m)$};

\draw [->, thick] (0.7,2) -- (4,2) node [midway, above] {$\Phi(\chi)$};
\draw [->, thick] (1,0) -- (3.7,0) node [midway, above] {$(\chi^{-1})^*$};

\draw [->, thick] (0,1.7) -- (0,0.3) node [midway, right] {$\sigma_m$}; 
\draw [->, thick] (5,1.7) -- (5,0.3) node [midway, right] {$\sigma_m$};
\end{tikzpicture}
$$ 
Now (i) is automatic since the filtration is defined locally and is preserved by $\Phi(\chi)$.

For (ii), the above commuting diagram implies that the identification $\sigma_X: Gr(\cE_X) \xrightarrow{\sim}  \cS_{\cO_X}(\Theta_X)$ can be glued to a global identification $\sigma: Gr(\cE_V^\prime) \simeq \left(\bigoplus_{m \in \Z} \cL^{m} \right)^+$. In fact, the object $Gr(\cE_V^\prime)$ is well-defined as an algebra (not merely as an algebroid) and $\sigma$ comes from an algebra isomorphism $\sigma: Gr(\cE_V^\prime) \simeq  \bigoplus_{m \in \Z} \cL^{m}$. This is because, on any triple overlap indexed by $\alpha$, $\beta$, $\gamma$, the composition $\Phi_{\alpha \gamma} \Phi_{\gamma \beta} \Phi_{\beta \alpha}$ is a filtered and symbol-preserving ring isomorphism on $\cE_\alpha$ and so it becomes the identity when taking $Gr$. (In particular, the ambiguity $P_{ijk}$ in  \cite[(0.1)]{kashiwara-quantization-contact}, which leads to the necessity of working with algebr\emph{oids}, disappear when passing to graded algebras.)

Passing now to (iii): on local charts, we have natural transformation of bimodules $$ \Omega_X^{-1} \otimes_\cO \Omega_X \rightarrow \cO_X$$ which clearly glue to a natural transformation of bimodules on $V$. This in turn induces the natural tarnsformation $*^2 \to \operatorname{id}$ by \Cref{prop: algebroid-eq-=-module-eq}. The fact that this natural transformation is a natural isomorphism can be checked on Darboux charts, where it is obvious.  

Finally, we exhibit the natural transformation $\delta_0$: in a Darboux chart $V \supset \cU \hookrightarrow P^*X$, $\delta_0$ is just given by the diagram
\begin{equation}\label{equation:polesello-compatibility}
\begin{tikzcd}
Gr_0(\cE_X)^+ \ar[r, "Gr_0(*)"] \ar[d, "\sigma_0"] & Gr_0(\cE_X^{op})^{+} \ar[d, "\sigma_0"]  \\
Pic_{\mathbb{P}^*X} \ar[r, "D"]& Pic_{\mathbb{P}^*X}^{op}
\end{tikzcd}
\end{equation}
That this diagram commutes is a restatement of \Cref{lemma:pic-O}. It follows from (ii) that these locally defined natural transformations induce a globally defined natural transformation. (We remark that a more intrinsic way to see this is to use \Cref{lem: compatibility-with-star} and trade $(\cE_X -mod)^{w_2(X)}$ for $\cE_X^{\sqrt{\Omega_X}} -mod$; under this identification, $Gr_0(*)$ is exactly $D$.)
This completes the proof. 
\end{proof}

\begin{theorem}[Global microlocal Riemann--Hilbert correspondence]\label{thm:gloabl mrh}
Let $V$ be a complex contact manifold. There is an equivalence
$$\mRH_V: \P erv_V \xrightarrow{\sim} \cE_V - \sMod_{rh}$$
between perverse microsheaves on $V$ and regular holonomic $\cE_V$-modules.
\end{theorem}

\begin{proof}
The sheaves $\Perv_V \subseteq \msh_V$ and $\cE_V- \sMod_{rh} \subseteq \cE_V -\sMod$ can be reconstructed by the gluing data
$$\left( \{ \P erv_\alpha^{w_2}\}, \{\cK_{\beta \alpha}\}, \{ c_{\gamma \beta \alpha} \} \right) \, \text{and} \, \left( \{ (\cE_\alpha - \sMod_{rh})^{w_2(\alpha)} \}, \{\cH_{\beta \alpha}\}, \{ c_{\gamma \beta \alpha} \} \right)$$
where we use \Cref{thm: altenative-construction-of-kashawara-E} to conclude for the latter.
But the compatibility statement in \Cref{thm: local mrh} implies that the local equivalence $\mRH_\alpha$ glues to an equivalence on $V$.
\end{proof}

\begin{remark} \label{not-identity}
To get a feel for Polesello's criterion \Cref{uniqueness-of-kashiwara's}, it is instructive to consider the example of a complex manifold $X$. Let $L= \Omega_X$ be the canonical bundle. Then Kashiwara's stack is $(\cD_X^{L/2})^+$ \emph{which is generally different from} $\cD_X^+$. So let us understand where Polesello's criterion fails for $\cD_X^+$.

Although $\cD_X \neq \cD_X^L$ in general, combining \Cref{lem: equivalence-modules} and Morita theory in the sheaf setting, e.g. \cite[Corollary 3.3.8]{dagnolo-polesello}, we see that the \emph{algebroids} $\cD_X^+$ and $\cD_X^{L,+}$ are equivalent.  (The bimodule is $L \otimes_{\cO} D_X$). As explained in \Cref{star-involution-D-modules}, $*$ gives an algebroid isomorphism $D_X^+-\mathfrak{M}od \to D_X^{L, op}-\mathfrak{M}od$. And by Morita theory (applied to the sheaf of categories $D_X^+-mod$), we know that  $D_X^+-mod$ and $(D_X^L)^+-mod$ are isomorphic. This allows us to construct an anti-involution on $D_X^+-\mathfrak{M}od$, which verifies (iii) in \Cref{uniqueness-of-kashiwara's}.

Consider further the fact that the associated graded $Gr(\cD_X) = \cS_\cO(\Theta)$ is given by symmetric products of holomorphic vector fields, i.e., functions of the cotangent bundle which are polynomials in the fiber direction. In particular, $Gr_0$ is simply $\cO_X$. Because $\cS_\cO(\Theta)$ is commutative, $Gr(\cD_X) = Gr(\cD_X^L)$ as algebras. This verifies (i) and (ii) in \Cref{uniqueness-of-kashiwara's}.

What breaks is the compatibility condition \eqref{eq: compatibility-1} . Namely, the equivalence $\cD_X^+ = \cD_X^{L,+}$ after applying $Gr_0$ is not the identity, since on the level of modules, it is given by
\begin{align*}
\cO_X -mod &\xrightarrow{\sim} \cO_X -mod \\
\cM &\mapsto L \otimes_\cO \cM.
\end{align*}
\end{remark}

\subsection{Contactomorphism group actions}\label{subsection:cont-actions}
In this subsection, we explain that the group of  complex contactomorphisms provides both $\cE_V$ and $p_* \msh_{V_0}$ a canonical equivariant structure. As a consequence, their global sections admit an action by contactomorphisms. Furthermore, we show that when restricting to the corresponding subcategories, the microlocal Riemann--Hilbert correspondence is compatible with contactomorphism group actions. 

Let $V$ be a complex contact manifold. Denote by $\Cont(V)$ the group of complex contactomorphisms on $V$ and let $u: \Cont(V) \times V \rightarrow V$ be the tautological action $u(g,p) \coloneqq g(p)$. As usual, we denote the symplectization of $V$ by $\tilde{V}$ and write $V_0:= \tilde{V}/\mathbb{R}_+$.

\begin{proposition}\label{prop: contact-action-mu}
The sheaf $p_* \msh_{V_0}$ admits a canonical $\Cont(V)$-equivariant structure. In particular, $\msh_{V_0}(V_0)$ admits a $\Cont(V)$-action. Furthermore, the $\Cont(V)$-equivariant structure restricts to the subsheaves, $p_* \msh_{V_0; \C-c}$, $\Pmsh_V$, and $\mathbb{P}erv_V$, so their global sections are preserved under the $\Cont(V)$-action.
\end{proposition}

\begin{remark} \label{rmk: equivraint-2-sheaf-description}
Similarly to the situation discussed in \Cref{rmk: equivraint-sheaf-description}, an equivariant sheaf with coefficients in $\Ab$ is given by a triple $(F, \phi, T)$ where, for $g_1, g_2 \in G$, there is an invertible $2$-morphism
$$T_{g_2 g_1}: ( u_{g_1}^* u_{g_2}^* F \xrightarrow{u_{g_1}^* \phi_{g_2} } u_{g_1}^* F \xrightarrow{\phi_{g_1}} F) = ( u_{g_2}^* u_{g_1}^* F =  u_{g_2 g_1}^* F \xrightarrow{ \phi_{g_2 g_1} } F  ),$$
and these satisfy the usual \v{C}ech ($3$-)cocycle conditions.
This an additional layer of structure is due to the existence of non-trivial $2$-morphisms in $\Ab$.

We also remark that, for an equivariant sheaf $(F, \phi, T)$, the global sections of $F$ admits a $G$-action given by
$$\Gamma(X;F) = \Gamma(X; u_g^* F)  \xrightarrow{\phi_g} \Gamma(X; F),$$
and the $T$'s provides the natural transformation which witness associativity etc.
\end{remark} 

\begin{proof}
As remarked after \Cref{def: pmsh}, the sheaf $\msh_{V_0}$ is the pullback of a universal sheaf $\msh_{B Pic(\cC)_0} (\xi)$ on $B Pic(\cC)_0(\xi)$ along the section $V_0 \rightarrow V_0^{B Pic(\cC)_0}$ which corresponds to the canonical orientation $o_{can}$. By \cite[Lemma 5.6]{CKNS}, for any $g \in \Cont(V)$, we have a canonical equivalence
$\phi_{\mu, g}: u_g^* \msh_{V_0}  =  \msh_{V_0}$. Indeed, this is one of the ingredient we needed to argue for the sheaf kernel $\cK_{\beta \alpha}$ in \Cref{thm: gluing-description-microsheaf}, and the existence of the $2$-morphism $T_{g_2 g_1}$, for $g_1, g_2 \in \Cont(V)$, and the fact that they satisfy the \v{C}ech cocycle conditions follow from the same argument. Lastly, the conditions of being in the subcategories mentioned above are invariant under the canonical equivalence $\phi_{\mu,g}$, and the induced $\Cont(V)$ on global sections follows from \Cref{rmk: equivraint-2-sheaf-description}.
\end{proof}

\begin{proposition} \label{prop: contact-action-E}
The sheaf of categories $\cE_V - \sMod$ admits a $\Cont(V)$-equivariant structure. Furthermore, the equivariant structure restricts to the subsheaves consisting of coherent, holonomic, and regular holonomic objects. Thus, the $\Cont(V)$-action on its global section $\cE_V -mod$ preserves, $\cE_V -mod_{coh}$, $\cE_V -mod_{hol}$, and $\cE_V -mod_{rh}$.
\end{proposition}

\begin{proof}
Kashiwara constructs the algebroid $\cE_V$ by choosing a cover ${U_\alpha}$ by Darboux charts, and identify the algebras $\cE_\alpha$ on double over laps with algebra isomorphisms $f_{\beta \alpha}$, which have to be glued as algebroids using the cocycles $P_{\gamma \beta \alpha}$ given by microdifferential operator \cite[(0.1)]{kashiwara-quantization-contact}. In other words, $\cE_V$ is glued from the data $(\cE_\alpha, f_{\alpha \beta}, P_{\gamma \beta \alpha})$. However given $g \in \Cont(V)$, for each $\alpha$, there is a canonical morphism by pullback
\begin{align*}
u_g ^* \cE_\alpha-mod &\rightarrow \cE_\alpha-mod \\
\cM_\alpha &\mapsto g_* \cM
\end{align*}
where we abuse the notation and omit the needed restriction of $g$ to the covering opens. Thus, there is an equivalence $\phi_{\cE,g}: g^* \cE_V = \cE_V$ where $g^* \cE_V$ is  the algebroid glued from $(g^* \cE_\alpha, g^* f_{\alpha \beta}, g^* P_{\gamma \beta \alpha})$. However, $g^* \cE_V$ is just $\cE_V$ since it satisfies the same uniqueness criterion \Cref{uniqueness-of-kashiwara's}. As in the previous \Cref{prop: contact-action-mu},  $\phi_{\cE,g}$ composes naturally on $g \in \Cont(V)$.
\end{proof}

\begin{proposition} \label{prop: compatibility-of-actions}
The microlocal Riemann-Hilbert correspondence (\Cref{thm:gloabl mrh}) is compatible with the $\Cont(V)$-equivariant structure of $\P erv_V$ and $\cE_V - \sMod_{rh}$. 
\end{proposition}

\begin{proof}
Recall that the key observation in the proof of \Cref{thm:gloabl mrh} is that gluing data for $p_* \msh_{V_0}$ provides a gluing data for $\cE_V$ as explained in \Cref{thm: gluing-description-microsheaf} and (\ref{for: gluing-description-E-module}). However, Maslov-to-kernel (\Cref{thm: maslov-to-kernel}) implies that the equivalence $\phi_{\mu,g}$, for $g \in \Cont(V)$, in \Cref{prop: contact-action-mu} is realized concretely by moving around the gluing data from \Cref{thm: gluing-description-microsheaf},  by $g$, in the same fashion as for the E-module case in the previous \Cref{prop: contact-action-E}. But \Cref{thm:gloabl mrh} is exactly obtained by matching these two sets of gluing data and we thus obtain the commuting diagram:

$$
\begin{tikzpicture}
\node at (0,2) {$\P erv_V $};
\node at (6,2) {$\cE_V - \sMod_{rh}$};
\node at (0,0) {$\P erv_V $};
\node at (6,0) {$\cE_V - \sMod_{rh}$};

\draw [double equal sign distance, thick] (0.8,2) -- (4.7,2) node [midway, above] {$\mRH_V$};
\draw [double equal sign distance, thick] (0.8,0) -- (4.7,0) node [midway, above] {$\mRH_V$};

\draw [double equal sign distance, thick] (0,1.7) -- (0,0.3) node [midway, right] {$\phi_{\mu,g}$}; 
\draw [double equal sign distance, thick] (6,1.7) -- (6,0.3) node [midway, right] {$\phi_{\cE,g}$}; 
\end{tikzpicture}
$$ 

\end{proof}

We now upgrade \Cref{thm:gloabl mrh} to the $\C^\times$-equivariant setting. We refer to \Cref{apn: equivariance} for detailed discussions on equivariance and only recall here that, in order to upgrade \Cref{prop: compatibility-of-actions} to an equivalence between equivariant objects, specify how the identification $\phi_g$ depends on $g$: The dependence should be continuous on the microsheaf side. On the E-module side, similar to the D-module situation, requiring the dependence to be holomorphic gives the notion of quasi-equivariant modules. Such a module $\cM$ admit a Lie algebra action by $\C = \Lie(\C^\times)$ by differentiating $\phi_r: r^* \cM = \cM$ for $r \in \C^\times$. However, $\C$ admits another Lie algebra action coming $\Lie(\C^\times) \rightarrow \cE_V$ by differentiating the $\C^\times$-action on the space $V$. The subcategory of equivariant modules consists of quasi-equivariant modules on which these two Lie algebra actions agree. We will highlight in the proof when the distinction becomes relevant.   

\begin{remark}
We note that the discussion holds more generally for actions by complex Lie group but we will not need such a generality and the formula which we will use in the following proof needs to much more complicated for that.
\end{remark}

We will need the following lemma.

\begin{lemma} \label{lem: complex-movie}
The map $a: \mathbb{C}^\times \times V \rightarrow V$ is an action by contactomorphisms if and only if it can be lifted to an action $\tilde{a}: \mathbb{C}^\times \times \tilde{V} \rightarrow \tilde{V}$ on its symplectization by homogeneous symplectomorphisms. Furthermore, the movie of $\tilde{a}$ 
$$\Lambda_a =\{\left(x, -\alpha, r, -\lambda_{V, \tilde{a}_r(x,-\alpha)}\left(\partial_r \tilde{a}_r(x, \alpha) \right), \tilde{a}_r(x,\alpha)| (x,\alpha) \right) \in V, r \in \C^\times\}  \subseteq \tilde{V} \times T^* \C^\times \times \tilde{V}$$
is a complex conic Lagrangian. Here we denote by $\lambda_{V}$ the Liouville form on $\tilde{V}$.  
\end{lemma}

\begin{proof}
 This is the complex version of \cite[Lemma A.1]{guillermou-kashiwara-schapira} but we follow the sign convention used in \cite[Proposition 3.12]{kuo-wrapped}. 
\end{proof}

\begin{theorem} \label{thm: equi-gloabl-mrh}
Let $a: \C^\times \times V \rightarrow V$ be a $\C^\times$-action by contactomorphisms. There is an equivalence between categories
$$\mRH_V^{a} \left(\P erv_V(V)\right)^a \xrightarrow{\sim} (\cE_V - mod_{rh})^a$$
between $G$-equivariant perverse microsheaves on $V$ and $G$-equivariant regular holonomic $\cE_V$-modules.
\end{theorem}

\begin{proof}
We begin with recalling that the right hand side $(\cE_V - mod_{rh})^a$ can be computed by the limit of the three-term sequence
$$ \cE_V - mod_{rh}  \rightrightarrows \cD_{\C^\times} \boxtimes \cE_V - mod_{rh} 
\mathrel{\substack{\textstyle\rightarrow\\[-0.6ex]
\textstyle\rightarrow \\[-0.6ex]
\textstyle\rightarrow}}
\cD_{\C^\times \times \C^\times} \boxtimes \cE_V - mod_{rh}.
$$
  
Here the arrows are given by, for example, $E p_2^*$ and $E a^*$, and $E a^*$ is the functor  
\begin{align}\label{for: E-functoriality}
E f^*: \cE_V -mod &\rightarrow  (\cD_{\C^\times} \boxtimes \cE_V) -mod    \\
\cM &\mapsto u^{*,hol} \cM = \cO_{\C^\times \times \cO_V} \otimes_{a^{*,top} \cO_V} a^{*,top} \cM \notag. 
\end{align}
We note that we use the fact that $\cO_{\C^\times}$ has a $\cD_{\C^\times}$-module structure.
We remark that, by forgetting the $\cD_{\C^\times}$-module structure, we could replace the second term with $(\cO_{\C^\times} \boxtimes \cE_V) -mod$ and it would give us the notion of $\C^\times$-quasi-equivariant $\cE_V$-modules.

To get the corresponding sequence on the microsheaf side is less straightforward and requires sheaf quantization. We consider the contact movie from the previous \Cref{lem: complex-movie}. As we are in the complex setting, the canonical orientation, as discussed in \Cref{sec: msh-complex}, induces a microsheaf kernel
$$\chi(a) \in \P \msh(\tilde{V} \times T^* \C^\times \times \tilde{V}),$$
see e.g. (\ref{Kchi}). When restricting to Darboux charts, it can be identified with microsheaf kernel of contactomorphisms between coprojective bundles, and Maslov-to-Kernel \Cref{thm: maslov-to-kernel} implies that, when fixing $r \in \C^\times$, we recover the action from \Cref{prop: contact-action-mu}. On the other hand, convolution from (\ref{convolution-microsheaf-kernel}) extends naturally to the family version and this provides us the three term sequence 
$$ \Perv_V(V)  \rightrightarrows Perv(T^* \C^\times \times V) 
\mathrel{\substack{\textstyle\rightarrow\\[-0.6ex]
\textstyle\rightarrow \\[-0.6ex]
\textstyle\rightarrow}}
\Perv(T^* \C^\times \times T^* \C^\times \times V),
$$
where the first two arrows are given by $\chi(a)$ and $\chi(p_2)$, quantization of the trivial action.
But we can then apply \Cref{thm:gloabl mrh} term by term, since we see from \Cref{prop: compatibility-of-actions} that the action is compatible under microlocal Riemann-Hilbert, and obtain $$\mRH_V^a: \left(\P erv_V(V)\right)^a \xrightarrow{\sim} (\cE_V - \mod_{rh})^a$$ by taking the limit.
\end{proof}

\section{The canonical sheaf of WKB operators}  \label{sec: w-module} 

In this section, we show that a choice of a primitive $1$-form\footnote{In the sequel, given a complex symplectic manifold $(X, \omega)$, we call a primitive of $\omega$ a \emph{Liouville form} or \emph{Liouville structure.}} along with a $\C^\times$-action of non-zero integer weight on a complex symplectic manifold $\fX$ induces canonically a ``quantized Liouville structure".  More precisely, following Kashiwara and Rouquier's terminology \cite{kashiwara-rouquier}, we will exibit a natural `$F$-action' $F_\fX$ on the canonical WKB algebroid $\cW_\fX$. 
Ultimately, we will identify the category of $F_\fX$-equivariant regular holonomic $\cW_\fX$-modules with microsheaves on $\fX$ by applying \Cref{thm:gloabl mrh intro}.

\subsection{$W$-modules}

Let $X$ be a complex manifold. Denote by $\rho: \P^*(X \times \C) \rightarrow T^* X$ the conification map defined by $\rho(z,\xi,s,\sigma) = (z, \sigma^{-1}\xi)$ where $(z, s)$ are the base coordinates and $(\xi,\sigma)$ are the cotangent coordinates. The algebra of WKB-differential operators on $T^* X$, defined by Polesello and Schapira in \cite[Section 8]{polesello-schapira} is a sheaf defined by the subring 
$$\cE_{\P^*(X \times \C), t} \coloneqq \{P \in \cE_{\P^*(X \times \C)}| [P, \partial_t] = 0\},$$
of microdifferential operators which commute with $\partial_t$ by pushing along $\rho$, i.e., $$\cW_{T^* X} \coloneqq \rho_* \cE_{\P^*(X \times \C),t}.$$

Unlike $\cE_{T^* X}$, the ring $\cW_{T^* X}$ is non-homogeneous: A section $P \in \cW_{T^* X}(U)$ has the form 
\begin{equation} \label{for: section-of-W}
P = \sum_{l \geq -m} f_l(z,w) \hbar^l, \ f_l \in \cO(U), \ m \in \N
\end{equation}
so it is in particular an algebra over the Laurent series ring $\C[[\hbar, \hbar^{-1}]$; we recall that this is by setting $w = \tau^{-1} \xi$ and $\hbar = \tau$.

For any complex symplectic manifold $\fX$, Polesello and Schapira in \cite[Section 9]{polesello-schapira} show that, similarly to the case of complex contact manifold studied in \cite{kashiwara-quantization-contact}, there exists a \emph{canonical} quantization $\cW_\fX$ on $\fX$: this is a DQ-algebroid which locally, on a Darboux chart $U \hookrightarrow T^* X$, is equivalent to  $\cW_{T^* X}|_U$. There is an analog of \Cref{uniqueness-of-kashiwara's} for $\cW_\fX$, characterizing it as the unique algebroid stack locally modeled on $\cW_{T^* X}$ and satisfying some axioms; see \cite[p.\ 4]{dagnolo-kashiwara}.\footnote{Similarly to the previously discussed fact that the canonical $E$-algebroid on $\mathbb{P}^*X$ is $\cE_X^{\sqrt{\Omega_X}}$, the canonical DQ-algebroid on $T^* X$ is $W_{T^* X}^{\sqrt{\Omega_X}}$.} There is also a more general classification theorem for DQ-algebroids due to Polesello in \cite[Section 5]{polesello-classification}.

The DQ-algebroid $\cW_\fX$ has the virtue that it exists on any complex symplectic manifold. However, for the conic symplectic manifolds of relevant geometric representation theory, the category of $\cW_\fX$ does not yet provide the correct category; practically because it is defined over $\hbar$ rather than $\C$, and philosophically because it does not incorporate the $\C^\times$-action.  Following \cite{kashiwara-rouquier}, it is desirable to introduce a $\C^\times$-equivariant structure. .

More precisely, the usual setup is as follows: let $(\fX, \lambda)$ be an exact complex symplectic manifold with a  $\C^\times$ action of weight $k$, $k \in \Z \setminus \{0\}$, i.e., a group action $f: \C^\times \times \fX \rightarrow \fX$ such that, if we denote $f_c(-) \coloneqq f(c,-)$, then
$c^* \lambda \coloneqq f_c^* \lambda = c^k \lambda$ for $c \in \C^\times$. Furthermore, there is usually a ``quantization" of this action $f$ which is an automorphism $F$ on $\cW_\fX$, a holomorphic quasi-equivariant structure $F$ on $\cW_\fX$, respecting the action $f$ and satisfying the Frobenius condition:
$$ F_c (\hbar) = c^k \hbar, \ c \in \C^\times.$$
We will thus refer this $F$ as the F-action. By adding a formal $k$-th root of $\hbar$, we can extend $F$ to $\cW_\fX[\hbar^{\frac{1}{k}}]$, which we will abuse the notation and denote it by the same notation, by $F_c(\hbar^\frac{1}{k}) = c \hbar^{\frac{1}{k}}$. To get the correct invariant, one thus first consider the $F_\fX$-equivariant modules $(\cW_\fX[\hbar^\frac{1}{k}], F ) - \mod_{coh}$, and then restrict the usual regular holonomic objects  $(\cW_\fX[\hbar^\frac{1}{k}], F ) - \mod_{coh, rh}$.  We refer the reader to \Cref{apn: equivariance} for details about equivariance and remark that the usual notion extends directly from the situation of algebra of sheaves to algebroids.

In order to relate this invariant with $E$-modules (to which microlocal Riemann-Hilbert is applicable), we need to use work of Petit \cite{petit-frobenius}. His main result \cite[Theorem 6.11]{petit-frobenius} holds in the following general setting: Let $\fX$ be a (complex) symplectic manifold with a $\C^\times$-action $f$ of weight $k$ which is free and proper, so that there is a projection to the contact quotient $p: \fX \rightarrow \fY \coloneqq \fX/\C^\times$. Furthermore, assume there is a DQ-algebroid $\cA_\fX$ equipped with an equivariant structure $F$, respecting the weight $k$ action. For such a $\cA_\fX$, one can define a conical algebra $\cB_\fY$ on $\fY$ by pushing along $p$ the subalgebra of $F$-invariant sections. 

\begin{theorem}[{\cite[Theorem 6.11]{petit-frobenius}}] \label{thm: petit-main}
There is an equivalence of categories 
$$\cB_\fY - \mod_{coh} = (\cA_\fX^{loc}, F) - \mod_{coh}.$$
\end{theorem}

Specializing to the case when $\fX \subseteq T^*X$ is a conic open subset \cite[Section 6.4]{petit-frobenius}, i.e., one which is invariant under the standard weight $1$ action
\begin{align*} 
f_c: \C^\times \times \fX &\rightarrow \fX \\
(c, z, w) &\mapsto (z, c w),
\end{align*}
one can take $\cA_\fX = \cW_{T^* X}(0) |_{\fX}$ to be the subalgebra of $\cW_\fX$ with degree less or equal to $0$. 
This DQ-algebra admits a F-action which is given concretely by  
\begin{align} \label{for: f-cotangent}
F_c: f_c^* \cW_\fX(0) &\xrightarrow{\sim} \cW_\fX(0) \\
\sum_{l \geq 0} f_l(z,w) \hbar^l &\mapsto \sum_{l \geq 0} f_l(z, c \cdot w) c^l \hbar^l.  \notag,
\end{align}
and it extends automatically to $\cW_\fX$.
In this case, $\cB_\fY = \cE_\fY$ on $\fY \subseteq \P^* X$, the algebra of micro-differential operators discussed in \Cref{sec: micro-diff}. In other words:  

\begin{corollary}[{\cite[Proposition 6.16]{petit-frobenius}}]\label{cor: petit-cotangent}
When $\fX \subseteq \dT^* X$ is a conic open set of the cotangent bundle away from the zero section, $\cW_\fX$ is the canonical W-algebra, and $F$ is the canonical action given by (\ref{for: f-cotangent}), there is an equivalence
\begin{equation} 
\cE_\fY - \mod_{coh} = (\cW_\fX, F)-  \mod_{coh}.
\end{equation}
\end{corollary}

The above identification in fact admits the following globalization:

\begin{proposition} \label{w-on-symplectization}
Let $V$ be contact complex manifold and $\tilde{V}$ its complex symplectization. Denote by $\cE_V$ Kashiwara's canonical quantization algebroid on $V$ \cite{kashiwara-quantization-contact} and by $\cW_{\tilde{V}}$ Polesello and Schapira's canonical quantization DQ-algebroid on $\tilde{V}$ \cite[Section 9]{polesello-schapira}. Then we have,
\begin{enumerate}
\item there is a F-action $F_V$  given by (\ref{for: f-cotangent}) on Darboux charts, and
\item there is an identification
\begin{equation}\label{equation:e-w-symplectization}\cE_V - \mod_{coh} = (\cW_{\tilde{V}}, F_V)-  \mod_{coh}.
\end{equation}
\end{enumerate}
\end{proposition}

The proof relies on the following observation: if one begins with a symplecization, then the effect of further contactizing and symplectizing amounts to multiplying with $\dot{T}^*\mathbb{C}$. More precisely:
\begin{lemma}\label{lemma:symp-symp}
Let $V$ be a contact manifold and let $(\tilde{V}, \lambda_V)$ be the symplectization. Let us temporarily denote by $C(-)$ the contactization of an exact symplectic manifold. Then we have \cite[(14)]{CKNS} \begin{align*}
    \widetilde{C\tilde{V}}= (\tilde{V} \times \mathbb{C} \times \C^\times, w(\lambda_V + dz))&= (\tilde{V} \times \dot{T}^*\mathbb{C}, \lambda_V + wdz) \\
    (x, z, w) &\mapsto (w \cdot x, z, w)
\end{align*} 
\qed
\end{lemma}

\begin{proof}[Proof of \Cref{w-on-symplectization}]

Fix a contact manifold $V$ and contemplate the construction of $\cE_V$ and $\cW_{\tilde{V}}$ by gluing \cite[Sec.\ 9]{polesello-schapira}. Fix Darboux charts $U_\alpha \subset \mathbb{P}^*M_\alpha$, $U_\beta \subset \mathbb{P}^*M_\beta$ and a contactomorphism $\phi_{\alpha \beta}: U_\alpha \to U_\beta$. Then we have

\begin{equation}
\begin{tikzcd}
    \dot{T}^*M_\alpha \times \dot{T}^*\mathbb{C} \supset \tilde{U}_\alpha \times \dot{T}^*\mathbb{C} \ar[r, "{(\phi_{\alpha \beta}, id)}"] \ar[d] &\dot{T}^*M_\beta \times \dot{T}^*\mathbb{C} \supset \tilde{U}_\beta \times \dot{T}^*\mathbb{C} \ar[d] \\
    \mathbb{P}^*M \supset U_\alpha  \ar[r, "\phi_{\alpha \beta}"] & \mathbb{P}^*M_\beta \supset U_\beta 
\end{tikzcd}
\end{equation}
We have isomorphisms of algebras $\Phi_{\alpha \beta}: \cE_{\tilde{U}_\alpha} \to \cE_{\tilde{U}_\beta}$. Since $\cE_{\tilde{U}_\alpha \times \dot{T}^*\mathbb{C}}= \cE_{\tilde{U}_\alpha} \hat{\boxtimes} \cE_{\dot{T}^*\mathbb{C}}$, we obtain isomorphisms $$\Phi_{\alpha \beta} \hat{\boxtimes} id: \cE_{\tilde{U}_\alpha \times \dot{T}^*\mathbb{C}} \to \cE_{\tilde{U}_\beta \times \dot{T}^*\mathbb{C}}.$$ The maps $\Phi_{\alpha \beta} \hat{\boxtimes} id$ commute with $\partial_t$, so they can be used to construct $\cW_{\tilde{V}}$. Moreover, if we undo the change of variable for (\ref{for: section-of-W}), we see that the $\C^\times$-action (\ref{for: f-cotangent}) is given by 
\begin{align} \label{for: f-cotangent-homog}
\rho_* \cE_{\P^*(X \times \C),t}(0) &\xrightarrow{\sim} \rho_* \cE_{\P^*(X \times \C),t}(0) \\
\sum_{j \leq 0} f_j(z, \xi, t) &\mapsto \sum_{l \geq 0} f_j(z, \xi, c \cdot t)\notag
\end{align} 
where the $f_j$'s are now homogeneous of degree $j$, and thus they leave the $\dot{T}^*\mathbb{C}$ component invariant and commute with the $\mathbb{C}^\times$-action \eqref{for: f-cotangent}. It follows that this action glues, which is (1). For (2), since both sides sheaves of categories which are are defined from the same gluing data, it is enough to verify the asserted equivalence on Darboux charts, which is \cite[Prop.\ 6.16]{petit-frobenius}.
\end{proof}

We end this subsection with the following auxiliary lemma, which follows immediately from the canonicity of $\cW$ which was already discussed above.

\begin{lemma} \label{lem: covering-W}
Let $\fX$ be a complex manifold and $p: \fY \rightarrow \fX$ a covering map. Then we have 
$$\cW_{\fY} = p^* \cW_{\fX}.$$
\qed
\end{lemma}

\subsection{Canonical F-action through stabilization} \label{subsection:stabilization} 
We begin this section with a lemma and an example which we will subsequently generalize.

\begin{lemma}[{\cite[2.3.2]{kashiwara-rouquier}}] \label{lem: weight-increasing}
Let $f: \C^\times \times \fX \rightarrow \fX$ be a $\C^\times$-action of weight $k$ with a F-action $F$ with and $n \in \N$ a positive integer. We denote by $f^{\circ n}$ the action which is given by $f^{\circ n}_r(x) \coloneqq f_{r^n}(x)$. Then $f^{\circ n}$ is of weight $nk$ with a F-action given by $F^{\circ n}$. Furthermore, there is an equivalence of categories
\begin{align*}
(\cW_\fX[\hbar^{\frac1k}], F)- \mod &\xrightarrow{\sim} (\cW_\fX [\hbar^{\frac{1}{nk} }], F^{\circ n}) - \mod \\
\cM &\mapsto \cW_\fX [\hbar^{\frac{1}{nk}}] \otimes_{\cW_\fX[\hbar^{\frac1k}]} \cM.
\end{align*} 
\end{lemma}

\begin{example} \label{eg: cotangent-affine-line}
Consider the case when $\fX = \dT^* \C$ with the weight $1$ action $f_t(s, \sigma) = (s, t \sigma)$ and the equivariant structure $F_\fX$ on $\cW_\fX$ which is given by
(\ref{for: f-cotangent}). Then \Cref{cor: petit-cotangent} implies that $(\cW_\fX, F_\fX)-\mod = \cE_{P^* X} - \mod$. When restricting to regular holonomic modules, we further have 
$$(\cW_{ \dT^* \C}, F_{ \dT^* \C})-\mod_{rh} = \cE_{P^* \C} - \mod_{rh} \overset{\mRH}{=} \Perv_{\P^* \C}(\P^* \C).$$
There is a $\C^\times$-action by homogeneous symplectomorphisms: 
\begin{align*}
u: \C^\times \times \dT^* \C &\rightarrow \dT^* \C \\
(r, s, \sigma) &\mapsto (r^{-1} s, r \sigma)
\end{align*}
The action is free and proper\footnote{An action of a topological group $G$ on a topological space $X$ is proper if the map $(g, x) \mapsto (x, g \cdot x)$ is proper.}, and it induces $\C^\times$-action by contactomorphisms $\bar{u}$ on $\P^* \C$, which, when identified $\P^* \C$ as $\C$ through
\begin{align*}
\C &\xrightarrow{\sim} \P^* \C \\
s &\mapsto [s,1],
\end{align*}
is given by $\bar{u}(r,s) = r^{-1}s$. Note further that $u$ is Hamiltonian with respect to $H(s,\sigma) = s \sigma$, and, as we will recall in \Cref{sec: hamiltonian-reduction}, taking symplectic reduction on the left hand side and $C^\times$-invariants on the right hand side induces
$$(\cW_{ pt}, F_{ pt})-\mod_{rh} = \C -\mod = \mcsh_{pt}(pt).$$
\end{example}

We also explain how one can reduce to the above case for the weight $k$ case, which through a stabilization trick, will be applied to general cases. For $k \in \Z \setminus \{0\}$, there exists a Liouville structure on $\dT^* \C$ which is given by $\sigma^k ds$. We will denote this Liouville manifold by $(\dT^* \C)^{[k]}$ since it can be obtained by pulling back the usual Liouville structure on $\dT^* \C$ along the $|k|$-fold covering map 
\begin{align*}
p_k: \dT^* \C &\rightarrow \dT^* \C \\
(s, \sigma) &\mapsto (s, \sigma^k).
\end{align*}
Note that we have $f_t^* (\sigma^k ds) = t^k (\sigma^k ds)$ and the same $\C^\times$-action on the underlying complex manifold becomes a weight $k$ action on $(\dT^* \C)^{[k]}$. Now, this action is compatible with the action in weight $1$ action discussed in \Cref{eg: cotangent-affine-line} in the following sense:

\begin{lemma} \label{lem: induced-C-cross-action-geometry}
The $\dT^* \C$, being a quotient of $(\dT^* \C)^{[k]}$ by the deck transformation group $C_k$, admits an action by $\C^\times/C_k$. When identifying $\C^\times/C_k$ as $\C^\times$ by the short exact sequence
$$1 \rightarrow C_{|k|} \rightarrow \C^\times \xrightarrow{t^k} \C^\times \rightarrow 1,$$
can be identified as $(f_t)^{\circ k}$, the $k$-times self-composition of  $f_t$, and so we have $p_k \circ f_t = (f_t)^{\circ k} \circ p_k$.
\end{lemma}

\begin{proof}
This is a special case of \Cref{lem: group-action-geometry}: If there is a short exact sequence of topological group
$$1 \rightarrow K \rightarrow G \rightarrow H \rightarrow 1,$$
with $K$ being discrete, and a covering space $X \rightarrow \overline{X}$ with deck transformation group $K$, then a $G$-action on $X$ induces canonical an $H$-action on $\overline{X}$. It is a direct computation to show that the induced action by $H$, which in our case is also $\C^\times$, is given by $k$-times self-composition.
\end{proof}

The above \Cref{lem: induced-C-cross-action-geometry} implies the following identifications between category of sheaves:

\begin{lemma} \label{lem: induced-C-cross-action-category}
The $*$-pullback $p_k^*:  sh(\dT^* \C) \rightarrow sh( (\dT^* \C)^{[k]})$ induces an equivalence,
$$ p_k^*: sh(\dT^* \C)^{(\C^\times/C_{|k|})} \xrightarrow{\sim} sh( (\dT^* \C)^{[k]})^{\C^\times},$$
between equivariant sheaves. Furthermore, this equivalence respects ring and module structures. In particular, we have the $F$-action $F_{\dT^* \C}^{\circ k}$ induces canonically an $F$-action $F_{(\dT^* \C)^{[k]}}$ such that 
$$ (\cW_{(\dT^* \C)^{[k]}}[\hbar^\frac1k], F_{(\dT^* \C)^{[k]}})-\mod = (\cW_{\dT^* \C}[\hbar^\frac1k], F_{\dT^* \C}^{\circ k})-\mod$$\end{lemma}

\begin{proof}
It is well-known that if $X \rightarrow \overline{X}$ is a covering space with deck transformation group $K$, the $*$-pullback induces an equivalence $sh(\overline{X}) = sh(X)^K$ between usual sheaves on the base and $K$-equivariant sheaves on the cover. As explained in \Cref{prop: group-action-category}. a diagram tracing upgrades it to the equivalence $sh(\overline{X})^H = sh(X)^G$, and the previous \Cref{lem: induced-C-cross-action-geometry} specializes it to our case. The fact that such equivalences respect ring and module structures are explained in \Cref{cor: group-action-category-with-ring}. One could then use either \Cref{lem: covering-W} to conclude that $\cW_{\dT^* \C}$ pulls back to $\cW_{(\dT^* \C)^{[k]}}$.
\end{proof}

\begin{proposition}
$$  (\cW_{(\dT^* \C)^{[k]}}[\hbar^\frac1k], F_{(\dT^* \C)^{[k]}})-\mod = \cE_{P^* \C} - \mod.$$
\end{proposition}

\begin{proof}
By composing with \Cref{lem: weight-increasing} and \Cref{eg: cotangent-affine-line}, we obtain
$$ (\cW_{(\dT^* \C)^{[k]}}[\hbar^\frac1k], F_{(\dT^* \C)^{[k]}})-\mod = (\cW_{\dT^* \C}[\hbar^\frac1k], F_{\dT^* \C}^{\circ k})-\mod = (\cW_{\dT^* \C}, F_{\dT^* \C})-\mod =  \cE_{P^* \C} - \mod.$$
\end{proof}

Let $\fX$ be a complex exact manifold with a $\C^\times$-action $f$ of weight $k$ and a F-action $F$. One obstruction to applying \Cref{thm: petit-main} directly, as in the case of $\dT^* \C$ discussed in \Cref{eg: cotangent-affine-line}, is that the $\C^\times$-action is seldom free. In this section, we will get around this obstruction by stabilizing, i.e., taking a product with $\dT^* \C$. 

The stabilization of $\fX$ by $\dT^* \C$ is the product exact symplectic manifold whose underlying space is simply $\fX \times \dT^* \C$ with its Liouville form given by $\lambda + \sigma^k d s$. The diagonal $\C^\times$-action on the stabilization, which we will denote it by $f_{st}$, is of weight $k$ and is given in coordinates by  
\begin{align*}
f_{st}: \C^\times \times \fX \times \dT^* \C &\rightarrow \fX \times \dT^* \C \\
 (t, x, s, \sigma) &\mapsto (f_t(x), s, t \sigma).
\end{align*}
Indeed, we can check that $t^* (\lambda + \sigma^k ds) = (t^* \lambda) + (t \sigma)^k ds = t^k (\lambda + \sigma ds)$.
We recall that $\cW_{\fX \times (\dT^* \C)^{[k]}} = \cW_{\fX} \uboxtimes \cW_{\dT^* \C}$ is given by the completion of the exterior product $\cW_{\fX} \boxtimes \cW_{\dT^* \C}$. Our construction of the F-action begins with the following change of coordinate:

\begin{lemma} \label{lem: splitting}
Denote by $\widetilde{\fX \times \C} \coloneqq \left(\fX \times \C \times \C^\times, \sigma (\lambda + ds) \right)$ the symplectization of the contact manifold $(\fX \times \C, \lambda + ds)$.
The exact symplectic manifold $\fX \times (\dT^* \C)^{[k]} $ is equivalent to the $|k|$-fold cover $(\widetilde{\fX \times \C})^{[k]}$ which can be explicitly given by
\begin{align*}
q_k: \fX \times \C \times \C^\times &\rightarrow X \times \C \times \C^\times \\
(x, s,\sigma) &\mapsto (x, s, \sigma^k).
\end{align*}
Furthermore, the $\C^\times$-action $f_{st}$ on $\fX \times \dT^* \C$ is identified with the lift along $q_k$ of the canonical $\C^\times$ action of a symplectization $f_{symp}$. 
\end{lemma}

\begin{proof}
Define the map $\varphi$ by  
\begin{align*}
\varphi: (\widetilde{\fX \times \C})^{[k]} &\rightarrow  \fX \times (\dT^* \C)^{[k]} \\
 (x, s, \sigma) &\mapsto (f(\sigma,x), s, \sigma).  
\end{align*}
Then, $\varphi^*( \lambda + \sigma^k ds) = f_\sigma^* \lambda + \sigma^k ds = \sigma^k(\lambda + ds)$.
For $t \in \C^\times$, we compute that
\begin{align*}
&f_{st,t} \circ \varphi(x, s, \sigma) = f_{st,t} (f_\sigma(x), s, \sigma) = (f_t \circ f_\sigma(x), s, t \sigma) \\
&= (f_{t \sigma}(x), s, t \sigma) = \varphi (x, s, t \sigma) = \varphi \circ f_{symp, t}(x, s, \sigma).
\end{align*}
\end{proof}

\begin{corollary}\label{corollary:change-variables} 
There is a canonical F-action $F_{st}$ on $\cW_{\fX \times \dT^* \C}$ which is defined by
$$\varphi^*(\cW_{\fX \times (\dT^* \C)^{[k]}}, F_{st}) = (\cW_{\widetilde{(\fX \times \C})^{[k]}}, q_k^* F_{\fX \times \C}).$$
\end{corollary}

\begin{remark}
The more suitable notation for the F-action should be $F_{(\fX, \lambda, f), st}$ since it depends both on the Liouville structure and the $\C^\times$-action. However, we choose the simpler notation $F_{st}$ since these choice will be fixed for this paper. 
\end{remark}

\begin{proof}
We recall that $\fX \times \C$ is a co-oriented contact manifold whose the contact form is given by $\lambda + ds$, and we thus have a canonical F-action on $\widetilde{\fX \times \C}$ by \Cref{w-on-symplectization}. Since $q_k$ is a covering map, the pullback $q_k^* F_{\fX \times \C}$ is a F-action on $q_k^* \cW_{\widetilde{\fX \times \C}}$, which by \Cref{lem: covering-W} is the same as $\cW_{\widetilde{(\fX \times \C})^{[k]}}$.
\end{proof}

Now we can proceed with the strategy used in \Cref{eg: cotangent-affine-line}. We begin with generalizing \Cref{lem: induced-C-cross-action-category}.

\begin{lemma} \label{lem: induced-C-cross-action-category-stabilizatoin}
There is an equivalence
$$ (\cW_{\widetilde{(\fX \times \C)}^{[k]}}[\hbar^\frac1k], q_k^* F_{\fX \times \C})-\mod =  (\cW_{\widetilde{\fX \times \C}}[\hbar^\frac1k],  F_{\fX \times \C}^{\circ k})-\mod$$
where we denote by $F_{\fX \times \C}^{\circ k} = F_{\fX \times \C} \circ \cdots \circ F_{\fX \times \C}$ its $k$-times self-composition.
\end{lemma}

\begin{proof}
Similar to the discussion in \Cref{lem: induced-C-cross-action-geometry} and \Cref{lem: induced-C-cross-action-category}, we begin with noticing that $ \widetilde{\fX \times \C}$ can be recovered as a quotient $\widetilde{(\fX \times \C)}^{[k]}/C_{|k|}$ and is thus acted by the quotient group $\C^\times/C_{|k|}$. This is a special case of \Cref{lem: group-action-geometry} where the group action $G$ on $X$ is given by $\C^\times$ acting on $\widetilde{(\fX \times \C)}^{[k]}$, $K$ as $C_{|k|}$, and the $k$-times self-composition appears when we identify $\C^\times/C_{|k|}$ with $\C^\times$ using the map $z \mapsto z^k$. 

But then \Cref{lem: induced-C-cross-action-geometry} implies that $\C^\times$-equivariant sheaves on $\widetilde{(\fX \times \C)}^{[k]}$ is identified with $\C^\times/C_k$-equivariant sheaves on $\widetilde{\fX \times \C}$, and \Cref{cor: group-action-category-with-ring} further explains that the identification respects ring and module structures. Thus we obtain the result again by $\C^\times/C_{|k|} = \C^\times$.
\end{proof}

\begin{theorem}\label{thm: F-W-to-E}
There is an equivalence of categories
$$(\cW_{\fX \times (\dT^* \C)^{[k]}}[\hbar^\frac1k], F_{st})-\mod_{coh} = \cE_{\fX \times \C} -\mod_{coh}.$$
When restricting to regular holonomic objects, we further have 
$$(\cW_{\fX \times (\dT^* \C)^{[k]}}[\hbar^\frac1k], F_{st})-\mod_{rh} = \Perv_{\fX \times \C}(\fX \times \C).$$
\end{theorem}

\begin{proof}
\Cref{corollary:change-variables} implies that 
$$\varphi^*: (\cW_{\fX \times (\dT^* \C)^{[k]}}[\hbar^\frac1k], F_{st})-\mod = (\cW_{\widetilde{(\fX \times \C})^{[k]}}[\hbar^\frac1k], q_k^* F_{\fX \times \C})-\mod.$$
The previous \Cref{lem: induced-C-cross-action-category-stabilizatoin} implies that,
$$ (\cW_{\widetilde{(\fX \times \C)}^{[k]}}[\hbar^\frac1k], q_k^* F_{\fX \times \C})-\mod =  (\cW_{\widetilde{\fX \times \C}}[\hbar^\frac1k],  F_{\fX \times \C}^{\circ k})-\mod$$ so we can descend from the cover to the contactization $\fX \times \C$.
To get rid of the $k$-th root, we apply \Cref{lem: weight-increasing} to obtain 
$$ (\cW_{\widetilde{\fX \times \C}}[\hbar^\frac1k],  F_{\fX \times \C}^{\circ k})-\mod =  (\cW_{\widetilde{\fX \times \C}},  F_{\fX \times \C}) -\mod $$
 Lastly, we recall that $F_{\fX \times \C}$ is the canonical F-action given by the symplectization \Cref{w-on-symplectization}, and the same proposition implies that 
$$(\cW_{\fX \times (\dT^* \C)^{[k]}}[\hbar^\frac1k], F_{st})-\mod_{coh} = \cE_{\fX \times \C} -\mod_{coh}.$$
The second statement follows from the fact that regular holonomic objects are preserved under the all the identifications and the main theorem of this paper \Cref{thm:gloabl mrh}. 
\end{proof}

\begin{remark}\label{remark:compatible-action-e-w}
    Recall from the discussion in \Cref{subsection:cont-actions} that the group of contactomorphisms $Cont(V)$ acts on the category $\cE_V - mod$, for $V$ a contact manifold. Since $\cW$-modules are induced from $\cE$-modules, the same considerations imply that the group of symplectomorphisms $Symp(\fX)$ acts on the category $\cW_\fX - mod$, for $\fX$ a symplectic manifold. 

    Suppose now that $\fX = \tilde{V}$. Then we have a canonical inclusion $\Cont(V) \hookrightarrow \Symp(\tilde{V})$. Crucially, the image consists of \emph{equivariant} symplectomorphisms, which commute with the canonical weight-$1$ $\C^\times$ action on the symplectization. Hence $\Cont(V)$ acts on the category $(\cW_{\tilde{V}}, F_V)-mod$ of $F$-equivariant $\cW_{\tilde{V}}$-modules introduced in \Cref{w-on-symplectization}, and the $\Cont(V)$ actions are compatible with \Cref{equation:e-w-symplectization}. In other words, there exists an equivariant structure $(\cW_{\tilde{V}}, \phi_{\tilde{V}})$ such that, for $g \in \Cont(V)$, the equivalence
$\phi_{\tilde{V}, g}: u_g^*   \cW_{\tilde{V}} \xrightarrow{\sim} \cW_{\tilde{V}} $ fixes the inclusion $p^* \cE_V \hookrightarrow \cW_{\tilde{V}}$.
\end{remark}

\subsection{Hamiltonian reduction} \label{sec: hamiltonian-reduction}

In \Cref{subsection:stabilization}, we considered the stabilization of an exact complex-symplectic manifold with a weight $k$ action. In this subsection, we explain how to undo the effect of stabilization by symplectic reduction \cite[2.5]{kashiwara-rouquier}. The general setting is the following:
Let $(M, \omega_M)$ be a complex symplectic manifold, $G$ be a complex lie group, $u: G \times M \rightarrow M$ a symplectic actions and $\mu: X \rightarrow \mathfrak{g}^*$ its moment map. Assume that $G$ acts properly and freely. Then there is a symplectic structure $\omega_Z$ on $p: Z \coloneqq \mu^{-1}(0)/G$ such $p^* \omega_Z = \omega_M |_{ \mu^{-1}(0)/G}$. This is a consequence of the slice theorem, which we will use later in the paper:

\begin{lemma}[{\cite[Lemma 2.7]{kashiwara-rouquier}}] \label{lem: local-form}
Locally on $Z$, the manifold $M$ is isomorphic to $T^* G \times Z$. More precisely, for any point $x \in \mu^{-1}(0)$, there exists a $G$-invariant open neighborhood $U$ of $x$ in $M$ and $a$ $G$-equivariant open symplectic embedding $U \rightarrow T^* G \times T^* \C^n$ compatible with the moment map.
\end{lemma}

Now we specialize to our situation and let $\fX$ be as in the previous subsection. Define a $\C^\times$-action on $\fX \times \dT^* \C$ by
\begin{align}\label{equation:rotation-action-ham}
u: \C^\times \times \fX \times (\dT^* \C)^{[k]} &\rightarrow \fX \times (\dT^* \C)^{[k]} \\
(r, x, s, \sigma) &\mapsto (x, r^{-k} s, r \sigma). \nonumber
\end{align}

\begin{lemma}\label{lemma:ham-rotation-moment}
    \eqref{equation:rotation-action-ham} is a free and proper action by isomorphism preserving strictly the $1$-form.  In fact, the action is Hamiltonian with moment map $\mu=(x, s, \sigma) = s \sigma^k$. Hence $\mu^{-1}(0)/\C^\times= \fX$.
\end{lemma}
\begin{proof}
    We have $u_r^* (\lambda + \sigma^k ds) = \lambda + (r \sigma)^k + d( r^{-k} s) ) = \lambda + \sigma^k ds$, which verifies that the action preserves the Liouville form. That the action is free and proper is elementary.

    To verify the remaining assertion, we may assume $\fX$ is a point since the action is trivial in that component. Let us first assume $k=1$. The moment map of the action of $\C^\times$ on $T^*\mathbb{C}$ induced by the tautological action on the base is $(z,w) \mapsto zw$. Since $(\dot{T}^*\mathbb{C})^{[1]}= \dot{T}^*\mathbb{C}$ is an open subset of $T^*\mathbb{C}$ subset, the moment map of the former is the restriction of the later. If now $k \neq 1$, recall that the symplectic form on $(\dot{T}^*\mathbb{C})^{[k]}$ is pulled back from the standard one under the map $(s, \sigma) \mapsto (s, \sigma^k)$. Hence the moment map is also pulled back (since it is defined by a local condition). Finally, note that $\mu^{-1}(0)= \{(x, 0, \sigma) \} = \fX \times \dT^*_0 \C \subseteq \fX \times (\dot{T}^* \mathbb{C})^{[k]}$, so the final claim follows. 
\end{proof}

As opposite to the F-actions, which is only quasi-equivariant, we will need to consider modules equivariant with respect to $u$ in order to perform symplectic reduction for W-modules; this is the special case, when $\lambda = 0$, of the notion of twisted $G$-equivariant W-modules defined in \cite[2.4.2]{kashiwara-rouquier}. In general, one has to fix a choice of a \textit{quantized moment map} $\mu_\cW$ \cite[2.4.1]{kashiwara-rouquier}, a way in which the Lie algebra $\fg$ acts on the W-algebra and equivariant modules are the full subcategory of  quasi-equivariant modules consisting of those whose structure differentiated to the $\mu_\cW$. 

\begin{example} \label{eg: lie-action-C}
Let $u: \C^\times \times \C \rightarrow C$ be the usual $\C^\times$-action which is given by $u(r,z) = r z$. Identifying $\Lie(\C^\times) = \C$, one can differentiate $u$ to obtain the Lie algebra homomorphism 
\begin{align*}
\C &\rightarrow \cD_{\C}(1) \\
 1 &\mapsto z \frac{\partial}{\partial z}
\end{align*}
In fact, there is a holomorphic $\C^\times$-equivariant structure on $\cD_{\C}$
\begin{align*}
\phi_{can,r}: u_r^* \cD_{\C} &\rightarrow \cD_{\C} \\
h(z) &\mapsto \left( z \mapsto h(rz) \right) \\
\frac{\partial}{\partial z} &\mapsto r^{-1} \frac{\partial}{\partial z}.
\end{align*}
and $\phi_{can,r}(z \frac{\partial}{\partial z}) = (rz) r^{-1} \frac{\partial}{\partial z} = z \frac{\partial}{\partial z}$ is indeed a $\C^\times$-invariant section. Now, we can compose $\mu$ with the standard inclusion $\pi^* \cD_\C \hookrightarrow \cW_{\dT^* \C}$, where $\pi: \dT^* \C \rightarrow \C$ is the projection to the base, and obtain a quantized moment map $\mu_{\cW}$. Passing to the covering by $q_k$, we obtain a quantized moment map  
\begin{align*}
\mu_{\cW}: \Lie(\C^\times) = \C &\rightarrow \cW_{(\dT^* \C)^{[k]}} \\
1 &\mapsto s \sigma^k \hbar^{-1}.
\end{align*}

\end{example}

In general, set
\begin{equation} \label{for: L-symplectic-reduction}
\cL \coloneqq \cW_M / \cW_M \langle \mu_\cW(A)| A \in \fg \rangle,
\end{equation}
i.e., the module obtained by quotient the left ideal generated by the image of quantized Hamiltonian. Then $\cL$ is a coherent $G$-equivariant modules with support contained in $\mu^{-1}(0)$, and definite the lattice $\cL(0)$ similarly. Let $\cW_Z \coloneqq \left( \left( (p_* \operatorname{\cE nd}_\cW(\cL) \right)^G \right)^{opp}$.

\begin{proposition}[{\cite[Proposition 2.8]{kashiwara-rouquier}}] \label{prop: general-symplectic-reduction}
There is a equivalence of categories
$$ \left( \cW_{M} -mod_{good} \right)^u = \cW_{Z} -mod_{good}.$$
Furthermore, assume that there is an F-action $F$ on $\cW_M$ which commutes with the $G$-equivariant structure. Then, there exists a F-action $\bar{F}$, such that
$$\left( (\cW_{M} , F)-\mod_{good} \right)^{u}= (\cW_Z, \bar{F})-mod_{good}.$$
\end{proposition}

As the WKB-algebra $\cW_{\fX \times (\dT^* \C)^{[k]}}$ is an exterior product and the action $u$ only actions on $(\dT^* \C)^{[k]}$ component, we can obtain a quantized moment map on $\cW_{\fX \times (\dT^* \C)^{[k]}}$.  Thus, we can perform symplectic reduction by the above \Cref{prop: general-symplectic-reduction}, which in our current special case states as the following:

\begin{proposition}\label{prop: stabilization-symplectic-reduction}
There is a equivalence of categories
$$ \left( \cW_{\fX \times \dT^* \C} -mod_{good} \right)^u = \cW_{\fX} -mod_{good}.$$
Furthermore, there exists a F-action $F_\fX$, such that
$$\left( (\cW_{\fX \times (\dT^* \C)^{[k]}} , F_{st})-\mod_{good} \right)^{u}= (\cW_\fX, F_\fX)-mod_{good}.$$
\end{proposition}

\begin{proof}
As explained in \Cref{prop: general-symplectic-reduction}, the equivalence in our case is given by considering the object 
$$\cL \coloneqq \cW_{\fX \times \dT^* \C}/ \cW_{\fX \times \dT^* \C}( s \sigma \hbar^{-1}),$$
(here $ \cW_{\fX \times \dT^* \C}( s \sigma \hbar^{-1})$ is the left-ideal generated by $( s \sigma \hbar^{-1})$) and the equivalence is given by
\begin{align*}
\left( \cW_{\fX \times \dT^* \C} -mod_{good} \right)^u &\xrightarrow{\sim} \cW_{\fX} -mod_{good} \\
\cM &\mapsto q_* \sHom_{\cW_{\fX \times \dT^* \C}}(\cL, \cM)
\end{align*}
where $q: \fX \times \dT^*_0 \C \rightarrow  \fX$ is the quotient map. Indeed, as remarked in the same discussion, $\supp(\cL) \subseteq \mu^{-1}(0)$. In other words, in this special case, the symplectic reduction simply kills off the stabilization factor. Since the F-action $F_{st}$ is compatible with the $\C^\times$-action $u$, it induces canonically the F-action action $F_{\fX}$.
\end{proof}

\begin{remark}
Again, the F-action $F_\fX$ depends on both the Liouville structure and the $\C^\times$-action $f$ but we keep the simpler notation $F_\fX$ since there is no risk of ambiguity.
\end{remark}

\begin{remark} \label{rmk: equi-support}
As noted in \cite[2.4.2]{kashiwara-rouquier}, if $\cM$ is a equivariant $W$-module, then $\supp(\cM) \subseteq \mu^{-1}(0)$. In our case, this implies that an $u$-equivariant $\cW_{\fX \times \dT^* \C}$-module is necessarily supported in $\fX \times \dT^*_0 \C$.  
\end{remark}

Now, we've recovered the desired category on the W-module side, we will spend the rest of the section matching the effect of taking invariance on the microsheaf side. As previously observed (see \Cref{lem: splitting}), our standing assumption that $\fX$ carries a $\C^\times$-action of weight $k$ allows us to identify the product $\fX \times (\dot{T}^*\mathbb{C})^{[k]}$ with the $k$-fold cover of the symplectization of the contactization of $\fX$, i.e.\ with $(\widetilde{\fX \times \mathbb{C}})^{[k]}$. Under this identification, \eqref{equation:rotation-action-ham} is induced by a $\C^\times$-action on $\fX \times \mathbb{C}$ by \emph{contactomorphism}. This is the content of the following two lemmas. 

\begin{lemma}\label{lemma:cont-action-ham}
Under $\phi$ in \Cref{lem: splitting}, the action $u$ corresponds to the action
\begin{align*}
v: \C^\times \times (\widetilde{\fX \times  \C})^{[k]} &\rightarrow  (\widetilde{\fX \times  \C})^{[k]} \\
(r, x, s, \sigma) &\mapsto (f_{r^{-1}}(x), r^{-k} s, r \sigma).
\end{align*}
Furthermore, define similarly 
\begin{align*}
\bar{v}: \C^\times \times \widetilde{\fX \times  \C} &\rightarrow \widetilde{\fX \times  \C} \\
(r, x, s, \sigma) &\mapsto (f_{r^{-1}}(x), r^{-k} s, r^k \sigma),
\end{align*}
and we have $q_k \circ v = \bar{v} \circ q_k$ where $q_k$ is the covering map defined in \Cref{lem: splitting}.
\end{lemma}

\begin{proof}
For $(x, s, \sigma) \in (\widetilde{\fX \times  \C})^{[k]}$, we have 
$$\varphi^{-1} \circ u_r \circ \varphi(x, s, \sigma) = \varphi^{-1} \circ u_r(f_{\sigma }(x), s, \sigma) = \varphi^{-1} (f_{\sigma }(x), r^{-k} s, t \sigma)  = (f_{r^{-1}}(x), r^{-k} s, r \sigma).$$
\end{proof}

\begin{lemma} \label{lem: induced-contact-action}
Since the action $\bar{v}$, by homogeneous symplectomorphisms, commutes with the weight $1$ action $f_{symp}$, it induces a $\C^\times$-action on $\fX \times \C$ by contactomorphisms
\begin{align*}
a: \C^\times \times \fX \times \C &\rightarrow \fX \times \C \\
(r, x, s) &\mapsto (f_{r^{-1}}(x), r^{-k} s).
\end{align*}
Here, we identify the quotient $\widetilde{\fX \times \C}/f_{symp}$ with $\fX \times \C$ by 
\begin{align*}
\widetilde{\fX \times \C}/f_{symp} &\rightarrow \fX \times \C  \\
[x,s, t] &\mapsto (x,s).
\end{align*}
\end{lemma}

To summarize the above discussion, the Hamiltonian  $\C^\times$-action \eqref{equation:rotation-action-ham} is induced by an action of  $\C^\times \subset Cont(\fX \times \C)$ and we have the following equivalence: 

\begin{proposition}\label{prop: equi-F-W-to-E}
There is an equivalence of categories
$$(\cW_\fX, F_\fX)-mod_{good} = (\cE_{\fX \times \C} -mod_{good} )^a.$$
Furthermore, any object $\cN \in (\cE_{\fX \times \C} -mod_{good} )^a$ is supported in $\fX \times \{0\}$.
\end{proposition}

\begin{proof}
The passage from W-modules to E-modules in \Cref{thm: petit-main} is given by taking $F$-invariant sections. Since the $\C^\times$-actions $f$ and $u$ commute, the first half of \Cref{thm: F-W-to-E} induces this equivariant version. For the statement about the support, as remarked in \Cref{rmk: equi-support}, objects in $\left( \cW_{\fX \times \dT^* \C} -mod_{good} \right)^u$ are supported in $\fX \times \dT^*_0 \C$. Under $\varphi$, the set corresponds to $\fX \times \{0\} \times \C^\times$, which is sent to $\fX \times \{0\}$ under the quotient map in \Cref{lem: induced-contact-action}.
\end{proof}

When further restricting to regular holonomic modules, we obtain the following result:

\begin{theorem} \label{thm: w-modules-to-mcsh}
Let $\fX$ be an exact complex symplectic manifold with a  $\C^\times$ action of weight $k$, $k \in \Z \setminus \{0\}$. There is an equivalence
\begin{equation}\label{equation:main-rh-symplectic} (\cW_\fX, F_\fX)-\mod_{rh} = (\mcsh_{\fX, \mathbb{C}-c}(\fX))^\heartsuit. 
\end{equation} 
\end{theorem}
\begin{proof}

By \Cref{thm: equi-gloabl-mrh} and \Cref{prop: equi-F-W-to-E}, we have the equivalence 
$$(\cW_\fX, F_\fX)-\mod_{rh} =\Perv_{\fX \times \C}(\fX \times \C)^a,$$ 
and objects on the right hand side are supported in $\fX \times \{0\}$. Following \cite[Thm.\ 6.3]{CKNS}, we have
$$\mu_\C sh_{\fX}:= (\mathbb{P}\mu sh_{\fX \times \C, \fX \times \{0\}})^{a} \hookrightarrow \mu sh_{\fX}.$$ \cite[Thm.\ 1.4]{CKNS} provides a $t$-structure on ${\mu_\C sh}_{\fX, \C-c}$ whose heart is denoted by $({\mu_\C sh}_{\fX, \C-c})^\heartsuit$. Since taking invariants commutes with taking hearts (both are limits), we have $$(\mathbb{P}erv_{\fX \times \{0\}, \C-c})^{a}= ({\mu_\C sh}_{\fX, \C-c})^\heartsuit.$$
\end{proof}

\appendix
\section{Recollections from Morita theory} 
 \label{app: morita}

To study the untwisting procedure in \Cref{lem: equivalence-modules} and \Cref{twisting-algebras-versus-twisting-modules}, we consider the following situation:
Let $R$ and $S$ be rings and $f: R \rightarrow S$ be a ring homomorphism. Viewing $S$ as a right $R$-module and a left $S$-module, we obtain a tensor-forgetful adjunction
$$S \otimes_R (-): R-mod \rightleftharpoons S-mod: \prescript{\ }{f}{(-)}$$
where for a  $S$-module $N$, we use $\forcoe{f}{N}$ to denote the $R$-module whose structure is given by $r \cdot_f n \coloneqq f(r)n$ for $r \in R$ and $n \in N$. We note that, if $\alpha: N_1 \rightarrow N_2$ is a $S$-module homomorphism, the functions $\forcoe{f}{\alpha}$ and$\forcoe{g}{\alpha}$ are both set-theoretically the same as $\alpha$.

\begin{lemma}
Let $f, g: R \rightarrow S$ be two ring homomorphisms. Any natural transformation $T: \forcoe{f}{(-)} \rightarrow \forcoe{g}{(-)}$ has the form $L_s$, for $s \in S$ such that $s f(r) = g(r) s$ for all $r \in R$, where $L_s$ is the family of function which is defined by left multiplication
\begin{align*}
L_s(N):N &\rightarrow N \\
n &\mapsto s n.
\end{align*}
In particular, the functors $\forcoe{f}{(-)}$ and $\forcoe{g}{(-)}$ are equivalent if and if only there is $s \in S^\times$ such that $g = \Ad(s) \circ f$.
\end{lemma}

\begin{proof}
The natural transformation $T: \forcoe{f}{(-)} \rightarrow \forcoe{g}{(-)}$ is determined by $T(N)(n)$ for any $N \in S-mod$ and $n \in N$.
Fix such an $n \in N$, the function 
\begin{align*}
\phi_n: S &\rightarrow N \\
s &\mapsto sn
\end{align*}
is an $S$-module homomorphism. Thus, we have he equality $T(N) \circ \forcoe{f}{\phi_n} = \forcoe{g}{\phi_n} \circ T(S)$. In other words, for any $s \in S$, we have $ T(N)(sn) =   \left(T(S)(s)\right) n$. Taking $s = 1$, we see that $T(N)(n) = \left(T(S)(1)\right)n$ is determined by the element $s_T \coloneqq T(S)(1) \in S$, or $T = L_{s_T}$. To see that $s_T$ satisfies the desired property, we take $N = S$ and $n = f(r)$ for $r \in R$. Then we see that
$$s_T f(r) = T(S)\left(f(r)\right) = T(S)\left( r \cdot_f 1\right) = r \cdot_g T(S)(1) = g(r) s_T.$$
One can similarly show that $T$ is a natural equivalence if $s_T \in S^\times$ and in this case $g(r) = s_T^{-1} f(r) s_T$, or $g = \Ad(s_T) \circ f$.
\end{proof}

\begin{corollary} \label{equivalence-to-identity}
For a ring $R$, a ring automorphism $f: R \rightarrow R$ induces a functor $\forcoe{f}{(-)}$ equivalent to the identity $\id_{R -mod}$ if and only if $f = \Ad(r)$ for some $r \in R^\times$. In this case, we have $$L_r: \id_{R-mod} \xrightarrow{\sim} \forcoe{\Ad(r)}{(-)}.$$
\end{corollary}

Another fact which we will use is that a $(R,S)$-bimodule with a common left $R$-module and right $S$-module generator $t$ gives rise to a anti-ring homomorphsm.
The fact is well-known but we spell out the details since we will need it in \Cref{star-involution-D-modules}. We also remark that this is the same procedure in\cite[Theorem 11.4.9]{kashiwara-schapira} to obtain the quantized contact transform. 
 
\begin{lemma} \label{lem: bimodule-ring-morhpism}
Let $M$ be a $(R,S)$-bimodule with a common generator $t$. That is the left $R$-module morphism $R \rightarrow M, r \mapsto r \cdot t$ and the right $S$-module homomorphism $S \rightarrow M, s \mapsto t \cdot s$ are both an isomorphism.
Then, $M$ defines an ring homomorphism $f_M: R \rightarrow S$.
\end{lemma}

\begin{proof}
For any $r \in R$, by the assumption on $t$, there exists a unique $f_M(r) \in S$ such that $r \cdot t = t \cdot f_M(r)$, which one can check satisfies $f_M(1_R) = 1_S$ and $f_M(r_1 r_2) = f_M(r_1) \cdot f_M(r_2)$.
\end{proof}

We also need to consider passing to `op' modules. That is, a $(R,S)$-bimodule $M$ can be viewed tautologically as an $(S^{op}, R^{op})$-bimodule, which we denote it by $M^{op}$.
For the sheaf-theoretical situation, since convolution has a prefer direction, the coordinates have to be swapped.
That is, we let $f: X \rightarrow Y$ be a homeomorphism between topological spaces and we denote by $v: X \times Y \xrightarrow{\sim} Y \times X$ the coordinate swapping map $v(x,y) = (y,x)$.
Let $S \in sh(X), R \in \sh(Y)$ be ring-valued sheaves and $M \in sh(X \times Y)$ be a $(p_2^* R, p_1^*S)$-bimodule supported on the graph $\Gamma_f$ so $v^* M$ is supported in $\Gamma_{f^{-1}} \subset Y \times X$.

\begin{lemma} \label{lem: opposite-sheaf-bimodule}
Convolution defines functors $$M \circ_S (-): S-mod \rightarrow R-mod$$ and  $$v^* M^{op} \circ_{R^{op}} (-): R^{op} -mod \rightarrow S^{op} - mod.$$
\end{lemma}

Finally, we mention that Morita theory generalizes mostly naturally in the algebroid set up. Our main reference for this discussion will be \cite{dagnolo-polesello}.
For an algebroid $\cA$ on $X$ (over $\C$), the category of modules is defined by 
\begin{equation}
\cA -mod \coloneqq \Fun_{\C}(\cA, \C_X -mod).
\end{equation}
Note that such an assignment naturally organizes to a sheaf $\cA- \sMod$.
Similarly to the case of algebras, for two algebroids $\cA$, $\cB$, one can form their tensor algebroid $\cA \otimes_{\C} \cB$
and a $(\cA, \cB)$-bimodule is an object $\cP \in \cA \otimes_{\C} \cB^{op}$. A $(\cA,\cB)$-bimodule is said to be invertible if there exists a $(\cB, \cA)$-bimodule $\cQ$ such that 
$$ \cP \otimes_{\cB} \cQ \cong \Delta_\cA \, \text{and} \, \cQ \otimes_{\cA} \cP \cong \Delta_\cB$$
where $\Delta_\cA$ is the diagonal bimodule which is given by $\left( (x,y) \mapsto \cA(x,y) \right)$.
An equivalence of bimodules $f: \cA \xrightarrow{\sim} \cB$ naturally induces an equivalence $f \circ (-): \cB -mod \xrightarrow{\sim} \cA -mod$.
In fact, this assignment organizes to a $2$-functor
\begin{equation}\label{for: algebroid-eq-to-module-eq}
\Equiv_{\C}(\cA, \cB) \rightarrow \Equiv_{\C}(\cB -mod, \cA -mod).
\end{equation}

\begin{definition}[{\cite[Definition 3.4.1]{dagnolo-polesello}}] \label{def: picard-good}
An algebroid $\cA$ is said to be Picard good if all invertible $\cA \otimes_\C \cA^{op}$ bimodules are locally free of rank one over $\cA$ (or, equivalently, over $\cA^{op}$).
\end{definition} 

\begin{proposition}[{\cite[Proposition 3.4.3]{dagnolo-polesello}}] \label{prop: algebroid-eq-=-module-eq}
When $\cA$ is Picard good and $\cB$ is locally equivalent to $\cA$, the \eqref{for: algebroid-eq-to-module-eq} is an equivalence.
\end{proposition}

\begin{proof}
D'Agnolo and Polesello prove the more general statement, in \cite[Proposition 3.3.7]{dagnolo-polesello}, that for linear stacks the image of \eqref{for: algebroid-eq-to-module-eq} consists of bimodules locally free of rank one over $\cA$.
In our case, the classification of the image is implied by the discussion at the beginning of this Appendix section, and \Cref{lem: bimodule-ring-morhpism} gives the explicit way to reconstruct the algebroid map from the bimodule on small open set when there is a section.
\end{proof}

\section{Recollections from algebraic topology} \label{apn: homotopy-section}

Here we trace through some relationships classifying certain null-homotopies and homotopies of null-homotopies etc. which 
we have needed in the text. Let us first recall the basic result on (de-)looping and commutativity.  It is due to Boardman-Vogt and May; we give references to the treatment in Lurie's books. 

\begin{theorem}[{\cite[Theorem 5.2.6.10]{lurie-ha}}] \label{delooping}
For $\infty > k \geq 0$, the functor of taking $k$-fold loop space
\begin{align*}
\Spc_*^{k \geq} &\xrightarrow{\sim} Mon_{\E_k}^{gp}(\Spc) \\
(X,*) &\mapsto \Omega_*^k (X)
\end{align*}
is an equivalence between $k$-connective space and group-like $\E_k$-monoid in spaces.
\end{theorem}

As a corollary, one obtain a similar statement for group-like $\E_\infty$-monoid in spaces; they are the abelian group objects in this setting:

\begin{corollary}[{\cite[Corollary 5.2.6.27]{lurie-ha}}]\label{infinite-loop-space}
There exists an equivalence
$$ \Omega^\infty: \Sp^{cn} \xrightarrow{\sim} Mon_{\E_\infty}^{gp}(\Spc)$$
between connective spectra, spectra with no non-trivial negative homotopy group, and group-like $\E_\infty$-monoid in spaces.
\end{corollary}

\begin{example}
Any abelian group $A$ can be seen as a discrete topological group and hence an abelian group object in spaces.
The theorem above thus implies that there exists a $k$-connective space $B^k(A)$ such that
$ \Omega_*^l B^k(A) = B^{k - l}(A)$, for $l \leq k$ where we use the notation $B^0 A = A $. 
In particular, $\Omega_* B(A) = A$ as a group in spaces.
\end{example}

\begin{example}
By Bott periodicity, $O$, $U$, and $U/O$ can be regarded as an object in $Mon_{\E_\infty}^{gp}(\Spc)$ and thus in $\Sp^{cn}$.
However, taking $\Sp^{cn} \hookrightarrow \Sp$ is invariant undertaking cofibers and we thus have fiber sequences
$U \rightarrow U/O \rightarrow BO$, $U/O \rightarrow BO \rightarrow B^2 U$, etc..
\end{example}

\begin{example} \label{kernel-of-space-group} 
Let $G \in \E_1^{gp}(\Spc)$ and $A \in \E_2^{gp}(\Spc)$. We mention in the last example that $Mon_{\E_\infty}^{gp}(\Spc)$ in general is only closed under taking cofibers. However, in case when there is a group homomorphism $\alpha: G \rightarrow BA$, its fiber $K \coloneqq \fib(\alpha)$ in fact admits a group structure. Indeed, passing to connective spaces by \Cref{delooping}, we need to argue that the fiber $F \coloneqq \fib(BG \rightarrow B^2 A)$ has a vanishing $\pi_0$. But this follows from the homotopy long exact sequence
$$ \cdots \rightarrow \pi_1(B^2 A)  \rightarrow \pi_0(F) \rightarrow \pi_0(BG) \rightarrow \cdots $$
and the fact that $\pi_1(B^2 A) = \pi_0(BA) = \{*\}$.
\end{example}

\begin{remark} \label{null-homotopies-versus-sections}
Let $G \in Mon_{\E_1}^{gp}(\Spc)$ be a group in spaces.
Consider a nice topological space $X$, e.g., a CW complex, and a classifying map $\eta: X \rightarrow BG$, classifying the principal $G$-bundle $\cP_\eta$.
We know that a null-homotopy $\tau$ of $\eta:X \rightarrow BG$, i.e., an equivalence $\tau: \eta \xrightarrow{\sim} *$ to the constant map,
corresponds to a section, which we will abuse the notation and denote it by $\tau: X \rightarrow \cP_\tau$.
For two of such sections $\tau_1$ and $\tau_2$, if we denote by $\tau_1^\rev: * \xrightarrow{\sim} \eta$ the reverse equivalence, then the concatenation $\tau_2 \# \tau_1^\rev$ will be a loop at the constant map $*$.
Thus, 
$$\tau_2 \# \tau_1^\rev \in \Omega_* Map(X,BG) = Map(X, \Omega_* BG) = Map(X, G)$$
is given by a ``$G$-valued function" $g: X \rightarrow G$. 
Now, the $k = 1$ case of \Cref{delooping} implies that $\Omega_* B G = G$ as groups, i.e., concatenation of loops in $B G$ is given the group multiplication of $G$. 
This implies that, as sections, $\tau_2 = g \cdot \tau_1$ where $\cdot$ comes from the $G$-action on $\cP_\tau$.
\end{remark}

Now let $(N,\xi)$ be a manifold $N$ with a symplectic/complex vector bundle $\xi$, or equivalently, a map $\xi: M \rightarrow B U(n)$.
Abuse the notation and denote again by $\xi$ the composition
$$M \xrightarrow{\xi} B U(n) \rightarrow BU \rightarrow B U/O \rightarrow B O^2 \xrightarrow{B w_2} B^3 (\Z/2).$$
Here, we note that the universal second Stiefel-Whitney, $w_2: BO \rightarrow B^2 \Z/2$, is induced by the exact sequence
$1 \rightarrow \Z/2 \rightarrow \Pin \rightarrow O \rightarrow 1$.
We note that since $BU \rightarrow B U/O \rightarrow B O^2$ is a fiber sequence, it is canonically null-homotopic,
and we will use the notation $\gamma_{can}$ to denote the induced null-homotopy for $\xi$.

Assume now that there exists an open cover $\{U_\alpha\}$ with an open embedding $f_\alpha: U_\alpha \hookrightarrow T^* M_\alpha$,
for some manifold $M_\alpha$ of dimension $n$, so that the vector bundle $\phi_\alpha: T^* M_\alpha \rightarrow BO(n)$ given by the fiber 
$\phi_\alpha(x_\alpha,\xi_\alpha) = T^*_{x_\alpha} M_\alpha$ is a polarization of $\xi$ on $U_\alpha$, ie., the following diagram commutes:

$$
\begin{tikzpicture}
\node at (0,2) {$T^* M_\alpha$};
\node at (4,2) {$BO(n)$};
\node at (0,0) {$U_\alpha$};
\node at (2,0) {$M$};
\node at (4,0) {$BU(n)$};

\draw [->, thick] (0.6,2) -- (3.3,2) node [midway, above] {$\phi_\alpha$};
\draw [right hook->, thick] (0.4,0) -- (1.7,0) node [midway, above] {$ $};
\draw [->, thick] (2.4,0) -- (3.3,0) node [midway, above] {$\xi$};

\draw [<-right hook, thick] (0,1.7) -- (0,0.3) node [midway, right] {$ $}; 
\draw [->, thick] (4,1.7) -- (4,0.3) node [midway, right] {$ $};
\end{tikzpicture}
$$ 
Passing to stabilization, this fiber polarization provides another null-homotopy for $\xi |_{U_\alpha}$ since $BO \rightarrow BU \rightarrow B U/O$ is a fiber sequence 
and we abuse the notation and still use $\phi_\alpha$ to denote it.
Denote by $\tau_{can}^\rev$ the reverse homotopy from the constant map to $\xi_{U_\alpha}$.

\begin{proposition} \label{second-stiefel-whitney}
The concatenation $\phi_\alpha \# \tau_{can}^\rev \in \Omega_* \Map \left(U_\alpha, B^3(\Z/2) \right) = \Map \left(U_\alpha, B^2(\Z/2) \right)$ 
is given by $\phi_\alpha \# \tau_{can}^\rev = f_\alpha^* w_2(M_\alpha)$.
\end{proposition}

\begin{proof}
We note that both null-homotopies are from composing with null-homotopies which happens in the sequence
$$BO \rightarrow BU \rightarrow B U/O \rightarrow B^2 O$$
and hence the following lemma implies that the corresponding point in $ \Map \left(U_\alpha, B^2(\Z/2) \right)$  is given by the composition
$$U_\alpha \xrightarrow{f_\alpha} T^* M_\alpha \xrightarrow{\phi_\alpha} BO \xrightarrow{w_2} B^2 (\Z/2),$$
i.e., $w_2(M_\alpha)$.
\end{proof}

\begin{lemma}
Let $\cC$ be a stable category, $A \xrightarrow{\alpha} B \xrightarrow{\beta} C$ 
be a fiber sequence and choose, up to a contractible ambiguity,
a homotopy $h: \beta \circ \alpha = 0$.
Extend the sequence and get $B \rightarrow \xrightarrow{\beta} C \xrightarrow{\gamma} A[1]$ with a similar choice of homotopy $g$.
Then under the identification $\Omega_0 \Map(A,A[1]) = \Map(A,A)$, the point $(g^\rev \bigcirc \alpha) \# (\gamma \bigcirc h)$ corresponds to $\id_A$. 
Here $`\bigcirc'$ denotes the horizontal composition between a $1$-morphism and a $2$-morphism.
\end{lemma}

\begin{proof}
The follows from tracing through universal properties and the fact that composition of pullbacks is a pullback, i.e., $\Omega_0 A[1] = A$ is given by both of the following diagrams:

$$
\begin{tikzpicture}


\node at (-10,2) {$A$};
\node at (-6,2) {$0$};
\node at (-10,0) {$0$};
\node at (-6,0) {$A[1]$};

\node at (-9.6,1.5) {$\ulcorner$};


\draw [->, thick] (-9.6,2) -- (-6.3,2) node [midway, above] {$ $};
\draw [->, thick] (-9.6,0) -- (-6.5,0) node [midway, above] {$ $};

\draw [->, thick] (-10,1.7) -- (-10,0.3) node [midway, right] {$ $}; 
\draw [->, thick] (-6,1.7) -- (-6,0.3) node [midway, right] {$ $}; 

\node at (-5,1) {$=$};

\node at (-4,2) {$A$};
\node at (0,2) {$B$};
\node at (4,2) {$0$};
\node at (-4,0) {$0$};
\node at (0,0) {$C$};
\node at (4,0) {$A[1]$};

\node at (-3.6,1.5) {$\ulcorner$};
\node at (0.4,1.5) {$\ulcorner$};


\draw [->, thick] (-3.6,2) -- (-0.3,2) node [midway, above] {$ $};
\draw [->, thick] (0.4,2) -- (3.7,2) node [midway, above] {$ $};
\draw [->, thick] (-3.6,0) -- (-0.3,0) node [midway, above] {$ $};
\draw [->, thick] (0.4,0) -- (3.5,0) node [midway, above] {$ $};

\draw [->, thick] (-4,1.7) -- (-4,0.3) node [midway, right] {$ $}; 
\draw [->, thick] (0,1.7) -- (0,0.3) node [midway, right] {$ $}; 
\draw [->, thick] (4,1.7) -- (4,0.3) node [midway, right] {$ $};
\end{tikzpicture}
$$ 

\end{proof}

Denote by $B^2(\Z/2)(\xi)$ the ($\infty$-)principal $B^2(\Z/2)$-bundle over $M$ classify by $\xi: M \rightarrow B^3(\Z/2)$.
Recall that null-homotopies of $\xi: M \rightarrow B^3(\Z/2)$ corresponds to sections of $B^2(\Z/2)(\xi) \rightarrow M$.
More precisely, since we have the pullback diagram

$$
\begin{tikzpicture}
\node at (-6,2) {$M \times B^2(\Z/2)$};
\node at (0,2) {$B^2(\Z/2)$};
\node at (5,2) {$*$};
\node at (-6,0) {$M$};
\node at (0,0) {$*$};
\node at (5,0) {$B^3(\Z/2)$};

\node at (-5.4,1.5) {$\ulcorner$};
\node at (0.6,1.5) {$\ulcorner$};


\draw [->, thick] (-4.5,2) -- (-1,2) node [midway, above] {$ $};
\draw [->, thick] (1,2) -- (4.7,2) node [midway, above] {$ $};
\draw [->, thick] (-5.5,0) -- (-0.3,0) node [midway, above] {$ $};
\draw [->, thick] (0.4,0) -- (4,0) node [midway, above] {$ $};

\draw [->, thick] (-6,1.7) -- (-6,0.3) node [midway, right] {$ $}; 
\draw [->, thick] (0,1.7) -- (0,0.3) node [midway, right] {$ $}; 
\draw [->, thick] (5,1.7) -- (5,0.3) node [midway, right] {$ $};
\end{tikzpicture},
$$
where we use $\Omega_* B^3(\Z/2) = B^2(\Z/2)$ for the right square,
a null-homotopy identifies $B^2(\Z/2)(\xi)$ with $M \times B^2(\Z/2)$, which admits an obvious section given by $M = M \times \{*\} \hookrightarrow M \times B^2(\Z/2)$.
In our case, denote by $\tau_{can}$ and $\phi_\alpha$, the sections which corresponds to the homotopies $h_{can}$ and $h_\alpha$.
Since $\Omega_* B^3(\Z/2) = B^2(\Z/2)$ as groups by Theorem \ref{delooping}, the previous Proposition \ref{second-stiefel-whitney} implies that

\begin{equation} \label{w2-twisted-section}
\phi_{\alpha} = \left(f_\alpha^* w_2(M_\alpha) \right) \cdot \tau_{can}
\end{equation}

where we view $f_\alpha^* w_2(M_\alpha) \in \Map \left(U_\alpha, B^2(\Z/2) \right)$ and the multiplication $\cdot$
is induced from the principal bundle structure $a: B^2(\Z/2) \times B^2(\Z/2)(\xi) \rightarrow B^2(\Z/2)(\xi)$.

\section{Recollections from symplectic and contact geometry} \label{appendix:contact-facts}

\begin{definition}
 A contact manifold $V=(V, \xi)$ is the data of a manifold $V$ of dimension $2n+1$ along with a codimension-$1$ sub-bundle $\xi \subset TV$ which is maximally non-integrable. Concretely, this means that we can locally write $\xi = \ker \alpha$, for some $1$-form $\alpha$ having the property that $\alpha \wedge (d\alpha)^n \neq 0$. 
\end{definition}

If $(V, \xi)$ is a contact manifold, then $\xi \to V$ has the structure of a conformally symplectic vector bundle. (Indeed, if $f\alpha = \alpha'$, then $d\alpha' = df \wedge \alpha + f d\alpha$, so $d\alpha'|_\xi = f d\alpha |_\xi$.)

\begin{definition}
Let $(V, \xi)$ be a contact manifold. A submanifold $W \subset V$ is said to be: 
\begin{itemize}
    \item a \emph{contact submanifold} if $\xi \cap TW \subset TV$ is a contact structure on $W$
    \item \emph{isotropic} if $TW \subset \xi \subset TV$  
    \item \emph{coisotropic} if $(TW \cap \xi)$ is an coisotropic subspace of $(\xi, d\alpha)$ at all points in $W$, for some (equivalently any) locally defined $1$-form $\xi= \ker \alpha$.
\end{itemize}
There are obvious notions of contact/isotropic/coisotropic embeddings. 
\end{definition}

If $W \subset V$ is coisotropic, the \emph{characteristic foliation} is the foliation integrating the distribution $(TW \cap \xi)^{d\alpha} \subset TW$.

Given a contact submanifold $(V, \xi= \eta \cap TV) \subset (W, \eta)$, we let $\xi^\perp \subset \eta|_V$ be the orthogonal complement of $\xi \subset \eta|_V$. Note that $\xi^\perp$ is inherits a canonical conformal symplectic structure and we have $\xi^\perp \oplus TV= TW|_V$ \cite[Sec.\ 2.5.3]{geiges2008introduction}. We call $\xi^\perp$ the \emph{conformal symplectic normal bundle} of $V \subset W$. 

The above definitions are equally valid in the category of real manifold with smooth maps and the category of complex manifolds with holomorphic maps.  \emph{However, for the remainder of this appendix, we will exclusively work in the category of real contact manifolds.} We refer to \cite[Section 2]{CKNS} for background on the complex case. 

Given a contact manifold $(Y, \xi)$, its \emph{symplectization} is the manifold $SY:= \{ \alpha \in T^*Y \mid \operatorname{ker} \alpha = \xi\}$ which is symplectic with respect to the restriction of the canonical $1$-form on $T^*Y$. We say that $Y$ is co-orientable iff $SY$ is disconnected; a co-orientation amounts to labeling one of the two components as positive. 

An exact symplectic manifold $(M, \lambda)$ is said to be \emph{Liouville} if it admits a proper embedding $S^+Y \hookrightarrow (M, \lambda)$ of the positive symplectization of a closed contact manifold. 
   
\begin{lemma}\label{lemma:exact-normal}
Let $\iota: N \hookrightarrow M$ be an embedding and consider the exact sequence of bundle maps $0 \to \operatorname{ker} \iota^* \to T^*M|_N \xrightarrow{\iota^*} T^*N \to 0$. The space of sections of $\iota^*$ is contractible. Any section $\sigma$ induces an inclusion of Liouville manifolds $T^*N \hookrightarrow T^*M$. 
\end{lemma}

\begin{definition}
    Let $(W, \xi)$ be a co-orientable contact manifold. A (possibly time-dependent) vector field $Z$ is called a \emph{contact vector field} it its flow preserves the contact structure. 
\end{definition}
\begin{lemma}\cite[Thm.\ 2.3.1]{geiges2008introduction}
    Let $(W, \xi)$ be a co-orientable contact manifold. A choice of contact form $\xi = \ker \alpha$ induces a bijection of sets
    \begin{equation}
        \{\text{contact vector fields on } W \} \leftrightarrow C^\infty(W)
    \end{equation}
This bijection takes $Z \mapsto \alpha(Z)$. The inverse sends a Hamiltonian $H$ to a contact vector field $Z$ uniquely defined by the equations $\alpha(Z_H)=H$ and $i_{Z_H} d\alpha =dH(R_{\alpha})\alpha-dH$, where $R_\alpha$ the Reeb vector field.
\end{lemma}

\begin{corollary}\label{corollary:extend-isotopy}
    Let $(W, \xi)$ be a co-orientable contact manifold and let $\phi_t: (V, \zeta)\hookrightarrow (W, \xi), t \in [0,1]$ be a $1$-parameter family of codimension zero contact embeddings. Then, after possibly shrinking $W$, $\phi_t$ extends to a family of contactomorphisms of $(V, \xi)$. More precisely, there exists a family of contactomorphisms $\tilde{\phi}_t: (W, \xi) \to (W, \xi), t \in [0,1]$ such that $\tilde{\phi}_t \circ \phi_0 = \phi_t$ holds on any compact subset of $W$. 
\end{corollary}
\begin{proof}
    Fix a contact form $\xi = \ker \alpha$ and consider the contact vector field $Z_t:= \frac{d}{dt} \phi_t$ defined on $\phi_t(W)$. Let $H_t:= \alpha(Z_t)$ be the corresponding family of Hamiltonians and let $\tilde{H}_t$ be an extension of $H_t$ to all of $W$. Now define $\tilde{\phi}_t$ to be the flow of $\tilde{H}_t$. 
\end{proof}

\begin{definition}\label{definition:thickening}
    Let $(V^{2n+1}, \xi\cap TV) \hookrightarrow (W^{2N+1}, \xi)$ be a contact embedding. A \emph{thickening} is a (germ of a) co-isotropic submanifold of dimension $(\dim W- \dim V)/2 + \dim V$ which contains $V$ and whose characteristic foliation is non-singular and transverse to $V$.  
\end{definition} 
We typically abuse language by identifying $V$ with its image. 
\begin{lemma}\label{lemma:thick-to-distribution}
With the notation of \Cref{definition:thickening}, suppose that $C$ is a thickening of $V$. Then the tangent space of the characteristic foliation defines a Lagrangian subbundle of the conformal symplectic normal bundle $(\xi \cap TV)^\perp \subset \xi|_V$. 
\end{lemma}
\begin{proof}
Given $p \in V \subset C$, let $L$ be the unique leaf of the foliation passing through $p$. Then $T_pL \subset (\xi \cap TV)^\perp_p$ is a half-dimensional subspace. It is also isotropic with respect to the conformal symplectic structure $d\alpha$ (because $d\alpha(V, -)$ vanishes on $\xi_p$ and hence on $T_pL$, for any $V \in T_pL$). 
\end{proof}

Let $\mathfrak{Thick}(V \subset W)$ be the space of thickenings of $V$. Let $\mathfrak{Lag}(V \subset W)$ be the space of Lagrangian distributions, that is, section of the bundle $LGr((\xi \cap TV)^\perp) \to V$. By \Cref{lemma:thick-to-distribution}, there is a forgetful map 
\begin{equation}\label{equation:forget-thick} 
\mathfrak{F}: \mathfrak{Thick}(V \subset W) \to \mathfrak{Lag}(V \subset W).
\end{equation}

We ultimately wish to prove that \eqref{equation:forget-thick} is a homotopy equivalence. This will require some preparatory lemmas. 

\begin{lemma}\label{lemma:local-contact-form} Let $(V, \xi)$ be a co-orientable contact manifold and let $(E, \omega) \to V$ be a symplectic vector bundle. Fix a section $\tau: V \to LGr(E)$ and let $L_\tau \subset E$ be the associated bundle of Lagrangian subspaces. There is a contact structure $\hat{\xi}$ defined on a neighborhood of the zero section with the following properties:
\begin{itemize}
\item[(i)] the restriction of $\hat{\xi}$ to the tangent space of the zero section agrees with $\xi$, i.e.\ $(\hat{\xi} \cap TV)= \xi$
\item[(ii)] the conformal symplectic normal bundle of the zero section is precisely $E$
\item[(iii)] $L_\tau$ is coisotropic
\end{itemize}
\end{lemma}

\begin{proof}
This is essentially proved in \cite[Thm.\ 5.4]{Avdek}, but we briefly summarize the construction.   Fix a contact form $\xi= \ker \alpha$ on $V$. 
Fix a complex structure on $E$ compatible with $\omega$, thus making $E$ into a unitary vector bundle. Fix unitary bundle trivializations $\pi^{-1}(U_i) \simeq U_i \times \mathbb{R}^{2n}$ corresponding to a covering $\{U_i\}$ of $V$, and where we denote the fiber coordinates by $(q,p)$. Since $U(n)$ acts transitively on Lagrangian subspaces, we may assume that the trivializations $\pi^{-1}(U_i) \simeq U_i \times \mathbb{R}^n$ are chosen to that $\pi^{-1}(U_i) \cap L_\tau = U_i \times L_{\tau, i}$, where $L_{\tau, i} \subset \mathbb{R}^{2n}$ is a linear Lagrangian subspace (so independent of the $U_i$ coordinates). Let $\{h_i\}$ be a partition of unity subordinate to the covering $\{U_i\}$ and define $\lambda_i= \frac{h_i}{2}(\sum_j p^j dq_j - q^jdp_j)$. Finally, set $\lambda= \sum_i \lambda_i$ and set $\hat{\alpha}= \pi^* \alpha + \lambda$ and let $\hat{\xi}:= \ker \hat{\alpha}$. Now (i) is simply \cite[Thm.\ 5.4(1)]{Avdek}. For (ii), observe that given any $p \in V$, any purely radial tangent vector $v \in T_pE$ pairs trivially with any vector $w$ tangent to the zero section (because $\pi^*d \alpha(v, -)=0$ and $d\lambda_i(w,-)=0$). Hence (ii) is simply \cite[Thm.\ 5.4(2)]{Avdek}. 

It remains to verify (iii). To begin with, note that $\hat{\alpha}$ is nonzero on $L_\tau$ near $V \subset L_\tau$ (since $\hat{\alpha}|_{TV} =\alpha \neq 0$). Hence $\xi \cap TL_\tau$ is a codimension-$1$ distribution in $TL_\tau$.  It is enough to check, for arbitrary $x \in V$, that any element in the tangent space of the fiber $(L_\tau)_x \subset E_x$ is contained in $TL_\tau \cap \xi$ and pairs to zero with any other vector in $TL_\tau \cap \xi$ with respect to $d\hat{\alpha}$ (indeed, for dimension reasons, the tangent space of the fibers $(L_\tau)_x$ must then be the whole orthogonal complement of $(TL_\tau \cap \xi, d\hat{\alpha})$).  

Suppose $v$ is tangent to $(L_\tau)_x$; then $\hat{\alpha}(v)= \pi^*\alpha(v)+ \lambda(v)=\lambda(v)= \sum_i \lambda_i(v)$. Similarly $d\hat{\alpha}(v, -) = d \pi^* \alpha(v, -) + d\lambda(v,-)= d\lambda(v,-)$, and $d\lambda= \sum d\lambda_i$. So it is enough to check for each $i$ that $\lambda_i(v)=0$ and $d\lambda_i(v, -)=0$. We may check this in the local chart $U_i \times \mathbb{R}^{2n} \simeq \pi^{-1}(U_i)$. 

By definition $\lambda_i=\frac{h_i}{2} \lambda_0$ on $U_i \times \mathbb{R}^{2n}$, where $\lambda_0 =(\sum_j p^j dq_j - q^jdp_j)$. So $d\lambda_i(v,-) = (dh_i \wedge \lambda_0)(v,-) + h_i d\lambda_0(v,-)$. Since $dh_i(v)=0$, it is enough to check that $\lambda_0(v)=0$ and $d\lambda_0(v,-)=0$. Both of these statements are clear, since $\pi^{-1}(U_i) \cap L_{\tau,i} = U_i \times L_{\tau,i}$, where $L_{\tau,i} \subset U_i$ is a (linear) Lagrangian subspace. 
    (To check that $\lambda_0(v)=0$, note that $\lambda_0 = i_Z \omega_0$ where $Z= \frac{1}{2}(p\partial_p + q\partial_q)$ is a radial vector field. Now $Z$ is tangent to any linear space, and since $L_i$ is Lagrangian, we conclude $\lambda_0(v)= i_Z\omega_0(v)= \omega(Z, v)= 0$. For the second statement, $d\lambda_0(v,-)= \omega_0(v,-)$, which vanishes identically on vectors tangent to $L_{i,\tau}$ and to $U_i$.)
\end{proof}

\begin{lemma}\label{lemma:adapted}
Let $(V, \xi)$ be a co-orientable contact manifold and let $(E, \omega) \to V$ be a symplectic vector bundle as in \Cref{lemma:local-contact-form} equipped with an auxiliary complex structure. Fix $f: V \to U(n)$ and define $\phi: E \to E$, $(v; (q, p)) \mapsto (v; f_v(q, p))$.  Then $\phi^* \hat{\alpha}$ is a contact form near the zero section. Moreover $\phi^* \hat{\alpha}$ and $d\phi^* \hat{\alpha}$ both agree with $\hat{\alpha}$ and $d \hat{\alpha}$ on $TE|_V$.
\end{lemma}
\begin{proof} By definition $\hat{\alpha} = \pi^* \alpha+ \sum_i \lambda_i$. Evidently $\phi^*\pi^*\alpha= (\pi \circ \phi)^*\alpha= \pi^* \alpha$, so we just need to check, for any $x \in V$, that $\phi^* \lambda_i$ and $d\phi^* \lambda_i$ both agree with $\lambda_i$ and $d \lambda_i$ on $TE|_x$ (in particular, this would imply that $\phi^*\hat{\alpha}$ is a contact form at $x$). 

Let $\{U_i\}$ be as in the proof of \Cref{lemma:local-contact-form}. If $x \notin U_i$, there is nothing to prove. If $x \in U_i$, we consider the trivialization $\pi^{-1}(U)= U_i \times \mathbb{R}^{2n}$. Then  $\phi(x; q,p)= (x, A_x(q,p))$, where $A_{(-)}: U_i \to U(n)$ is induced by $f$ and the trivialization.  So we compute $\phi^*\lambda_i= \phi^*(h_i(pdq-qdp))= h_i \phi^*(pdq-qdp)$. Now $\phi^*(pdq-qdp) = A_x(p) dq + pd(A_x(q)) - (A_x(q) dp + q d(A_x p)) = A_x( p dq - qdp) + pq dA_x- qpdA_x= (pdq - qdp) + (pq dA_x - qp dA_x)$.  Since $(pq dA_x - qp dA_x)$ vanishes to second order the zero section, the claim follows.
\end{proof}

In the sequel, we use the following notation introduced by Gromov: given a space $X$ and a subspace $Z \subset X$, $Op(Z)$ denotes some unspecified open neighborhood of $Z$ in $X$ which may vary from line to line.

\begin{lemma}\label{lemma:contact-coisotropic-neighborhood}
Let $(W, \xi)$ be a contact manifold and let $(V, TV \cap \xi) \subset W$ be a contact submanifold. Let $S$ be a manifold and let $0 \in S$ be a point. Let $(C_s)_{s \in S}$ be a family of coisotropic submanifolds containing $V$, with the property that (a) the characteristic foliation is transverse to $V$, for all $s \in S$; (b) $(TC_s)|_V$ is constant (in $s$) and defines a constant Lagrangian distribution.

Then there is a family of contact isotopies $(\phi_s)_{s \in S}$, $\phi_s: Op(C_0) \to Op(C_s)$ with the following properties:
\begin{itemize}
\item[(i)] $\phi_s(C_0)= C_s$ and preserves the characteristic foliation
\item[(ii)] $\phi_s$ restricts to the identity map on $V$ and $(d\phi_s)$ is the identity on $T (Op(C_0))|_V \subset TE|_V$
\item[(iii)] $\phi_0: Op(C_0) \to Op(C_0)$ is the identity
\end{itemize}
\end{lemma}

\begin{proof}
Fix an auxiliary metric on $C$. Using the exponential map of the restriction of this metric to each leaf, we can define a family of diffeomorphisms $C_0 \to C_s$ which are constant on $V$, and such that the tangent map $df_s: (TC_0)|_V \to (TC_s)|_V$ is the identity. 

We now appeal to \cite[Cor.\ 2.26]{sackel2019getting} and the discussion directly following. We must argue that we can extend $f_s$ to a bundle map $TM|_{C_0} \to TM_{C_s}$ which preserves $\xi$ and the conformal symplectic structure.  We can then appeal to \cite[Thm.\ 2.5]{sackel2019getting} to extend to a contactomorphism, as desired (while it is not stated that way, \cite[Thm.\ 2.5]{sackel2019getting} works equally well in families).

To this end, the tangent space of $E$ along any co-isotropic $C$ decomposes as $$TM|_C= (TC \cap \xi)^{d\alpha} \oplus (TC \cap \xi) / (TC \cap \xi)^{d\alpha}  \oplus \xi|_C /(TC \cap \xi) \oplus TC/  (TC \cap \xi).$$

It follows from \Cref{lemma:char-fol-determined} that $f_s$ preserves $\xi$ and also preserves the conformal symplectic structure on $(TC \cap \xi) / (TC \cap \xi)^{d\alpha}$. It remains to define a map $F_s: \xi|_{C_0} \to \xi|_{C_s}$ preserving the conformal symplectic structure, and extending 
\begin{equation}\label{equation:extension-problem}
    df_s: (TC_0 \cap \xi)^{d\alpha} \oplus (TC_0 \cap \xi) / (TC_0 \cap \xi)^{d\alpha} \to (TC_s \cap \xi)^{d\alpha} \oplus (TC_s \cap \xi) / (TC_s \cap \xi)^{d\alpha}.
\end{equation}

Consider the bundle $Q \to S \times C_0$ whose fiber over $(s, x)$ is the space of maps $F_s(x): \xi_x \to \xi_{f_s(x)}$ extending $df_s(x)$. The fibers are (non-canonically) isomorphic to the space of symplectomorphisms of $(\mathbb{R}^{2n}, \omega_0)$ fixing the co-isotropic $\{p_{n-k}= \dots= p_n=0 \}$ where $2k$ is the codimension of the $C_s$.  This is a consequence of the linear algebra fact (Witt's theorem), that given $(V_1, \omega_1), (V_2, \omega_2)$ two symplectic vector spaces, given subspaces $F_i \subset V_i$ and a map $\phi: F_1 \to F_2$ such that $\phi^* (\omega_2)|_{F_2}= (\omega_2)|_{F_1}$, then we can extend $\phi$ to a linear symplectomorphism $V_1 \to V_2$. On $S \times V$, the identity map defines a section of $Q$, because $TC_s|_V$ is independent of $s \in S$. Similarly, the identity map defines a section of $Q$ over $\{0\} \times C_0$, because $f_0=id$. Now just extend the section to a neighborhood of $S \times V \cup \{0\} \times C_0 \subset S \times C_0$. 
\end{proof}

The following lemma, which is used in the proof of \Cref{lemma:contact-coisotropic-neighborhood}, is an elaboration of the discussion in \cite[Sec.\ 2.4]{sackel2019getting}.
\begin{lemma}\label{lemma:char-fol-determined}
Let $\xi= \ker \alpha, \xi'= \ker \alpha'$ be a co-orientable contact structures on a manifold $E$ of dimension $n+1+2k$. Let $C$ be a submanifold of dimension $n+1+k$ endowed with a nonsingular nonsingular foliation $\cF$ transverse to a submanifold $V^{n+1} \subset C$. Suppose that $C$ is co-isotropic with respect to both $\xi, \xi'$ and $\cF$ is the characteristic foliation. Suppose moreover that:
\begin{itemize}
\item $(\xi \cap TC)|_V= (\xi' \cap TC)|_V$ 
\item the restriction to $V$ of the conformal cosymplectic normal bundles $(TC \cap \xi) / (TC \cap \xi)^{d\alpha}$ and $(TC \cap \xi') / (TC \cap \xi')^{d\alpha'}$ coincide
\end{itemize} 
Then, after possibly replacing $C$ with a open subset containing $V$, we have:
\begin{itemize}
\item[(i)] $\xi \cap TC= \xi' \cap TC$
\item[(ii)] the conformal cosymplectic normal bundles $(TC \cap \xi) / (TC \cap \xi)^{d\alpha}$ and $(TC \cap \xi') / (TC \cap \xi')^{d\alpha'}$ are naturally isomorphic (i.e.\ the identity is an isomorphism of symplectic vector bundles)
\end{itemize}
\end{lemma}
\begin{proof}
Choose a point $p \in C$ and choose a vector field $V$ tangent to the leaves of the foliation whose time $1$ flow takes some point $q \in V$ in to $p$. 
Following \cite[Sec.\ 2.4]{sackel2019getting}, we have that $\mathcal{L}_V j^*\alpha = i_V j^*d\alpha$ vanishes on $\xi \cap TC$. Hence $\mathcal{L}_V j^*\alpha= \nu j^*\alpha$ for some $\nu: C \to \mathbb{R}$. Let $\Phi_{(-)}(-)$ be the flow of $V$. We compute $$\frac{d}{dt}(\Phi_t^* j^*\alpha) = \Phi_t^*(\mathcal{L}_V j^*\alpha) = \Phi_t^*( \nu j^*\alpha)= \nu(\Phi_t) \Phi_t^*(j^*\alpha).$$ In particular, writing $F(t)=(\Phi_t^*(j^*\alpha))_q \in TC_q$ and $f(t) = \nu(\Phi_t(q))$, we have the differential equation $F'(t)= f(t) F(t)$. The unique solution is $F(t)= e^{\int_0^t f(s)} F(0)$, where $F(0)= j^*\alpha_q$.  The consequence is that $\Phi_t^* j^*\alpha$ is determined by $(j^*\alpha)_q$; hence $(j^*\alpha)_p$ agrees up to a scaling factor with $\Phi_{-1}^* (j^*\alpha)_q$. Since $(j^*\alpha)_q= (j^*\alpha')_q$, this proves (i).

The proof of (ii) is similar. We compute $\mathcal{L}_V j^*d\alpha =i_V d (\nu j^* \alpha)= i_V (d\nu \wedge j^* \alpha+  \nu j^*d\alpha)$. Since $j^*\alpha$ vanishes on $\xi \cap TC$, we have $\mathcal{L}_V j^*d\alpha = \nu j^*d\alpha$ (where the equality is an equality of forms on $TC \cap \xi$ and $\nu: C \to \mathbb{R}$).  Consider $$\frac{d}{dt} (\Phi_t^* d\alpha)=\Phi_t^*(\mathcal{L}_Vd\alpha)= \nu(\Phi_t(q)) \Phi_t^*(d\alpha),$$ which is a differential equation on $(TC \cap \xi)_q$. Writing $F(t) = (\Phi_t^* d\alpha)$, $f(t)=\nu(\Phi_t(q))$, we have $F'(t)= f(t) F(t), F(0):= (j^*d\alpha)_q$. The unique solution is $F(t)= e^{\int_0^t f(s)} F(0)$. 
In particular, the conformal type of $F(t)$ depends only on $F(0)$ (in particular, it is independent of $\nu$). 
\end{proof}

\begin{corollary}\label{collary:thickening-extend}
    Any family of thickenings $(\sigma_s)_{s \in [0,1]}$ is induced by a global contact isotopy $(\Phi_s)_{s \in [0,1]}$.
\end{corollary}
\begin{proof}
Combine \Cref{lemma:contact-coisotropic-neighborhood} and \Cref{corollary:extend-isotopy}.
\end{proof}

We also have:
\begin{corollary}\label{corollary:fill-in-family-coisotropic}
In the notation of \Cref{lemma:contact-coisotropic-neighborhood}, suppose that $S=S^n$. The family of diffeotopies $S^n \times Op(C_0) \to E$ extends to a family $D^{n+1} \times Op(V) \to E$. 
\end{corollary}
\begin{proof}
This follows from the parametric neighborhood theorem \cite[Thm.\ 2.10]{sackel2019getting}. More precisely, we have a family of contact embeddings of $Op(C_0)$ whose derivatives on the contact submanifold $V \subset op(C_0)$ are fixed, so \Cref{lemma:char-fol-determined} along with loc.\ cit.\ implies that this extends. 
\end{proof}
\begin{corollary}\label{corollary:contractible-fibers}
Let $(W, \xi)$ be a contact manifold and let $(V, TV \cap \xi) \subset W$ be a contact submanifold. Given $\tau \in \mathfrak{Lag}(V \subset W)$, let $\mathfrak{Thick}(V \subset W)_\tau$ be the fiber of the map $\mathfrak{F}$ in \eqref{equation:forget-thick}. Then $\mathfrak{Thick}(V \subset W)_\tau$ is (weakly-)contractible. 
\end{corollary}
\begin{proof}
Combine \Cref{lemma:contact-coisotropic-neighborhood} and \Cref{corollary:fill-in-family-coisotropic}. 
\end{proof}

\begin{lemma}\label{lemma:triv-family-moser}
Let $(W, \xi_s)$ be a family contact manifolds, $s \in S$. Let $V \subset W$ be a submanifold. Suppose that $\alpha_s$ is a family of contact forms for $\xi_s$ with the property that $\alpha_s, d\alpha_s$ are independent of $s$ on $TW|_V$, and $(V, TV \cap \xi_s) \subset W$ is a contact submanifold. Then the family is trivial. More precisely, there is a map $\psi: Op(V \times S) \hookrightarrow W$ such that $\psi(-, s): (Op(V), \xi_s) \to (W, \xi_0)$ is a contact embedding and $d \psi(-, s)$ is the identity map on $TW|_V$. 
\end{lemma}

\begin{proof} 
If $S = pt$, then this follows immediately from the contact neighborhood theorem \cite[Thm.\ 2.5]{sackel2019getting}. The lemma follows by running the same argument with an $S$ parameter. Here is a brief sketch. Let $\beta_{s,t}:= t\alpha_s + (1-t) \alpha$. We want to build a family of diffeomorphisms $\phi_{s,t}: Op(V) \to Op(V)$ so that $\phi_{s,t}^*\beta_{s,t}= \lambda_{s,t} \beta_{s,0}$, where $\lambda_{s,t}$ is a positive function, such that $\phi_{s,t}|_V=id$ and $d\phi_{s,t}|_{TW|_V}$ is also the identity. By the usual Moser-type argument, we define $\phi_{s,t}$ as the time-$1$ flow of a vector field $V_{s,t}$. The construction on $V_{s,t}$ in the proof of  \cite[Thm.\ 2.5]{sackel2019getting} is canonical, so varies smoothly with $s$. Note that $V_{s,t}$ which vanishes identically on $V$, so the time one flow of this vector field is well-defined on a sufficiently small neighborhood of $V$.
\end{proof}

\begin{lemma}\label{lemma:fibration-hur}
The forgetful map $\mathfrak{F}$ in \eqref{equation:forget-thick} is a fibration. 
\end{lemma}
\begin{proof}
We may as well assume that $W=E$ is the total space of a unitary vector bundle over $V$, with a contact form $\hat{\alpha}$ defined in \Cref{lemma:local-contact-form}.  It is enough to prove that, for every point on the base, there is an open ball $\cO$ so that the preimage of $\cO$ is a fibration \cite[p.\ 51]{may1999concise} (indeed, the space of Lagrangian sections $\mathfrak{L}(V)$ is second countable, so any open cover admits a countable subcover, allowing us to apply loc.\ cit.\ ).  Note that $LGr(2k)= U(k)/O(k)$; choose a local slice $\Sigma$ near $e \in U(n)$.

Given $\tau \in \mathfrak{Lag}(V \subset E)$, let $$\mathcal{O}_\tau:= \{\sigma \in \mathfrak{Lag}(V \subset E) \mid \sigma= g \cdot \tau, g: V \to \Sigma \}.$$ 

We consider a diagram 
\begin{equation}
    \begin{tikzcd}
        Y \times \{0\} \ar[d] \ar[r] & \mathfrak{F}^{-1}(\mathcal{O}_\tau) \ar[d] \\
        Y \times [0,1] \ar[r, "f"] & \mathcal{O}_\tau
    \end{tikzcd}
\end{equation}

We map $(y, s)$ to the element $\tilde{f}(x, s) \in \tilde{U} \subset Maps(V \to \Sigma)$ which takes $f(y, 0)$ to $f(y,s)$ (i.e.\ $\tilde{f}(x, s)\cdot f(y,0)= f(y,s)$). Let $\psi_{(y,s)}: E \to E$ be induced by $\tilde{f}(y,s)$ (i.e.\ $\psi_{(y, s)}:= \tilde{f}(y,s) \cdot - $. Note that $\psi_{(y,0)}$ is the identity.

By \Cref{lemma:adapted}, $\psi_{(y,s)}^* \hat{\alpha}$ is still a contact form and its restriction to $V$, as well as the restrition of $d \psi_{(y,s)}^* \hat{\alpha}$ to $V$, is independent of $s$. So we have a map $Y \times I \to \mathcal{A}$, where $\mathcal{A}$ (temporarily) denotes the set of contact forms $\hat{\beta}$ on $Op(V) \subset E$ such that $\hat{\alpha}|_{TE|_V} = \hat{\beta}|_{TE|_V}$ and $d\hat{\alpha}|_{TE|_V} = d\hat{\beta}|_{TE|_V}$.  By \Cref{lemma:triv-family-moser} (setting $S:= Y\times I$ in the notation of that lemma), we have a map $Y \times I \to \mathfrak{F}^{-1}(\mathcal{O}_\tau)$ lifting the map $f: Y \times I \to \mathcal{O}_\tau$. 
\end{proof}

\begin{corollary}\label{corollary:forget-homotopy-equivalence}
\eqref{equation:forget-thick} is a (weak-)homotopy equivalence. 
\end{corollary}
\begin{proof}
By \Cref{lemma:fibration-hur} this is a fibration, and \Cref{corollary:contractible-fibers} shows that the fibers are (weakly-) contractible.
\end{proof}

\section{Equivariant sheaves} \label{apn: equivariance}

We recall the notion of equivariant sheaves in this section and will follow the conventions of \Cref{section:review-sheaf-theory}. In particular, given a topological space $X$, $sh(X)$ denotes the category of sheaves valued in some ambient stable category $\mathcal{C}$.

A continuous action on a suitable space $u: G \times X \rightarrow X$ induces an action of $sh(G)$ on $sh(X)$ through $!$-pushforward, $u_!: sh(G \times X) \rightarrow sh(X)$. We have $sh(G \times X) = sh(G) \otimes sh(X)$, where $\otimes$ is the tensor product on presentable categories, which naturally restricts to the stable ones. We define the co-invariants to be the geometric realization
$$ sh(X)_u \coloneqq \colim \limits \left( \cdots
\mathrel{\substack{\textstyle\rightarrow\\[-0.6ex]
\textstyle\rightarrow \\[-0.6ex]
\textstyle\rightarrow \\[-0.6ex]
\textstyle\rightarrow}}
sh(G \times G \times X) \mathrel{\substack{\textstyle\rightarrow\\[-0.6ex]
\textstyle\rightarrow \\[-0.6ex]
\textstyle\rightarrow}} sh(G \times X) \rightrightarrows
 sh(X) \right).
 $$
We remark that, as usual, we are omitting the degeneracy maps. We will denote this category by $sh(X)_G$ when there is no ambiguity about the action. Our first observation is that, when restricting to manifolds and Lie groups, invariants and co-invariants define the same notion in this setting and so there is an unambiguous notion of the category of $G$-equivariant sheaves:

\begin{lemma} \label{lem: inv=coinv-sheaf}
Let $M$ be a manifold and $u: G \times M \rightarrow M$ a Lie group action. Then,
there is a canonical equivalence $sh(M)^u = sh(M)_u$ between invariants and co-invariants.
\end{lemma}

\begin{proof}
Colimits in $\PrLst$ are computed through the equivalence $\PrLst = \PrRst^{op}$ by passing to right adjoints. In our case, the right adjoint of $u_!$ is $u^!$. However, since Lie groups are parallelizable, $u^! = u^* [\dim G]$. That is, $sh(X)_u$  is computed by the totalization

$$   \lim \limits \left( sh(M) \rightrightarrows sh(G \times M) 
\mathrel{\substack{\textstyle\rightarrow\\[-0.6ex]
\textstyle\rightarrow \\[-0.6ex]
\textstyle\rightarrow}} sh(G \times G \times M) 
\mathrel{\substack{\textstyle\rightarrow\\[-0.6ex]
\textstyle\rightarrow \\[-0.6ex]
\textstyle\rightarrow \\[-0.6ex]
\textstyle\rightarrow}} \cdots
 \right),
$$
where the first two arrows are $u^* [\dim G]$ and $p_2^* [\dim G]$ and similarly for the higher degrees. But since limits in $\PrLst$ and $\PrLst$ can both be computed in simply $\Cat$ and the fact that $u^*$ and $p^*$ are both left adjoints, this same diagram can be viewed as a limit diagram in $\PrLst$ as well. Up to a shifting $[\dim G]$ at each degree $\geq 2$ this is the same diagram which computes $sh(M)_u$.
\end{proof}

\begin{remark} \label{rmk: equivraint-sheaf-description}
An object in $sh(X)^G$ is a pair $(F, \phi)$, where $F \in sh(X)$ is a sheaf and $\phi: u^* F \xrightarrow{\sim} p_2^* F$ is an isomorphism. Indeed, the two arrows $sh(X) \rightrightarrows sh(G \times X)$, $u^*$ and $p_2^*$ give the map $\phi$, which, when restricting along the inclusion $i_g: X \xrightarrow{\sim} \{g\} \times X \hookrightarrow G \times X$, gives an identification $\phi_g: u_g^* F \xrightarrow{\sim} F$. Let $e \in G$ be the unit. Since $i_e$ is the degeneracy map from $X$ to $G \times X$, we have $\phi_e = \id_F$. Similarly, the arrows between $sh(G \times X)$ and $sh(G \times G \times X)$
impose the condition that the identification $\phi$ is associative, and respects the unit map.
\end{remark}

\begin{remark} \label{rmk: loc-equivraint-sheaf}
Following \Cref{lem: inv=coinv-sheaf} and \Cref{rmk: equivraint-sheaf-description}, we note that one can similarly consider an action by $\Loc(G)$, the category of local systems on $G$. However, in this setting, the notion of invariants and co-invariants differ from each other, and they are a priori different from the corresponding notion for a $sh(G)$-action. 
That is, one can consider the composition 
$$\Loc(G) \otimes sh(M) = sh_{0_G \times T^* M}(G \times M) \hookrightarrow sh(G \times M) \xrightarrow{u_!} sh(M).$$
Since $sh_{0_G \times T^* M}(G \times M) \hookrightarrow sh(G \times M)$ preserves both limits and colimits, it admits a left adjoint $\iota^*$ and a right adjoint $\iota^!$. 
Thus to compute the co-invariants $sh(M)_{u,loc}$, one passes to right adjoints as in \Cref{rmk: equivraint-sheaf-description}. For ${p_2}_!$ restricted on $\Loc(G) \otimes sh(M)$, the right adjoint is $p_2^*[\dim G]$ and $p_2^* F [\dim G] \in sh_{0_G \times T^* M}(G \times M)$ for any $F \in sh(M)$. However, $u^*F$ in general does not have the correct microsupport condition so the right adjoint of $u_!$ is the composition $\iota^! u^*[\dim G]$. Thus, a local-system $G$-equivariant sheaf is the data $(F, \phi, \cdots)$ where $F \in sh(M)$ is a sheaf and $\phi: \iota^! u^* F \xrightarrow{\sim} p_2^* F$ is a equivalence with higher coherence data. A similar discussion shows that objects in $sh(M)^{u,loc}$ admits the description $(F, \phi, \cdots)$ where $F \in sh(M)$ is a sheaf but the equivalence is now $\phi: \iota^* u^* F \xrightarrow{\sim} p_2^* F$, with $\iota^!$ replacing $\iota^*$.
\end{remark}

Ultimately, we will consider a setting where $M$ is a complex manifolds and where the natural objects to consider are (holomorphic) coherent sheaves and D-modules. Since coherent sheaves can be defined as topological sheaves which are modules over the ring of holomorphic functions, it is convenient to henceforth restrict our discussion to sheaves valued in the the $1$-category of abelian groups. We further assume that our coefficient category is the (abelian) category of $R$-modules for some discrete ring $R$ which is always understood to be the ring of complex numbers whenever we discuss coherent sheaves.  

\begin{definition} \label{def: top-equivaraint-sheaf}
Let $X$ be a topological space, $G$ a topological group, and $u: G \times X \rightarrow X$ a continuous group action. Then the category of $G$-equivariant abelian sheaves on $X$ is the  $(2,1)$-limit in $\Ab$ of the three-term truncated simplicial diagram  
$$ sh(X; R)^u \coloneqq \lim \limits \left( sh(X; R) \rightrightarrows sh(G \times X; R) 
\mathrel{\substack{\textstyle\rightarrow\\[-0.6ex]
\textstyle\rightarrow \\[-0.6ex]
\textstyle\rightarrow}} sh(G \times G \times X; R) \right),
$$
where the arrows are induced by the $*$-pullback from the corresponding diagram in topological spaces,
$$ 
G \times G \times X
\mathrel{\substack{\textstyle\rightarrow\\[-0.6ex]
\textstyle\rightarrow \\[-0.6ex]
\textstyle\rightarrow}}
G \times X  \rightrightarrows X.$$
Again, we are omitting the degeneracy maps. We will denote this category by $sh(X; R)^G$ when there is no ambiguity about the action.
\end{definition}

We see from \Cref{rmk: equivraint-sheaf-description} that there is a symmetric monoidal structure on $sh(X; R)^G$ given by
\begin{equation} \label{for: equivariant-tensor}
 (F,\phi) \otimes (G, \psi) \coloneqq (F \otimes G, \phi \otimes \psi).
 \end{equation}
The unit is the constant sheaf $1_X$ with the canonical equivariant structure $\phi_{G,X}$ by the composition 
$$u^* 1_X = 1_{G \times X} = p_2^* 1_X.$$ 

In the next section, we compare equivariant sheaves on a covering space.
\begin{lemma} \label{lem: group-action-geometry}
Let $u: G \times X \rightarrow X$ be a free and proper action in topological spaces where Galois theory of covering spaces holds, e.g., when both $G$ and $X$ are manifolds. Assume that the deck transformation $N$ of the covering map $p: X \rightarrow \overline{X}$ acts through the $G$-action, i.e., there exists an inclusion $N \hookrightarrow G$ as a normal subgroup so that the homeomorphism $u_n: X \xrightarrow{\sim} X$, $n \in N$, is given by the same action by viewing $N \trianglelefteq \pi_1(X,x)$, for some fixed $x \in X$. Then, there exists an $q: G \rightarrow H \coloneqq G/N$ action $\overline{u}$ on $\overline{X}$ and, for any $g \in G$, the following diagram commute:

$$
\begin{tikzpicture}
\node at (0,2) {$X$};
\node at (4,2) {$X$};
\node at (0,0) {$\overline{X}$};
\node at (4,0) {$\overline{X}$};

\draw [->, thick] (0.3,2) -- (3.7,2) node [midway, above] {$u_g$};
\draw [->, thick] (0.3,0) -- (3.7,0) node [midway, above] {$\overline{u}_{\bar{g}}$};

\draw [->, thick] (0,1.7) -- (0,0.3) node [midway, right] {$p$}; 
\draw [->, thick] (4,1.7) -- (4,0.3) node [midway, right] {$p$};
\end{tikzpicture}
$$ 
\end{lemma}

\begin{proof}
The question is to descend the composition $G \times X \xrightarrow{u} X \xrightarrow{p} \overline{X}$, i.e., we argue that the dotted line $\overline{u}$ in the commuting diagram exists:
$$
\begin{tikzpicture}
\node at (0,2) {$G \times X$};
\node at (5,2) {$X$};
\node at (0,0) {$H \times \overline{X}$};
\node at (5,0) {$\overline{X}$};

\draw [->, thick] (0.7,2) -- (4.7,2) node [midway, above] {$u$};
\draw [->, thick, dashed] (0.7,0) -- (4.7,0) node [midway, above] {$\overline{u} $};

\draw [->, thick] (0,1.7) -- (0,0.3) node [midway, right] {$(q \times p)$}; 
\draw [->, thick] (5,1.7) -- (5,0.3) node [midway, right] {$p$};

\draw [->, thick] (0.7,1.6) -- (4.6,0.3) node [midway, above] {$ $}; 
\node at (3.7,1.35) {$\circlearrowleft$};
\end{tikzpicture}
$$ 

Since the left vertical map is a quotient, it suffices to check that
$$p \circ u \left( (n_1, n_2) \cdot (g, x) \right) = p \circ u (g, x)$$
for all $n_1, n_2 \in N$, $g \in G$, and $x \in X$. We thus compute that 
\begin{align*}
p \circ u \left( (n_1, n_2) \cdot (g, x) \right) 
&= p \circ u (n_1 g, u_{n_2}(x) ) = p \circ u_{n_1 g n_2}(x) = p \circ u_{n_1} \circ u_{g n_2 g^{-1}} \circ u_g(x) \\
&= p \circ u_g (x) = p \circ u (g,x).
\end{align*}

\end{proof}

\begin{proposition} \label{prop: group-action-category}
In the situation of \Cref{lem: group-action-geometry}, we have equivalences of symmetric monoidal categories
\begin{enumerate}
\item $sh(X)^N = sh(\overline{X})$, and 
\item $sh(X)^G = sh(\overline{X})^H$.
\end{enumerate}
\end{proposition}

\begin{proof}
(1) is a standard result regarding covering spaces, and the main observation is that there is a pullback diagram
$$
\begin{tikzpicture}
\node at (0,2) {$N \times X$};
\node at (4,2) {$X$};
\node at (0,0) {$X$};
\node at (4,0) {$\overline{X}$};

\node at (0.6,1.5) {$\ulcorner$};

\draw [->, thick] (1,2) -- (3.7,2) node [midway, above] {$u_o$};
\draw [->, thick] (0.4,0) -- (3.7,0) node [midway, above] {$p$};

\draw [->, thick] (0,1.7) -- (0,0.3) node [midway, right] {$p_2$}; 
\draw [->, thick] (4,1.7) -- (4,0.3) node [midway, right] {$p$};
\end{tikzpicture},
$$
where $u_o$ is the restriction of the action $u$ to $N$. Thus, for any $\bar{F} \in sh(X)$, since $p \circ u = p \circ p_2$, we have a canonical equivariant structure
$$ u_o^* (p^* \bar{F}) = (p \circ u_o)^* \bar{F} = (p \circ p_2)^* \bar{F} = p_2^* (p^* \bar{F}).$$ 
In fact, the same statement for the stable category setting is also true, and its proof is a standard application of the Lurie-Barr-Beck monadicity theorem and Beck-Chevalley. The statement in the derived algebraic geometry setting can be found in \cite[Proposition 7.2.2]{gaitsgory-rozenblyum-DAG-I} and a detailed exposition can be found in \cite[Proof of Prop.\ 2.3.1]{hchen}. We note that this general version recovers our case since $N$ is discrete and we work over the abelian categories. 

For (2), as remarked, the equivalence for (1) is induced by
$p^*: sh(\overline{X}) \rightarrow sh(X)$. An object in $sh(\overline{H})^H$ is a sheaf $\bar{F} \in \sh(\overline{X})$ with an $H$-equivariant structure $\phi: \overline{u}^* \bar{F} = \overline{p}_2^* \bar{F}$ where we use $\overline{p}_2: H \times \overline{X} \rightarrow \overline{X}$ to denote the projection to the second component. For a such an $\bar{F}$, we have 
$$ u^* p^* \bar{F} =  (q \times p)^* \overline{u}^* \bar{F} \overset{\phi}{=}  (q \times p)^* \overline{p}_2^* \bar{F} = p_2^* \overline{p}^* F $$
and thus a $G$-equivariant structure of $\overline{p}^*F$.
To construct the inverse, we note that $\phi: u^* F = p_2^* F$ in particular implies that $u|_{K \times X}^* F = p_2|_{K \times X}^* F$, and so  $F = p^* \bar{F}$ for some $\bar{F} \in sh(\overline{X})$ by (1). Furthermore, $\phi$ equip canonically an equivalence
$$(q \times p)^* \overline{u}^* \bar{F} = (q \times p)^* \overline{p}_2^* \bar{F},$$
which restricts to the identity $\id_F$ since $(q \times p) |_{K \times X} = (1 \times p)$, where we abuse the notation and use $1$ to denote the unit $1 \in G$ as well as the constant map induced by it. Thus, we conclude that $\phi$ is lifted from 
$\bar{\phi}: \overline{u}^* \bar{F} = \overline{p}_2^* \bar{F}$ and we obtain an $H$-equivariant structure of $\bar{F}$. 
\end{proof}

\begin{remark} \label{rmk: equivariant-subs}
Consider the case when $X$ is a smooth manifold and $G$ a Lie group. The notion of equivariance and \Cref{prop: group-action-category} naturally restricts to the subcategory of real constructible sheaves $sh_{\R-c}(X)$ by \cite[Proposition 8.4.6]{kashiwara-schapira}. Similarly, when $X$ is a complex manifold and $G$ a complex Lie group. They restrict to the subcategory of complex real constructible sheaves with perfect stalks $sh_{\C-c}(X)^b$ by \cite[Proposition 8.5.7]{kashiwara-schapira}. Similarly, when there is a $G$-invariant subanalytic (resp. complex) Legendrian, we can further restrict to the category $sh_\Lambda(X)$ (resp. $sh_\Lambda(X)^b$).
\end{remark}

We now turn our attention to the holomorphic setting and assume $X$ is a complex manifold and $G$ a complex Lie group. Recall that in this setting, when $f: X \rightarrow Y$ is a holomorphic map, there is a natural ring homomorphism
\begin{align*}
f^{*,top} \cO_Y &\rightarrow \cO_X \\
h &\mapsto h \circ f
\end{align*}
which induces the $*$-pullback for analytic sheaves
\begin{align*}
f^*: \cO_X - mod &\rightarrow \cO_Y -mod \\
\cM &\mapsto \cO_X \otimes_{f^{*,top} \cO_Y} f^{*,top} \cM.
\end{align*}
Here we use the notation $f^{*,top}$ to emphasize that this is the $*$-functoriality for topological sheaves $sh(X)$. In particular, we have the following holomorphic version of \Cref{def: hol-equivaraint-sheaf}:

\begin{definition} \label{def: hol-equivaraint-sheaf}
Let $X$ be a complex manifold, $G$ a complex Lie group, and $u: G \times X \rightarrow X$ a holomorphic action. Then the category of holomorphic $G$-equivariant sheaves on $X$ is defined to be the  $(2,1)$-limit in $\Ab$ of the three-term truncated simplicial diagram  
$$ (\cO_X -mod)^u \coloneqq \lim \limits \left( \cO_X -mod \rightrightarrows \cO_{G \times X}-mod 
\mathrel{\substack{\textstyle\rightarrow\\[-0.6ex]
\textstyle\rightarrow \\[-0.6ex]
\textstyle\rightarrow}} \cO_{G \times G \times X} -mod \right),
$$
where the arrows are induced by the (holomorphic) $*$-pullback from the corresponding diagram in complex manifolds,
$$ 
G \times G \times X
\mathrel{\substack{\textstyle\rightarrow\\[-0.6ex]
\textstyle\rightarrow \\[-0.6ex]
\textstyle\rightarrow}}
G \times X  \rightrightarrows X.$$
\end{definition}

\begin{remark}
Equivalently, this means that the identification $\phi_g: u_g^* \cM \xrightarrow{\sim} \cM$ depends on $g$ holomorphically. 
\end{remark}

We further incorporate equivariance with module structures. 
Let $\sR \in sh(X;\Ring)$ be a ring-valued sheaves. One can check directly that this coincide with the notion of ring-objects in sheaves $sh(X)$, $\Ring\left(sh(X)\right)$. Similarly, the notion of $\cO_X$-algebras coincide with the notion of ring-objects in $\cO_X-mod$. Because the discussion of these two settings are parallel to each other, we will use the notation $\sS$ to mean either the functor $sh(-)$ or $\cO_{(-)}-mod$, and $*$-functoriality in the respective setting. The $*$-pullback are symmetric, and, as a result, the diagram restricts to ring-objects and we obtain a well-defined category of equivariant sheaves of rings
$$ \left( \Ring(\sS(X)\right)^G  \coloneqq \lim \limits \left( \Ring(\sS(X)) \rightrightarrows \Ring(\sS(G \times X))
\mathrel{\substack{\textstyle\rightarrow\\[-0.6ex]
\textstyle\rightarrow \\[-0.6ex]
\textstyle\rightarrow}} \Ring(\sS(G \times G \times X)) \right).
$$

However, because $\sS(X)^G$ has a symmetric monoidal structure (\ref{for: equivariant-tensor}), the notion of ring-objects in equivariant sheaves, $\Ring( \sS(X)^G)$, also makes sense. One can check easily that the two notions coincide.
\begin{lemma} \label{lem: top-equivaraint-ring-sheaf}
For a $G$-action on $X$, we have
$$ \Ring( \sS(X)^G) =  \left(\Ring(\sS(X))\right)^G.$$
\end{lemma}

\begin{proof}
A ring object $(\cR,\phi)$ in $\sS(X)^G$ is a equivariant sheaf with a multiplication and a unit 
$$ (\cR,\phi) \otimes (\cR,\phi) \rightarrow (\cR, \phi), \ (1_X, \phi_{G,X}) \rightarrow (\cR,\phi).$$
By writing out the definition, one sees that $\phi$ provides $G$-equivalence that equips $\cR \in \Ring(\sS(X))$ an equivariant structure.
\end{proof}

Fix $(\cR,\phi) \in \Ring(\sS(X)^G)$, one can consider the category of modules $(\cR,\phi) - \mod \left( \sS(X)^G \right)$. Similarly to the previous Lemma \ref{lem: top-equivaraint-ring-sheaf}, we have a version for modules.

\begin{lemma} \label{lem: top-equivaraint-module-sheaf}
The category $(\cR,\phi) - \mod \left( \sS(X)^G \right)$ can be computed as the limit
$$
\lim \  \left( \cR - \mod  \rightrightarrows (u^* \cR) - \mod
\mathrel{\substack{\textstyle\rightarrow\\[-0.6ex]
\textstyle\rightarrow \\[-0.6ex]
\textstyle\rightarrow}} \left((u \times \id)^* u^* \cR\right) - \mod \right).
$$
Here we implicitly use $\phi: u^* \cR \xrightarrow{\sim} p_2^* \cR$ to identify $u^* \cR- \mod = p_2^* \cR -\mod$ and similarly for the categories over $G \times G \times X$. 
\end{lemma}

\begin{remark}
Roughly speaking, \Cref{lem: top-equivaraint-ring-sheaf} and \Cref{lem: top-equivaraint-module-sheaf} say that a ring object $(\cR,\phi) \in \Ring(\sS(X)^G)$ is a sheaf of rings, such that identification $\phi_g: u_g^* \cR = \cR$, which is continuous in $g$, is compatible with the ring structure. That is, the two different morphisms
$$ (u_g^* \cR) \otimes (u_g^* \cR) = u_g^* (\cR \otimes \cR) \rightarrow u_g^* \cR \xrightarrow{\phi_g} \cR$$
and 
$$(u_g^* \cR) \otimes (u_g^* \cR) \xrightarrow{\phi_g \otimes \phi_g} \cR \otimes \cR \rightarrow \cR$$
are the same, and so for other parts of the structures, and similarly for the module structure.
\end{remark}

We now upgrade \Cref{prop: group-action-category} to the version for modules.
\begin{corollary} \label{cor: group-action-category-with-ring}
In the situation of \Cref{prop: group-action-category}, fix $(\sR,\phi)$ in $\Ring\left(\sS(\overline{X})\right)$, then
\begin{enumerate}
\item $p^*\sR$ admits a canonical $N$-equivariant structure $\phi_{p^* \sR}$ such that, 
$$ (p^*\sR, \phi_{p^* \sR})  - \mod  =   \sR  - \mod,$$ 
\item the morphism $(q \times p)^* \phi$ equips $p^* \sR$ with a $G$-equivariant structure, which we abuse notation in denoting by $p^* \phi$, such that 
$$\left(p^* \sR, p^* \phi \right)- \mod  =   (\sR,\phi) - \mod.$$
\end{enumerate}

\end{corollary}

\begin{proof}
We can directly apply \Cref{prop: group-action-category} once we expand the categories in the statement as limit diagrams in \Cref{lem: top-equivaraint-ring-sheaf} and \Cref{lem: top-equivaraint-module-sheaf}.
\end{proof}

\bibliographystyle{plain}
\bibliography{refs}

\end{document}